\documentclass[12pt,a4paper,reqno]{amsart}
\usepackage{amsmath}
\usepackage{amsthm}
\usepackage{amsfonts}
\usepackage{amssymb}
\usepackage{color}
\usepackage{mathrsfs}
\usepackage{graphicx}

\numberwithin{equation}{section}

\theoremstyle{plain}
\newtheorem{theorem}{Theorem}

\theoremstyle{plain}
\newtheorem{lemma}{Lemma}[section]

\theoremstyle{remark}
\newtheorem*{remark}{Remark}

\theoremstyle{remark}

\theoremstyle{definition}
\newtheorem*{no-red}{Restriction on Red Coloring}

\theoremstyle{definition}
\newtheorem*{step1}{Step 1}

\theoremstyle{definition}
\newtheorem*{step2}{Step 2}

\theoremstyle{definition}
\newtheorem*{step3}{Step 3}

\theoremstyle{definition}
\newtheorem*{step4}{Step 4}

\theoremstyle{definition}
\newtheorem*{step5}{Step 5}

\theoremstyle{definition}
\newtheorem*{step6}{Step 6}

\theoremstyle{definition}
\newtheorem*{step7}{Step 7}

\theoremstyle{definition}
\newtheorem*{step8}{Step 8}

\theoremstyle{definition}
\newtheorem*{step9}{Step 9}

\theoremstyle{definition}
\newtheorem*{part1}{Part 1}

\theoremstyle{definition}
\newtheorem*{part2}{Part 2}

\theoremstyle{definition}
\newtheorem*{part3}{Part 3}

\def\bfe{\mathbf{e}}

\def\bfm{\mathbf{m}}
\def\bfn{\mathbf{n}}

\def\bfs{\mathbf{s}}
\def\bft{\mathbf{t}}

\def\bfv{\mathbf{v}}
\def\bfw{\mathbf{w}}
\def\bfx{\mathbf{x}}
\def\bfy{\mathbf{y}}
\def\bfz{\mathbf{z}}

\def\bfD{\mathbf{D}}

\def\bfG{\mathbf{G}}

\def\bfS{\mathbf{S}}

\def\bfdelta{{\boldsymbol\delta}}

\def\bzero{\mathbf{0}}

\def\T{\mathbf{T}}

\def\dd{\mathrm{d}}
\def\ee{\mathrm{e}}
\def\ii{\mathrm{i}}

\def\eps{\varepsilon}

\def\Rr{\mathbb{R}}

\def\Zz{\mathbb{Z}}

\def\AAA{\mathcal{A}}
\def\BBB{\mathcal{B}}
\def\CCC{\mathcal{C}}

\def\FFF{\mathcal{F}}

\def\III{\mathcal{I}}
\def\JJJ{\mathcal{J}}

\def\LLL{\mathcal{L}}
\def\MMM{\mathcal{M}}

\def\SSS{\mathcal{S}}

\def\ZZZ{\mathcal{Z}}

\def\frakB{\mathfrak{B}}

\def\frakG{\mathfrak{G}}

\def\frakI{\mathfrak{I}}

\def\GGGG{\mathscr{G}}

\def\RRRR{\mathscr{R}}

\DeclareMathOperator{\meas}{meas}

\renewcommand{\le}{\leqslant}
\renewcommand{\ge}{\geqslant}

\title{Uniformity in cube-covering systems}

\author[Beck]{J. Beck}

\address{Department of Mathematics, Hill Center for the Mathematical Sciences, Rutgers University, Piscataway NJ 08854, USA}

\email{jbeck@math.rutgers.edu}

\author[Chen]{W.W.L. Chen}

\address{School of Mathematical and Physical Sciences, Faculty of Science and Engineering, Macquarie University, Sydney NSW 2109, Australia}

\email{william.chen@mq.edu.au}

\author[Yang]{Y. Yang}

\address{School of Science, Beijing University of Posts and Telecommunications, Beijing 100876, China}

\email{yangyx@bupt.edu.cn}

\begin{document}



\begin{abstract}
We establish various analogs of the Kronecker--Weyl equidistribution theorem that can be considered higher-dimensional versions of results
established in our earlier investigation in \cite{BCY1} of the discrete $2$-circle problem studied in 1969 by Veech~\cite{veech69}.
Whereas the Veech problem can be viewed as one of geodesic flow on a $2$-dimensional flat surface,
here we study geodesic flow in higher-dimensional flat manifolds.
This is more challenging, as the overwhelming majority of the available proof techniques for non-integrable flat systems are based
on arguments in dimension~$2$.
For higher dimensions, we need a new approach.
\end{abstract}

\maketitle

\thispagestyle{empty}

%
%

\section{Analog of the Veech $2$-circle problem}\label{sec1}

We extend the idea of Veech~\cite{veech69} and our earlier work \cite{BCY1} to higher dimension.
First we $2$-color the unit torus $[0,1)^2$ red and green in such a way that each of the red and green parts is the union of finitely many polygons.
Figure~1.1 shows two examples, where the shaded part represents red and the white part represents green.
In particular, we assume that the green part has positive area.

\begin{displaymath}
\begin{array}{c}
\includegraphics[scale=0.8]{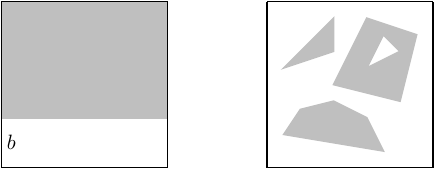}
\\
\mbox{Figure 1.1: examples of $2$-colorings of the torus $[0,1)^2$}
\end{array}
\end{displaymath}

Observe that in the picture on the right, one of the red (shaded) parts does not look like a polygon, but it is the union of finitely many polygons.
A similar remark applies to the green (white) part.

Next, we consider a $2$-torus system as shown in Figure~1.2, where each square represents the unit torus $[0,1)^2$, with identical $2$-coloring.

\begin{displaymath}
\begin{array}{c}
\includegraphics[scale=0.8]{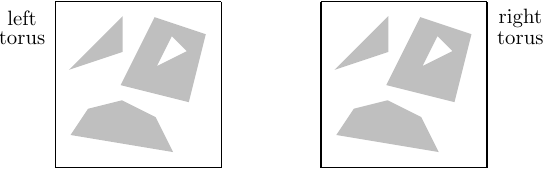}
\\
\mbox{Figure 1.2: a $2$-torus system with identical coloring}
\end{array}
\end{displaymath}

Let $\bfv=(1,\alpha_1,\alpha_2)$ be a Kronecker vector, and let $\bfv_0=(\alpha_1,\alpha_2)$.

Let $\bfs_0\in[0,1)^2$ be an arbitrary starting point, and consider the $\bfv_0$-shift sequence
\begin{displaymath}
\bfs_n=\bfs_0+n\bfv_0,
\quad
n=0,1,2,3,\ldots,
\end{displaymath}
in the unit torus $[0,1)^2$; in other words, \textit{modulo one}.
Assume that the point $\bfs_0$ is on the left torus.
If $\bfs_1$ is in the red (shaded) part, then we keep it on the left torus.
If $\bfs_1$ is in the green (white) part, then we move it to the corresponding point on the right torus.
In general, $\bfs_n$ is on a particular torus.
If $\bfs_{n+1}$ is in the red (shaded) part, then we keep it on the same torus.
If $\bfs_{n+1}$ is in the green (white) part, then we move it to the corresponding point on the other torus.
Thus the sequence $\bfs_0,\bfs_1,\bfs_2,\bfs_3,\ldots$ moves from one torus to the other whenever it hits the green part.
The problem is then to describe the distribution of this sequence in the union of the two tori,
clearly a \textit{parity} problem motivated by the Kronecker--Weyl equidistribution theorem.

We can visualize this discrete $2$-torus system on the plane as a simple continuous system in $3$-space.
Figure~1.3 illustrates this observation in the case of the simpler $2$-coloring in the picture on the left in Figure~1.1.

\begin{displaymath}
\begin{array}{c}
\includegraphics[scale=0.8]{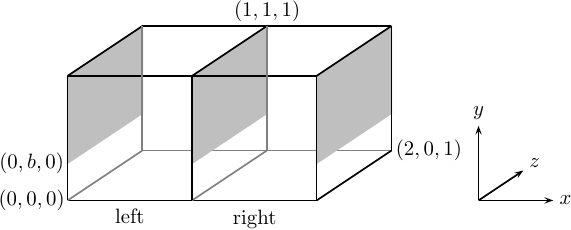}
\\
\mbox{Figure 1.3: $2$-cube-$b$ $3$-manifold with repeated barriers}
\end{array}
\end{displaymath}

Here there are three $yz$-parallel square faces of the $2$-cube solid, each of which is in part a barrier, colored red (shaded),
and in part a gate, colored green (white).
The latter is non-empty, and permits travel between the two cubes.
For these two cubes to form a $3$-manifold, we have to guarantee that it is boundary-free.
We use boundary identification which is a modification of the boundary identification for the torus $[0,1)^3$.
The two $xy$-parallel square faces with $z=0$ are identified with the two $xy$-parallel square faces with $z=1$ by trivial perpendicular translation.
The two $xz$-parallel square faces with $y=0$ are identified with the two $xz$-parallel square faces with $y=1$ by trivial perpendicular translation.
The right side of the red (shaded) rectangle on the $yz$-parallel square face with $x=0$
is identified with the left side of the red (shaded) rectangle on the $yz$-parallel square face with $x=1$,
while the right side of the red (shaded) rectangle on the $yz$-parallel square face with $x=1$
is identified with the left side of the red (shaded) rectangle on the $yz$-parallel square face with $x=2$.
Finally, the green (white) rectangle on the $yz$-parallel square face with $x=0$
is identified with the green (white) rectangle on the $yz$-parallel square face with $x=2$.
For convenience, we refer to this as the $2$-cube-$b$ $3$-manifold.

We thus have a flat $3$-manifold, with euclidean metric almost everywhere, and with boundary identification via perpendicular translation.
Thus geodesic flow in this $3$-manifold is $1$-direction linear flow.
It moves rather like $1$-direction geodesic flow on the torus $[0,1)^3$, and the novelty comes from the effect of the barriers.

There is clearly an equivalence between the discrete $2$-dimensional $2$-torus system and this new continuous $3$-dimensional $2$-cube system.
An infinite $\bfv_0$-shift sequence is equidistributed on the $2$-torus with the $2$-coloring given in the picture on the left of Figure~1.1
if and only if the corresponding half-infinite $1$-direction geodesic with direction vector $\bfv$ is equidistributed in the $2$-cube-$b$ $3$-manifold.

Assume now that $\bfv=(1,\alpha_1,\alpha_2)\in\Rr^3$ is a Kronecker vector.
Is it true that every half-infinite $1$-direction geodesic with direction vector $\bfv$ is equidistributed in any $2$-cube-$b$ $3$-manifold with $0<b<1$?

It turns out that for every Kronecker vector~$\bfv$, there are infinitely many values of the parameter $b$ for which equidistribution fails.
To explain this, we need to look at the corresponding problem in lower dimension.
The projection of the $2$-cube-$b$ $3$-manifold to the $xy$-plane gives rise to the $2$-square-$b$ surface
which arises from the work of Veech~\cite{veech69}.
Some of the anti-equidistribution results on such surfaces obtained recently by the authors in \cite{BCY1}
can be converted to anti-equidistribution results on $2$-cube-$b$ $3$-manifolds.

For instance, let $\alpha_1\in(0,1/2)$ be irrational, and let $b=2\alpha_1$.
Then for every $\alpha_2$ for which $\bfv=(1,\alpha_1,\alpha_2)\in\Rr^3$ is a Kronecker vector,
every half-infinite $1$-direction geodesic with direction vector $\bfv$ violates equidistribution in the $2$-cube-$b$ $3$-manifold.
For more details, see \cite[Theorem~2.1]{BCY1}.

The papers \cite{BCY1} and \cite{veech69} contain some equidistribution results on the $2$-square-$b$ surface.
These, unfortunately, do not immediately lead to corresponding results on the $2$-cube-$b$ $3$-manifold.
Nevertheless, using a different approach, we can establish equidistribution for most half-infinite geodesics in the $2$-cube-$b$ $3$-manifold.
Furthermore, we can generalize the result to any $2$-coloring of the unit torus $[0,1)^2$ where each of the red and green parts
is the union of finitely many polygons, and where the green part has positive area.
The richness of the possibilities to fix such $2$-colorings is particularly interesting.

Indeed, we can consider an $n$-torus system, with $n$ copies of the unit torus $[0,1)^2$, where $n\ge2$ is an integer.
This then leads to a flat $3$-manifold, with euclidean metric almost everywhere, and with boundary identification via perpendicular translation.
For instance, if we take $n=4$ and use the $2$-coloring of the torus $[0,1)^2$ as shown in the picture on the right in Figure~1.1,
then we have the $4$-cube $3$-manifold with repeated barriers as shown in Figure~1.4.

\begin{displaymath}
\begin{array}{c}
\includegraphics[scale=0.8]{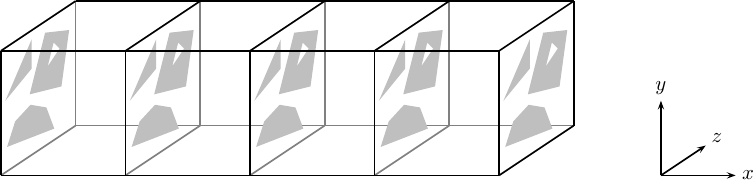}
\\
\mbox{Figure 1.4: a $4$-cube $3$-manifold with repeated barriers}
\end{array}
\end{displaymath}

\begin{theorem}\label{thm1}
Let $n\ge2$ be an integer, and let $\MMM$ be any $n$-cube $3$-manifold with barriers,
where the $yz$-parallel square faces have identical $2$-coloring
such that each of the red and green parts is the union of finitely many polygons,
and where the green part has positive area.
Then for almost every starting point and almost every direction $\bfv=(1,\alpha_1,\alpha_2)\in\Rr^3$,
the corresponding half-infinite $1$-direction geodesic is equidistributed in~$\MMM$.
\end{theorem}

As a trivial corollary, we deduce that the half-infinite $1$-direction geodesic spends asymptotically the same amount of time
in each one of the $n$ cubes of the $n$-cube $3$-manifold.

We remark that any polygon in the given $2$-coloring can be replaced by a circle, an ellipse, or any other \textit{piecewise smooth closed curve}.
It requires an extra analytic discussion in the proof that we postpone to Section~\ref{sec6}.

An immediate question that arises is whether we can extend Theorem~\ref{thm1} to include every Kronecker direction $\bfv\in\Rr^3$.

By a grid type $2$-coloring of the torus $[0,1)^2$, we mean dividing the torus $[0,1)^2$ into $m^2$ subsquares in the natural way,
where $m$ is a positive integer, and coloring at least one of the subsquares green and the remainder red, as illustrated in Figure~1.5.

\begin{displaymath}
\begin{array}{c}
\includegraphics[scale=0.8]{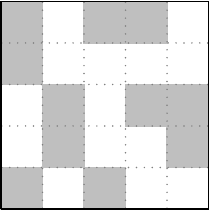}
\\
\mbox{Figure 1.5: grid type $2$-coloring of the torus $[0,1)^2$}
\end{array}
\end{displaymath}

Using a similar method, the authors can establish the following stronger conclusion in this setting; see \cite[Theorem~4]{BCY2}.

\begin{theorem}\label{thm2}
Let $n\ge2$ be an integer, and let $\MMM$ be any $n$-cube $3$-manifold with barriers,
where the $yz$-parallel square faces have identical grid type $2$-coloring,
and where the green part has positive area.
Then every half-infinite $1$-direction geodesic with a Kronecker direction $\bfv=(1,\alpha_1,\alpha_2)\in\Rr^3$ is equidistributed in~$\MMM$.
\end{theorem}

Figure~1.6 shows the $2$-cube box with a $2$-coloring on the middle $yz$-parallel square face
such that each of the red (shaded) and green (white) parts is the union of finitely many polygons,
and where the green part has positive area.
Consider billiard in this $2$-cube box, where in addition to the square faces on the surface of the box,
there are additional barriers on the middle $yz$-parallel square face given by the parts colored red.
As usual, we consider the ideal case of a point billiard that bounces back at any barrier, following the well-known rules of optical reflection.

\begin{displaymath}
\begin{array}{c}
\includegraphics[scale=0.8]{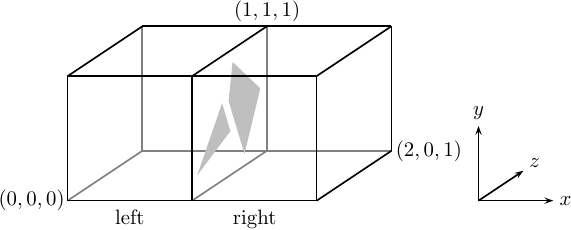}
\\
\mbox{Figure 1.6: $2$-cube box with barriers in the middle}
\end{array}
\end{displaymath}

We are interested in the long term bahavior of the billiard orbit.
We shall show that Theorem~\ref{thm1} contributes to our understanding of such questions.

To establish equidistribution for such billiard orbits, we extend the idea of K\"{o}nig and Sz\"{u}cs \cite{KS13} and apply $3$-dimensional \textit{unfolding}.
This converts the billiard orbit in this $2$-cube box with barriers into a $1$-direction geodesic in a boundary-free flat $3$-manifold.
The latter system is an $8$-copy construction involving $16$ cubes, and results from three consecutive reflections across a plane.

The original $2$-cube box with barriers in the middle is highlighted in bold in Figure~1.7.
We reflect it across the plane $x=2$, then reflect the $2$-copy union across the plane $y=1$,
and finally reflect the $4$-copy union across the plane $z=1$ to obtain an $8$-copy union.
Thus the original $2\times1\times1$ box becomes a $4\times2\times2$ box with two repeated sets of barriers
on the $yz$-parallel squares $[0,2)^2$ on the faces $x=1$ and $x=3$.

\begin{displaymath}
\begin{array}{c}
\includegraphics[scale=0.8]{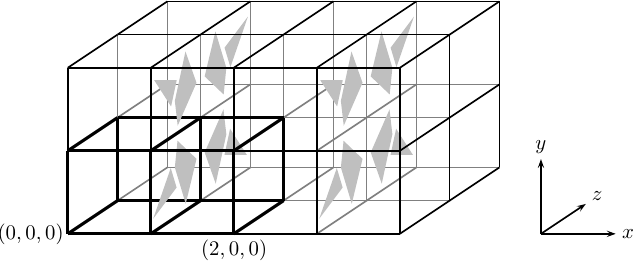}
\\
\mbox{Figure 1.7: unfolding the billiard orbit in a $2$-cube box}
\\
\mbox{with barriers in the middle}
\end{array}
\end{displaymath}

This box has boundary, and we turn it into a boundary-free flat $3$-manifold with boundary identification.
First of all, the square faces on the boundary of this box are identified by perpendicular translation.
Next, let $L$ and $R$ denote respectively the barriers on the square faces on the plane $x=1$ and $x=3$ respectively.
The left side of $L$ is identified with the right side of~$R$, while the right side of $L$ is identified with the left side of~$R$.
For convenience, we refer to this special $3$-manifold as the \textit{$2$-cube-billiard $3$-manifold}.

Clearly $1$-direction geodesic flow in the $2$-cube-billiard $3$-manifold is an $8$-fold cover of billiard flow in the $2$-cube box with barriers.

We now remove the part $[0,1)\times[0,2)\times[0,2)$ on the left and join it instead to the right to become $[4,5)\times[0,2)\times[0,2)$,
as shown in Figure~1.8.

\begin{displaymath}
\begin{array}{c}
\includegraphics[scale=0.8]{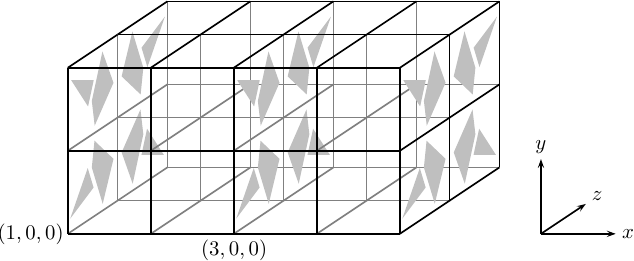}
\\
\mbox{Figure 1.8: cutting and pasting}
\end{array}
\end{displaymath}

Contracting the resulting $3$-manifold by a factor $1/2$ in each of the three directions,
we then obtain a $2$-cube $3$-manifold with repeated barriers that we have studied in Theorem~\ref{thm1}.
The following result is then a corollary of Theorem~\ref{thm1} in the special case $n=2$.

\begin{theorem}\label{thm3}
Consider billiard in a special $2$-cube box with barriers in the middle square face joining the cubes given by a $2$-coloring
such that each of the red and green parts is the union of finitely many polygons, and where the green part has positive area.
Then for almost every starting point and almost every initial direction $\bfv=(1,\alpha_1,\alpha_2)\in\Rr^3$
with $\bfv_0=(\alpha_1,\alpha_2)\in[-1,1]^2$, the corresponding half-infinite billiard orbit is equidistributed in this special $2$-cube box.
\end{theorem}

Unfortunately, Theorem~\ref{thm1} does not seem to help in the case of more complicated billiards with barriers.

Of course, there is no reason why the $2$-coloring on the distinct $yz$-parallel square faces
of the $n$-cube $3$-manifold should be the same, apart from possibly making the problem a little simpler.

We know that if a geodesic hits an $xy$-parallel square face of~$\MMM$,
then it jumps to the corresponding point on the identified $xy$-parallel square face and continues in the same direction,
and if a geodesic hits an $xz$-parallel square face, then it jumps to the corresponding point on the identified $xz$-parallel square face
and continues in the same direction.

Suppose now that a geodesic with direction $\bfv=(1,\alpha_1,\alpha_2)$ hits a $yz$-parallel square face at a point~$P$.
Then the continuation of the geodesic depends on the coloring of the intersection point~$P$.
If $P$ is green, then the geodesic continues on its way in the same direction.
If $P$ is red, then we consider a directed line starting from $P$ in the direction $(-1,0,0)$.
This line will hit a red point $P'$ for the first time.
Then the geodesic continues from $P'$ in the same direction, as shown in Figure~1.9.

\begin{displaymath}
\begin{array}{c}
\includegraphics[scale=0.8]{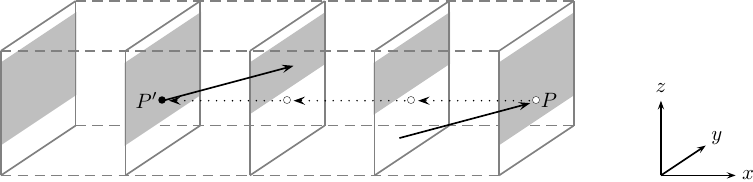}
\\
\mbox{Figure 1.9: when a geodesic hits a red barrier}
\end{array}
\end{displaymath}

There is the pathological case that $P'=P$, so that the geodesic continues on its way as if $P$ were green.
To avoid such situations, we deem the point $P$ to be colored green.
This is some kind of automatic recoloring, but the following rule is a little simpler.

\begin{no-red}
On any line perpendicular to any given $yz$-parallel square face of~$\MMM$,
there are either no points colored red or at least $2$ distinct points colored red.
Here we use the standard convention that the left-most $yz$-parallel square face
and the right-most $yz$-parallel square face are the same.
\end{no-red}

Instead of requiring perfect repetition of the $2$-coloring on all the $yz$-parallel square faces as in Theorem~\ref{thm1},
we now impose the substantially weaker condition of \textit{local repetition}.
More precisely, we require a small \textit{local repetition color-split neighborhood},
in the form of a line segment with the same local $2$-coloring of red and green in the two opposite side-neighborhoods.
For illustration, see Figure~1.10, where the three highlighted rectangles are in the same position within the square torus.

\begin{displaymath}
\begin{array}{c}
\includegraphics[scale=0.8]{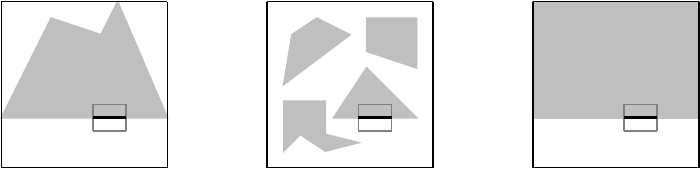}
\\
\mbox{Figure 1.10: different $2$-colorings with local repetition}
\end{array}
\end{displaymath}

We emphasize that the local repetition color-split neighborhood must be present on all the $yz$-parallel square faces of the $n$-cube $3$-manifold.
Naturally, we still need $2$-colorings on each $yz$-parallel square face such that each of the red and green parts is the union of
finitely many polygons, and where the green part has positive area.
Furthermore, we also require that the Restriction on Red Coloring holds.
Since the $2$-colorings on the different $yz$-parallel squares faces can now be different,
this represents substantially more freedom for the $2$-colorings.
In Figure~1.11, we have a $4$-cube $3$-manifold with local repetition color-split provided by the triangular red (shaded) regions
at the corners of the $yz$-parallel squares faces.
The positions of the local repetition color-split neighborhood on the different $yz$-parallel squares faces are indicated by the short thick lines.

\begin{displaymath}
\begin{array}{c}
\includegraphics[scale=0.8]{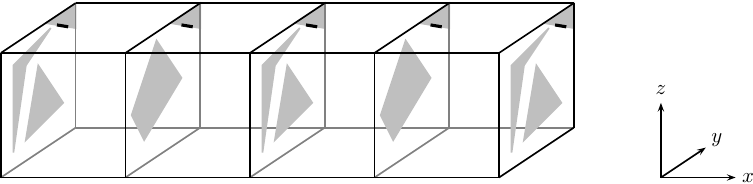}
\\
\mbox{Figure 1.11: a $4$-cube $3$-manifold with local repetition}
\end{array}
\end{displaymath}

The result that we can prove in this more general setting is not really weaker.

\begin{theorem}\label{thm4}
Let $n\ge2$ be an integer, and let $\MMM$ be any $n$-cube $3$-manifold with barriers,
where each $yz$-parallel square face has a $2$-coloring such that each of the red and green parts is the union of finitely many polygons,
and where the green part has positive area.
Suppose further that the Restriction on Red Coloring holds,
and that there is a local repetition color-split neighborhood on the $yz$-parallel square faces.
Then for almost every starting point and almost every direction $\bfv=(1,\alpha_1,\alpha_2)\in\Rr^3$,
the corresponding half-infinite $1$-direction geodesic is equidistributed in~$\MMM$.
\end{theorem}

As in Theorem~\ref{thm1}, any polygon of the $2$-colorings can be replaced by a circle, an ellipse, or any other piecewise smooth closed curve.

%
%

\section{Preparation for the proof of Theorem~\ref{thm4}}\label{sec2}

We work with an equivalent discrete form of the problem.

For every integer $i=0,1,\ldots,n-1$, let
\begin{displaymath}
U_i=\{i\}\times[0,1)^2
\end{displaymath}
denote the $i$-th $yz$-parallel square face of~$\MMM$, and let
\begin{displaymath}
X_0=\bigcup_{i=0}^{n-1}U_i.
\end{displaymath}
Since each $yz$-parallel square face has a $2$-coloring, for every integer $i=0,1,\ldots,n-1$, there exist two sets $\GGGG_i$ and $\RRRR_i$
such that
\begin{equation}\label{eq2.1}
U_i=\GGGG_i\cup\RRRR_i
\quad\mbox{and}\quad
\GGGG_i\cap\RRRR_i=\emptyset.
\end{equation}

Since the collection of non-Kronecker vectors $\bfv=(1,\alpha_1,\alpha_2)\in\Rr^3$ has measure zero,
we may therefore start our discussion by assuming that $\bfv$ is a Kronecker vector.
Let $\bfv=(1,\alpha_1,\alpha_2)\in\Rr^3$ be a Kronecker vector, and let $\bfv_0=(\alpha_1,\alpha_2)$.
We define an invertible transformation $\T=\T_\bfv:X_0\to X_0$ as follows.

For any point $P=(i,\bfy)\in U_i\subset X_0$ where $\bfy=(y,z)\in[0,1)^2$, let
\begin{equation}\label{eq2.2}
\T(P)=\left\{\begin{array}{ll}
(i+1,\{\bfy+\bfv_0\}),
&\mbox{if $P+\bfv=(i+1,\{\bfy+\bfv_0\})\in\GGGG_{i+1}$},\\
(i^*,\{\bfy+\bfv_0\}),
&\mbox{if $P+\bfv=(i+1,\{\bfy+\bfv_0\})\in\RRRR_{i+1}$},
\end{array}\right.
\end{equation}
where the addition in $P+\bfv$ is modulo~$n$ for the first coordinate and modulo~$1$ for the remaining coordinates.
Thus
\begin{equation}\label{eq2.3}
\{\bfy+\bfv_0\}=(\{y+\alpha_1\},\{z+\alpha_2\}),
\end{equation}
where $0\le\{\beta\}<1$ denotes the fractional part of a real number~$\beta$.
Furthermore, the value of $i^*$ is determined by
\begin{equation}\label{eq2.4}
i^*=\max\{j<i+1:(j,\{\bfy+\bfv_0\})\in\RRRR_j\},
\end{equation}
with the convention that
\begin{equation}\label{eq2.5}
i+1>i>i-1>\ldots>1>0>n-1>n-2>\ldots>i+2.
\end{equation}
Note that \eqref{eq2.4} and \eqref{eq2.5} are motivated by Figure~1.9, and Figure~2.1 illustrates the special case when $n=4$.
If $P+\bfv$ is a red point, we then move from $P+\bfv$ in the direction of the vector $(-1,0,0)$ until we hit the first red point,
and this red point lies on the $yz$-parallel square face~$U_{i^*}$.

\begin{displaymath}
\begin{array}{c}
\includegraphics[scale=0.8]{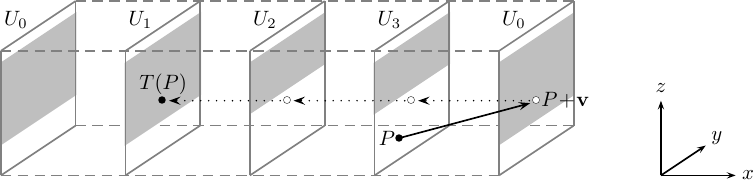}
\\
\mbox{Figure 2.1: when $P+\bfv$ is red}
\end{array}
\end{displaymath}

Projecting the transformation $\T$ by ignoring the first coordinate, we obtain the invertible transformation
\begin{displaymath}
\T_0=\T_{\bfv_0}:[0,1)^2\to[0,1)^2:\bfy\mapsto\{\bfy+\bfv_0\}.
\end{displaymath}

It is easy to see that $\T$ preserves the $2$-dimensional Lebesgue measure~$\lambda_2$.
Our goal is to establish that for almost every vector $\bfv_0\in[0,1]^2$, the transformation $\T=\T_\bfv:X_0\to X_0$ is ergodic, where $\bfv=(1,\bfv_0)$.
The basic idea is quite surprising, as we prove ergodicity for this non-integrable system by taking advantage of the split singularities.

\begin{remark}
We are going to use Birkhoff's well known pointwise ergodic theorem on measure preserving transformations.
Since we simply apply ergodic theory, we do not expect the reader to have any serious expertise in the subject.
Thus knowledge of Lebesgue integral and basic measure theory suffices.
The theorem concerns a measure-preserving system $(X,\AAA,\mu,\T)$.
Here $(X,\AAA,\mu)$ is a measure space, where $X$ is the underlying space, $\AAA$ is a $\sigma$-algebra of subsets of~$X$
and $\mu$ is a non-negative $\sigma$-additive measure on $X$ with $\mu(X)<\infty$,
while $\T:X\to X$ is a measurable map which is measure-preserving, so that $\T^{-1}A\in\AAA$ and $\mu(\T^{-1}A)=\mu(A)$ for every $A\in\AAA$.

Let $L^1(X,\AAA,\mu)$ denote the space of measurable and integrable functions in the measure space $(X,\AAA,\mu)$.
Then the general form of Birkhoff's pointwise ergodic theorem says that for every function $f\in L^1(X,\AAA,\mu)$, the limit
\begin{equation}\label{eq2.6}
\lim_{m\to\infty}\frac{1}{m}\sum_{j=0}^{m-1}f(\T^j\bfx)=f^*(\bfx)
\end{equation}
exists for $\mu$-almost every $\bfx\in X$, where $f^*\in L^1(X,\AAA,\mu)$ is a $\T$-invariant measurable function satisfying the condition
\begin{displaymath}
\int_Xf\,\dd\mu=\int_Xf^*\,\dd\mu.
\end{displaymath}

A particularly important special case is if $\T$ is \textit{ergodic}, when every measurable $\T$-invariant set $A\in\AAA$ is \textit{trivial}
in the precise sense that $\mu(A)=0$ or $\mu(A)=\mu(X)$.
This is equivalent to the assertion that every measurable $\T$-invariant function is constant $\mu$-almost everywhere.

If $\T$ is ergodic, then \eqref{eq2.6} simplifies to
\begin{equation}\label{eq2.7}
\lim_{m\to\infty}\frac{1}{m}\sum_{j=0}^{m-1}f(\T^j\bfx)=\int_Xf\,\dd\mu,
\end{equation}
and the right-hand side of \eqref{eq2.6} is the same constant for $\mu$-almost every $\bfx\in X$.

The remarkable intuitive interpretation of \eqref{eq2.7} is that the \textit{time average} on the left hand side
is equal to the \textit{space average} on the right hand side.
\end{remark}

%
%

\section{Proof of Theorem~\ref{thm4}}\label{sec3}

We focus on the particular measure-preserving system $(X_0,\AAA,\lambda_2,\T)$,
where $\AAA$ is the family of Borel sets in~$X_0$, $\lambda_2$ is $2$-dimensional Lebesgue measure and $\T=\T_\bfv$.
We shall establish ergodicity by contradiction.

\begin{step1}
Suppose on the contrary that $\T$ is not ergodic.
Then there exists a non-trivial measurable $\T$-invariant subset $S_0\subset X_0$ such that $0<\lambda_2(S_0)<n$.
We try to derive a contradiction.

Removing possibly a set of $\lambda_2$-measure zero, we may assume that for every point $\bfx\in X_0$,
the point $\T^j(\bfx)$ is well defined for every integer $j=1,2,3,\ldots.$

\begin{lemma}\label{lem31}
Consider the measure-preserving system $(X_0,\AAA,\lambda_2,\T)$, where $\T=\T_\bfv$ and $\bfv=(1,\bfv_0)$ is a Kronecker vector.
For any $\T$-invariant subset $S_0\subset X_0$, let the multiplicity function $\widetilde{\chi}_{S_0}$ of $S_0$
be defined for every point $P\in[0,1)^2$ by
\begin{displaymath}
\widetilde{\chi}_{S_0}(P)=\vert\{i=0,1,\ldots,n-1:(i,P)\in S_0\}\vert.
\end{displaymath}
Suppose further that $S_0$ is a proper subset of $X_0$, so that $S_0\ne\emptyset$ and $S_0\ne X_0$.
Then there exists an integer $k_0=1,\ldots,n-1$ such that $\widetilde{\chi}_{S_0}(P)=k_0$ for almost every point $P\in[0,1)^2$,
so that $\lambda_2(S_0)=k_0$.
\end{lemma}

\begin{proof}
Since $\bfv$ is a Kronecker vector, it follows that the $\bfv_0$-shift on the unit torus $[0,1)^2$ is ergodic.
Meanwhile, it is easy to check that the multiplicity function $\widetilde{\chi}_{S_0}$ is $\T_0$-invariant.
Thus Birkhoff's ergodic theorem implies that $\widetilde{\chi}_{S_0}$ is constant almost everywhere.
Note that $\widetilde{\chi}_{S_0}$ is integer valued and cannot be equal to $0$ or~$n$.
This completes the proof.
\end{proof}

\end{step1}

\begin{step2}
Given a point $\bfz\in X_0$ and a radius $0<r<1/2$, let $D(\bfz;r)$ denote the circular disk of radius $r$ and center~$\bfz$.
Clearly $D(\bfz;r)$ has area~$\pi r^2$.
Note that $D(\bfz;r)\subset X_0$, due to the fact that $X_0$ is a compact flat surface.

Since the non-trivial $\T$-invariant subset $S_0\subset X_0$ is measurable, it follows from Lebesgue's density theorem
that for almost every $\bfz\in S_0$,
\begin{displaymath}
\lim_{r\to0}\frac{\lambda_2(S_0\cap D(\bfz;r))}{\pi r^2}=1,
\end{displaymath}
whereas for almost every $\bfz\in S_0^c=X_0\setminus S_0$,
\begin{displaymath}
\lim_{r\to0}\frac{\lambda_2(S_0\cap D(\bfz;r))}{\pi r^2}=0.
\end{displaymath}

Let $M$ be a large integer, and divide each of $U_0,U_1,\ldots,U_{n-1}$ into $M^2$ congruent squares of area $(1/M)^2$ in the standard way.
We refer to these small squares as special $(1/M)$-squares.
Thus there are precisely $nM^2$ special $(1/M)$-squares in~$X_0$.

In view of Lebesgue's density theorem, we formulate here and prove in Section~\ref{sec4}
the following lemma for the hypothetical non-trivial measurable $\T$-invariant subset $S_0\subset X_0$.

\begin{lemma}\label{lem32}
Let the real number $\eps\in(0,1)$ be arbitrarily small and fixed, and let the real number $\eps_1>0$ be fixed.
There exists a finite threshold $m_0=m_0(S_0;\eps;\eps_1)$ such that for every integer $M\ge m_0$, there exist at least $(1-\eps_1)nM^2$ special
$(1/M)$-squares $Q$ in $X_0$ such that either
\begin{displaymath}
\frac{\lambda_2(S_0\cap Q)}{(1/M)^2}>1-\eps
\quad\mbox{or}\quad
\frac{\lambda_2(S_0\cap Q)}{(1/M)^2}<\eps.
\end{displaymath}
\end{lemma}

Let $N$ be a large even integer.
Let $\FFF(N/2)$ denote the standard decomposition of the unit torus $[0,1)^2$ into $(N/2)^2$ axis-parallel congruent small squares
of common side length~$2/N$ such that the origin $(0,0)$ is the vertex of a small square.
For $\bfdelta=(\delta_1,\delta_2)\in\{0,1\}^2$, let $\FFF_\bfdelta(N/2)$ denote the translation of $\FFF(N/2)$ modulo one
such that the vertex $(0,0)$ moves to $(\delta_1/N,\delta_2/N)$.
We refer to the small squares in the four partitions $\FFF_\bfdelta(N/2)$, $\bfdelta\in\{0,1\}^2$, as basic $(2/N)$-squares.
It is not difficult to see that any axis-parallel square $B$ of side length $1/N$ in the unit torus is contained in a basic $(2/N)$-square.

For every $i=0,1,\ldots,n-1$, we replicate the families $\FFF_\bfdelta(N/2)$, $\bfdelta\in\{0,1\}^2$,
on the $yz$-parallel square face $U_i=\{i\}\times[0,1)^2$.
In other words, we write
\begin{displaymath}
\FFF_\bfdelta(N/2;i)=\{i\}\times\FFF_\bfdelta(N/2),
\quad
\bfdelta\in\{0,1\}^2,
\quad
i=0,1,\ldots,n-1.
\end{displaymath}
For every $\bfdelta\in\{0,1\}^2$, Lemma~\ref{lem32} with $M=N/2$ then gives the following.

Let the real number $\eps\in(0,1)$ be arbitrarily small and fixed, and let the real number $\eps_1>0$ be fixed.
There exists a finite threshold $m_0=m_0(S_0;\eps;\eps_1)$ such that for every $\bfdelta\in\{0,1\}^2$ and every even integer $N\ge m_0$,
there exist at least $(1-\eps_1)n(N/2)^2$ basic $(2/N)$-squares
\begin{equation}\label{eq3.1}
Q\in\bigcup_{i=0}^{n-1}\FFF_\bfdelta(N/2;i)
\end{equation}
such that either
\textcolor{white}{xxxxxxxxxxxxxxxxxxxxxxxxxxxxxx}
\begin{displaymath}
\frac{\lambda_2(S_0\cap Q)}{(2/N)^2}>1-\eps
\quad\mbox{or}\quad
\frac{\lambda_2(S_0\cap Q)}{(2/N)^2}<\eps.
\end{displaymath}
We say that these sets $Q$ satisfy the \textit{$\eps$-nearly zero-one law}.
The conclusion of this step is that we have a lower bound for the number of such basic $(2/N)$-squares of the form \eqref{eq3.1}.
\end{step2}

\begin{step3}
To apply a version of the splitting method, first introduced in \cite{BCY2}, we need to make use of the
local repetition color-split neighborhood~$\SSS$, with green part $\SSS_\GGGG$ and red part $\SSS_\RRRR$.
As this color-split neighborhood is present on all the $yz$-parallel square faces of the $n$-cube $3$-manifold~$\MMM$,
we assume, for simplicity, that $\SSS$ lies on the unit torus $[0,1)^2$.

We can clearly find within $\SSS$ a buffer zone $\BBB_N$ as shown in Figure~3.1.
The length of $\BBB_N$ is a constant $c_1$ depending only in $\SSS$, while the width is $1/2N$, where the integer $N$ is sufficiently large,
with the color-split boundary splitting $\BBB_N$ into green and red strips of width~$1/4N$.

\begin{displaymath}
\begin{array}{c}
\includegraphics[scale=0.8]{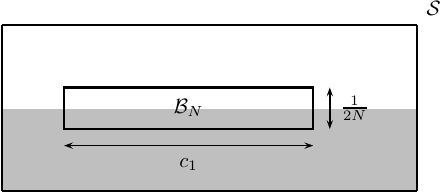}
\\
\mbox{Figure 3.1: buffer zone $\BBB_N$ within the local repetition color-split neighborhood $\SSS$}
\end{array}
\end{displaymath}

The following result follows from simple geometric considerations.

\begin{lemma}\label{lem33}
Let $A_0\subset[0,1)^2$ be an arbitrary axis-parallel square with side length $1/N$ and center $c(A_0)$.
Then for every $\bft\in[0,1)^2$ such that $c(A_0)+\bft\in\BBB_N$, the set $A_0+\bft$ has substantial color-split in the sense that
\begin{displaymath}
\lambda_2((A_0+\bft)\cap\SSS_\GGGG)\ge\frac{1}{14N^2}
\quad\mbox{and}\quad
\lambda_2((A_0+\bft)\cap\SSS_\RRRR)\ge\frac{1}{14N^2}.
\end{displaymath}
\end{lemma}

The extreme case of the above takes place when $\BBB_N$ is tilted at $45$ degrees, with the center $c(A_0)$ of $A_0$ at the corner,
as shown in Figure~3.2.

\begin{displaymath}
\begin{array}{c}
\includegraphics[scale=0.8]{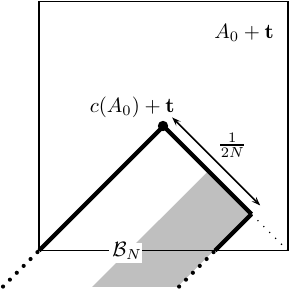}
\\
\mbox{Figure 3.2: extreme case of substantial color-split}
\end{array}
\end{displaymath}

We need to show that such substantial color-splits occur quite frequently.
We shall establish in Section~\ref{sec5} the following result.

\begin{lemma}\label{lem34}
Suppose that the integer $N$ is even and sufficiently large.
For every pair $\bfs_0\in[0,1)^2$ and $\bfv_0\in[0,1]^2$, let $F(\bfs_0;\bfv_0;N^2)$ denote the number of integers $j=0,1,\ldots,N^2-1$
such that $\bfs_0+j\bfv_0\in\BBB_N$.
Then for every $\eps_2>0$, there exists a constant $c_2>0$ such that
\begin{equation}\label{eq3.2}
\lambda_2(\{\bfv_0\in[0,1)^2:F(\bfs_0;\bfv_0;N^2)\ge c_2N\})\ge1-\eps_2.
\end{equation}
\end{lemma}

\end{step3}

\begin{step4}
Consider a $(1/N)$-square $A$ on some $yz$-parallel square face $U_i$  of~$\MMM$.
We see from \eqref{eq2.1}--\eqref{eq2.5} that the image $\T(A)$ is then given by $A+\bfv$
followed by appropriate bounce-backs on those parts that hit red.
If $A+\bfv$ hits both green and red, then clearly the splitting of the image $\T(A)$ is caused by a \textit{color-split}.
For instance, in Figure~3.3, the white square on the $yz$-parallel square face $U_1$ denotes~$A$,
the square $A+\bfv$ experiences a color-split on the $yz$-parallel square face~$U_2$,
and the image $\T(A)$, indicated in black, is split between the $yz$-parallel square faces $U_1$ and~$U_2$.
Indeed, if the $2$-coloring on the different $yz$-parallel square faces of $\MMM$ are identical,
then the only splitting of image is caused by a color-split.

\begin{displaymath}
\begin{array}{c}
\includegraphics[scale=0.8]{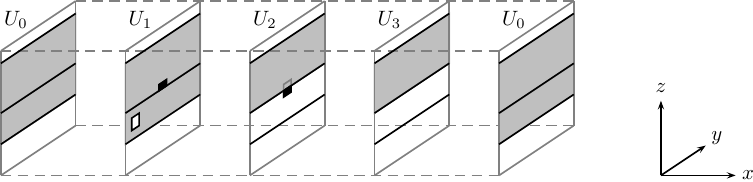}
\\
\mbox{Figure 3.3: a color-split}
\end{array}
\end{displaymath}

On the other hand, if the different $yz$-parallel square faces of $\MMM$ can have different $2$-colorings,
then there are other instances that cause splitting of the image $\T(A)$.
For instance, in Figure~3.4, for the white square $A$ on the $yz$-parallel square face~$U_3$,
the square $A+\bfv$ is on the $yz$-parallel square face~$U_0$, but then different parts have different bounce-backs,
and the image $\T(A)$ is split between the $yz$-parallel square faces $U_1$ and~$U_3$.
Although a color-split always leads to an image-split, this example shows that the converse is not true.

\begin{displaymath}
\begin{array}{c}
\includegraphics[scale=0.8]{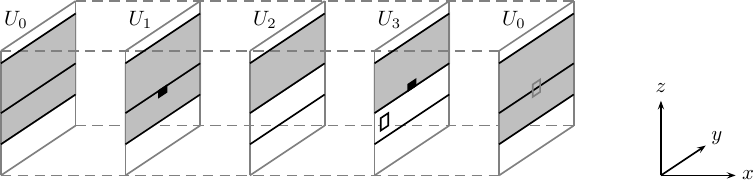}
\\
\mbox{Figure 3.4: an image-split that is not a color-split}
\end{array}
\end{displaymath}

It turns out that it is advantageous to consider another kind of splitting which may or may not result in an image-split.
For every $i=0,1,\ldots,n-1$, let $\{i\}\times\Gamma_i$ denote the color-split boundaries of the $2$-coloring on~$U_i$,
so that $\Gamma_i\subset[0,1)^2$ is a collection of straight edges.
Consider the union
\begin{displaymath}
\Gamma=\bigcup_{i=0}^{n-1}\Gamma_i,
\end{displaymath}
and, for every $i=0,1,\ldots,n-1$, replace $\{i\}\times\Gamma_i$ on $U_i$ by $\{i\}\times\Gamma$,
so that the different $yz$-parallel square faces of $\MMM$ now have identical \textit{$\Gamma$-splits}.
Although an image-split is always a $\Gamma$-split, the example in Figure~3.5 shows that the converse is not true.
Here, and also in Figures 3.3 and~3.4, the $\Gamma$-splits are indicated by the black line segments.

\begin{displaymath}
\begin{array}{c}
\includegraphics[scale=0.8]{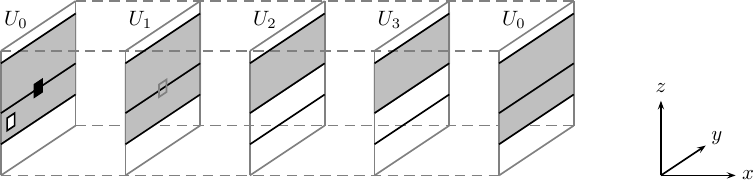}
\\
\mbox{Figure 3.5: a $\Gamma$-split that is not an image-split}
\end{array}
\end{displaymath}

The important point here is that if there is no $\Gamma$-split, then there is no splitting of any kind and in particular no color-split.

For the splitting method to work, we need a substantial color-split in the local repetition color-split neighborhood
and a \textit{long} color-split-free chain either side of it.
In view of the above observation, it therefore suffices to find a substantial color-split in the local repetition color-split neighborhood
and a \textit{long} $\Gamma$-split-free chain either side of it.
As the $\Gamma$-split is identical on all the $yz$-parallel square faces of~$\MMM$, this allows us to consider $\bfv_0$-flow
on the unit torus $[0,1)^2$ instead of geodesic flow in the direction $\bfv$ in~$\MMM$.

To establish this, we shall show that the number of \textit{short} $\Gamma$-split-free chains is \textit{small}.
More precisely, we have the following lemma which we prove in Section~\ref{sec4}.
Before we can state the result, we first need some definitions.

Let $A_0\subset[0,1)^2$ be an arbitrary axis-parallel square of side length~$1/N$.
Consider the $N^2$ $\bfv_0$-shift images
\textcolor{white}{xxxxxxxxxxxxxxxxxxxxxxxxxxxxxx}
\begin{equation}\label{eq3.3}
A_0+j\bfv_0,
\quad
j=0,1,\ldots,N^2-1,
\end{equation}
in $[0,1)^2$.
We say that a set $\{A_0+j\bfv_0:j\in J\}$, where $J$ is a subset of consecutive integers in $\{0,1,\ldots,N^2-1\}$,
is $\Gamma$-split-free if $(A_0+j\bfv_0)\cap\Gamma=\emptyset$ for every $j\in J$.
Furthermore, we say that the set $\{A_0+j\bfv_0:j\in J\}$ is a $\Gamma$-split-free chain if it is $\Gamma$-split-free
and not contained in a bigger $\Gamma$-split-free set.
It is convenient to define the length of the $\Gamma$-split-free chain to be $\vert J\vert+1$.

Thus the sequence \eqref{eq3.3} decomposes into a subsequence of $\Gamma$-split members,
with any two consecutive members of this subsequence possibly separated by a $\Gamma$-split-free chain in between.

For every axis-parallel square $A_0\subset[0,1)^2$ of side length $1/N$, we can clearly write $A_0=A_0(\bfs_0)$,
where $\bfs_0=c(A_0)$ is the center of~$A_0$.
Also, for every vector $\bfv=(1,\bfv_0)$, where $\bfv_0\in\Rr^2$, we can identify
\begin{displaymath}
\bfv,
\quad
\bfv_0
\quad\mbox{and}\quad
\bfv^*=\frac{\bfv}{\vert\bfv\vert}
\end{displaymath}
with each other.
Clearly the collection of all vectors $\bfv^*$, where $\bfv_0\in\Rr^2$,
forms the subset $\bfS^2_+=\{(x,y,z)\in\bfS^2:x>0\}$ of the unit ball~$\bfS^2$.

\begin{lemma}\label{lem35}
Let the real number $\eps_2\in(0,1)$ be arbitrarily small and fixed,
and let $A_0\subset[0,1)^2$ be an arbitrary axis-parallel square of side length $1/N$ that runs uniformly over $[0,1)^2$,
in the sense that the center $\bfs_0$ is uniformly distributed on $[0,1)^2$.
Then there exists a positive absolute constant $c_3$ such that for at least $(1-\eps_2)$-proportion of the pairs
\textcolor{white}{xxxxxxxxxxxxxxxxxxxxxxxxxxxxxx}
\begin{displaymath}
(\bfs_0,\bfv^*)\in[0,1)^2\times\bfS^2_+,
\end{displaymath}
the sequence \eqref{eq3.3} in $[0,1)^2$ contains at most $\eps_2N$ $\Gamma$-split-free chains with length at most $c_3\eps_2^2N$.
The $(1-\eps_2)$-proportion is in terms of the product of the $2$-dimensional Lebesgue measure on $[0,1)^2$
and the normalized surface area measure on~$\bfS^2_+$.
\end{lemma}

\end{step4}

\begin{step5}
The motivation for this step is not immediately obvious at this stage, but we include it here for the sake of convenience.
This will become clear in Step~7.

We establish in Section~\ref{sec4} the following result on \textit{clustering}.
It demonstrates that the overwhelming majority of arithmetic progressions are not clustered.

\begin{lemma}\label{lem36}
Let the real number $\eps_2\in(0,1)$ be arbitrarily small and fixed.
There exists a finite constant $C^*=C^*(\eps_2)$ such that for any starting point $\bfs_0$ in the unit torus $[0,1)^2$
and for at least $(1-\eps_2)$-proportion of the vectors $\bfv_0\in[0,1]^2$,
every axis-parallel square $Q$ of side length $2/N$ in $[0,1)^2$ contains modulo one at most $C^*$ elements of the arithmetic progression
\begin{displaymath}
\bfs_0+j\bfv_0,
\quad
j=0,1,\ldots,N^2-1,
\end{displaymath}
of $N^2$ terms.
In particular, we can take $C^*=1+16/\eps_2$.
\end{lemma}

\end{step5}

\begin{step6}
Here we combine our observations in Steps 3--5.
Combining Lemmas \ref{lem33}, \ref{lem34}, \ref{lem35} and~\ref{lem36}, we see that there exists a set $\bfG\bfD\subset[0,1]^2$
of \textit{good} directions $\bfv_0\in[0,1]^2$ with $\lambda_2(\bfG\bfD)\ge1-3\eps_2$
such that for every $\bfv_0\in\bfG\bfD$, the conclusion of Lemma~\ref{lem36} holds for any $\bfs_0\in[0,1)^2$,
and there exists an axis-parallel square $A_0^*$ in $[0,1)^2$ of side length $1/N$ such that the following two conditions hold:

(i)
The sequence
\textcolor{white}{xxxxxxxxxxxxxxxxxxxxxxxxxxxxxx}
\begin{equation}\label{eq3.4}
A_0^*+j\bfv_0,
\quad
j=0,1,\ldots,N^2-1,
\end{equation}
in $[0,1)^2$ contains at least $c_2N$ members which exhibit substantial color-split in the local repetition color-split neighborhood~$\SSS$,
in the sense that
\begin{displaymath}
\lambda_2((A_0^*+j\bfv_0)\cap\SSS_\GGGG)\ge\frac{1}{14N^2}
\quad\mbox{and}\quad
\lambda_2((A_0^*+j\bfv_0)\cap\SSS_\RRRR)\ge\frac{1}{14N^2}.
\end{displaymath}

(ii)
The sequence \eqref{eq3.4} has at most $\eps_2N$ \textit{short} $\Gamma$-split-free chains with length at most $c_3\eps_2^2N$.

Combining (i) and (ii), we conclude that the sequence \eqref{eq3.4} has at least
\begin{equation}\label{eq3.5}
\frac{c_2N-2\eps_2N-2}{2}\ge\frac{c_2N}{3}
\end{equation}
pairs of consecutive \textit{long} $\Gamma$-split-free chains of length at least $c_3\eps_2^2N$,
where each pair is separated by a member of the sequence with substantial color-split in the local repetition color-split neighborhood~$\SSS$.
Here we require that $\eps_2>0$ is chosen to be arbitrarily small so that the inequality \eqref{eq3.5} holds.
\end{step6}

\begin{step7}
To facilitate the use of the splitting method, we need to first demonstrate the existence of a pair of consecutive long $\Gamma$-split-free chains
separated by $A_0^*+j_0\bfv_0$ with substantial color-split in the local repetition color-split neighborhood~$\SSS$ and such that
both $A_0^*+(j_0-1)\bfv_0$ and $A_0^*+(j_0+1)\bfv_0$ are $\Gamma$-split-free and satisfy the $\eps$-nearly zero-one law.
More precisely, we return to the $n$-cube $3$-manifold~$\MMM$.
Then $A_0^*+j_0\bfv_0$ corresponds to an axis-parallel square $A^{**}$ of side length $1/N$ on some $yz$-parallel square face~$U_i$,
while $A_0^*+(j_0-1)\bfv_0$ and $A_0^*+(j_0+1)\bfv_0$ correspond to the squares $A^{**}-\bfv$ and $A^{**}+\bfv$ respectively.
Then we claim that either
\begin{equation}\label{eq3.6}
\frac{\lambda_2(S_0\cap(A^{**}-\bfv))}{\lambda_2(A^{**}-\bfv)}>1-\eps
\quad\mbox{or}\quad
\frac{\lambda_2(S_0\cap(A^{**}-\bfv))}{\lambda_2(A^{**}-\bfv)}<\eps,
\end{equation}
and either
\textcolor{white}{xxxxxxxxxxxxxxxxxxxxxxxxxxxxxx}
\begin{equation}\label{eq3.7}
\frac{\lambda_2(S_0\cap(A^{**}+\bfv))}{\lambda_2(A^{**}+\bfv)}>1-\eps
\quad\mbox{or}\quad
\frac{\lambda_2(S_0\cap(A^{**}+\bfv))}{\lambda_2(A^{**}+\bfv)}<\eps.
\end{equation}

To establish our claim, suppose, on the contrary, that for every pair of consecutive long $\Gamma$-split-free chains
separated by $A_0^*+j_0\bfv_0$ with substantial color-split in the local repetition color-split neighborhood~$\SSS$,
either \eqref{eq3.6} or \eqref{eq3.7} fails.
Without loss of generality, suppose that \eqref{eq3.7} fails, so that $A^{**}+\bfv$ fails the $\eps$-nearly zero-one law.
Note that $A^{**}+\bfv$ corresponds to the first member of a $\Gamma$-split-free chain of at least $c_3\eps_2^2N-2$ terms.
Since the subset $S_0\subset X_0$ is $\T$-invariant, it follows that every set
\begin{equation}\label{eq3.8}
\T^j(A^{**}+\bfv),
\quad
0\le j\le c_3\eps_2^2N-3,
\end{equation}
is an axis-parallel square of side length $1/N$ on a $yz$-parallel square face of $\MMM$
that is $\Gamma$-split-free and fails the $\eps$-nearly zero-one law.
Taking into consideration all pairs of consecutive long $\Gamma$-split-free chains,
we see that there are at least
\begin{displaymath}
\frac{c_2N}{3}\cdot(c_3\eps_2^2N-2)
\end{displaymath}
such axis-parallel squares of the form \eqref{eq3.8} which fail the $\eps$-nearly zero-one law.
Each is contained in a basic $(2/N)$-square $Q$ of the form \eqref{eq3.1} for some $\bfdelta\in\{0,1\}^2$.
On the other hand, in view of Lemma~\ref{lem36}, any basic $(2/N)$-square $Q$ contains at most $C^*=C^*(\eps_2)$
such axis-parallel squares of the form \eqref{eq3.8} which fail the $\eps$-nearly zero-one law.
We conclude therefore that there are at least
\begin{displaymath}
\frac{c_2N(c_3\eps_2^2N-2)}{3C^*}
\end{displaymath}
distinct basic $(2/N)$-squares $Q$ of the form \eqref{eq3.1} for some $\bfdelta\in\{0,1\}^2$ which fail the $\eps$-nearly zero-one law.
Meanwhile, the conclusion from Step~2 is that the number of such basic $(2/N)$-squares $Q$ is bounded above by
\begin{displaymath}
4\eps_1n\left(\frac{N}{2}\right)^2.
\end{displaymath}
For any given $\eps_2>0$, we can now choose $\eps_1>0$ sufficiently small in terms of the other parameters to obtain a contradiction.
This establishes the claim.
\end{step7}

\begin{step8}
We now apply the splitting method.
Assume that $A^{**}$ is an axis-parallel square of side length $1/N$ on some $yz$-parallel square face of~$\MMM$,
and that \eqref{eq3.6} and \eqref{eq3.7} both hold.
For every $i=0,1,\ldots,n-1$, let $A_i\subset U_i$ be the axis-parallel square of side length $1/N$
that is in the same general position as $A^{**}$, and let
\begin{displaymath}
B_{i-1}=A_i-\bfv
\quad\mbox{and}\quad
C_{i+1}=A_i+\bfv.
\end{displaymath}
Note that each $A_i$ falls within the local repetition color-split neignborhood~$\SSS$, with local $2$-coloring $\SSS_\GGGG$ and $\SSS_\RRRR$,
and that each $B_{i-1}$ and each $C_{i+1}$ is $\Gamma$-split-free and so monochromatic,
and also $\eps$-nearly in $S_0$ or $\eps$-nearly outside~$S_0$.

For each $i=0,1,\ldots,n-1$, let $B_{i-1}(+)\subset B_{i-1}$ be defined by
\begin{displaymath}
\T(B_{i-1}(+))=A_i\cap\SSS_\GGGG,
\end{displaymath}
so that its image under $\T$ is precisely the green part of~$A_i$, and let
\begin{displaymath}
B_{i-1}(-)=B_{i-1}\setminus B_{i-1}(+).
\end{displaymath}
Suppose that $C_{i+1}\subset\GGGG_{i+1}$, so that it is on the green part of~$U_{i+1}$.
Then it follows from the definition of~$\T$, as given by \eqref{eq2.1}--\eqref{eq2.5}, that
\begin{displaymath}
C_{i+1}=\T^2(B_{i-1}(+))\cup \T^2(B_i(-)).
\end{displaymath}
Suppose that $C_{i+1}\subset\RRRR_{i+1}$, so that it is on the red part of~$U_{i+1}$.
Corresponding to \eqref{eq2.4} and using the convention \eqref{eq2.5}, let
\begin{displaymath}
i^*=\max\{j<i+1:C_j\subset\RRRR_j\}.
\end{displaymath}
Then it follows from the definition of~$\T$, as given by \eqref{eq2.1}--\eqref{eq2.5}, that
\begin{displaymath}
C_{i^*}=\T^2(B_{i-1}(+))\cup \T^2(B_i(-)).
\end{displaymath}
In either case, the two sets $B_{i-1}(+)$ and $B_i(-)$, and so also the two sets $B_{i-1}$ and~$B_i$,
are either both $\eps$-nearly in $S_0$ or both $\eps$-nearly outside~$S_0$.
Thus the sets
\begin{displaymath}
B_0,\ldots,B_{n-1}
\quad\mbox{and}\quad
C_0,\ldots,C_{n-1}
\end{displaymath}
are either all $\eps$-nearly in $S_0$ or all $\eps$-nearly outside~$S_0$.
This contradicts the assertion in Step~1 that $\widetilde{\chi}_{S_0}(P)=k_0$ for almost every point $P\in[0,1)^2$
for some integer $k_0<n$.

We can guarantee $\eps_2\to0$ by taking $N\to\infty$.
This establishes that $\T=\T_\bfv$ is ergodic, where $\bfv=(1,\bfv_0)$, for almost every $\bfv_0\in[0,1]^2$.
\end{step8}

\begin{step9}
We can now use the standard technique of extending ergodicity to unique ergodicity using functional analysis and Borel measures.
This completes the proof of Theorem~\ref{thm4}.
\end{step9}

%
%

\section{Proof of Lemmas \ref{lem32}, \ref{lem35} and~\ref{lem36}}\label{sec4}

\begin{proof}[Proof of Lemma~\ref{lem32}]
Since $S_0$ is Lebesgue measurable, given any $\delta>0$, there exists a finite set of disjoint axis-parallel rectangles
such that their union $V$ satisfies
\begin{displaymath}
\lambda_2(V\setminus S_0)+\lambda_2(S_0\setminus V)<\delta.
\end{displaymath}
Suppose that $0<\lambda_2(S_0)=\tau<n$.
Then
\begin{displaymath}
\lambda_2(V)>\lambda_2(S_0)-\delta=\tau-\delta
\quad\mbox{and}\quad
\lambda_2(S_0^c\cap V)<\delta,
\end{displaymath}
where $S_0^c=X_0\setminus S_0$.
Since $V$ is a finite union of disjoint axis-parallel rectangles, there clearly exists a threshold $t_1=t_1(V;\delta)$
such that the union $V_1$ of the special $(1/t)$-squares $Q$ contained in $V$ has measure
\begin{displaymath}
\lambda_2(V_1)>\lambda_2(V)-\delta>\tau-2\delta,
\end{displaymath}
provided that the integer $t\ge t_1$.
Let $\BBB$ denote the set of special $(1/t)$-squares $Q$ in $V_1$ that satisfy
\textcolor{white}{xxxxxxxxxxxxxxxxxxxxxxxxxxxxxx}
\begin{displaymath}
\frac{\lambda_2(S_0^c\cap Q)}{(1/t)^2}\ge\eps.
\end{displaymath}
Then, provided that $t\ge t_1(V;\delta)$, we have
\begin{displaymath}
\delta
>\lambda_2(S_0^c\cap V)
\ge\lambda_2 (S_0^c\cap V_1)
=\sum_{Q\subset V_1}\lambda_2(S_0^c\cap Q)\
\ge\sum_{Q\in\BBB}\lambda_2(S_0^c\cap Q)
\ge\frac{\eps\vert\BBB\vert}{t^2},
\end{displaymath}
so that
\textcolor{white}{xxxxxxxxxxxxxxxxxxxxxxxxxxxxxx}
\begin{displaymath}
\vert\BBB\vert\le\frac{\delta t^2}{\eps}=\frac{\eps_1t^2}{6},
\end{displaymath}
if we choose $\delta=\eps\eps_1/6$.
Deleting the special $(1/t)$-squares $Q\in\BBB$, we see that $V_1$ contains at least
\textcolor{white}{xxxxxxxxxxxxxxxxxxxxxxxxxxxxxx}
\begin{equation}\label{eq4.1}
\left(\tau-2\delta-\frac{\eps_1}{6}\right)t^2\ge\left(\tau-\frac{\eps_1}{2}\right)t^2
\end{equation}
special $(1/t)$-squares $Q$ such that
\textcolor{white}{xxxxxxxxxxxxxxxxxxxxxxxxxxxxxx}
\begin{displaymath}
\frac{\lambda_2(S_0^c\cap Q)}{(1/t)^2}<\eps.
\end{displaymath}
It follows that, as long as the integer $t\ge t_1(V;\delta)$, the number of special $(1/t)$-squares in $X_0$ that satisfy
\textcolor{white}{xxxxxxxxxxxxxxxxxxxxxxxxxxxxxx}
\begin{equation}\label{eq4.2}
\frac{\lambda_2(S_0\cap Q)}{(1/t)^2}>1-\eps
\end{equation}
is bounded below by \eqref{eq4.1}.
Repeating the same argument but replacing $S_0$ by $S_0^c$, we obtain another threshold $t_2=t_2(V;\delta)$ such that,
as long as the integer $t\ge t_2(V;\delta)$, the number of special $(1/t)$-squares in $X_0$ that satisfy
\begin{equation}\label{eq4.3}
\frac{\lambda_2(S_0\cap Q)}{(1/t)^2}<\eps
\end{equation}
is bounded below by
\textcolor{white}{xxxxxxxxxxxxxxxxxxxxxxxxxxxxxx}
\begin{equation}\label{eq4.4}
\left(n-\tau-\frac{\eps_1}{2}\right)t^2.
\end{equation}
Combining the lower bounds \eqref{eq4.1} and \eqref{eq4.4}, we see that, provided that an integer $M\ge\max\{t_1,t_2\}$,
the number of special $(1/M)$-squares in $X_0$ that satisfy \eqref{eq4.2} or \eqref{eq4.3}
is bounded below by $(n-\eps_1)M^2$, and this completes the proof.
\end{proof}

\begin{proof}[Proof of Lemma~\ref{lem35}]
Throughout this proof, the parameters $C_1,C_2,C_3,\ldots$ represent positive absolute constants.

Let $\eta=\eta(\eps_2)>0$, to be fixed later.

Suppose that an integer~$j^*$, satisfying $0<j^*\le\eta N$, is such that both $A_0$ and $A_0+j^*\bfv_0$ exhibit $\Gamma$-split,
but $A_0+j\bfv_0$ for any integer $j$ satisfying $0<j<j^*$ is $\Gamma$-split-free.
This gives rise to a $\Gamma$-split-free chain of length~$j^*$, and we say that $\bfv_0$ is a \textit{bad} direction for~$A_0$.

We now go to $3$-space as follows.
A discrete point sequence
\begin{displaymath}
\bfs_0+j\bfv_0,
\quad
j=0,1,\ldots,N^2-1,
\end{displaymath}
in the torus $[0,1)^2$ defines a straight line
\begin{displaymath}
\bfs+t\bfv\in\Rr^3,
\end{displaymath}
where $\bfv=(1,\bfv_0)$.
It follows that the sequence \eqref{eq3.3} leads to a sequence
\begin{equation}\label{eq4.5}
A+j\bfv,
\quad
j=0,1,\ldots,N^2-1,
\end{equation}
in $\Rr^3$, where $A=\{x_0\}\times A_0$ is an axis-parallel square of side length $1/N$ on some plane $x=x_0$ where $x_0$ is an integer.
The $\Gamma$-split-free chain under consideration is then characterized by a vector from the center of $A$ to the centre of $A+j^*\bfv$ in~$\Rr^3$.
Clearly the direction of any such vector is given by $\bfv$, and the length is at most $\sqrt{3}\eta N$, where we assume that $\sqrt{3}\eta N>1$.

Let the integer $m$ satisfy
\textcolor{white}{xxxxxxxxxxxxxxxxxxxxxxxxxxxxxx}
\begin{equation}\label{eq4.6}
2^{m-1}<\sqrt{3}\eta N\le2^m.
\end{equation}

Our basic idea is straightforward.
We consider all $yz$-parallel unit squares arising from the lattice~$\Zz^3$,
and extend $\Gamma$ to each in the usual way.
Assume that $A$ is a $\Gamma$-split $(1/N)$-square on some $yz$-parallel unit square.
Consider some other $\Gamma$-split $(1/N)$-square $A^*$ on some other $yz$-parallel unit square.
If the center of $A^*$ is reached from the center of $A$ via a vector in the direction of $\bfv$ and of length at most~$2^m$,
then $A^*$ determines a bad direction $\bfv$ for~$A$.
Theoretically, we can determine all possible $\Gamma$-split $(1/N)$-squares $A^*$ on other $yz$-parallel unit squares,
the centers of which can be reached from the center of $A$ via vectors of length at most~$2^m$,
leading to a collection of bad directions for~$A$.
Carrying this out, however, is impossible without more care.

Clearly every $yz$-parallel unit square has $N^2$ special $(1/N)$-squares obtained by dividing the unit square into congruent squares
with side length $1/N$ in the standard way.
Since the set $\Gamma$ is defined by the boundaries of a finite set of polygons,
there are at most $C_1N$ special $(1/N)$-squares that exhibit $\Gamma$-split in its interior or on part of its boundary.
We refer to these as $\Gamma$-split special $(1/N)$-squares.
We also say that a special $(1/N)$-square on some $yz$-parallel unit square is \textit{exceptional} if it has $\Gamma$-split
or it is in a $3\times3$ array of $9$ special $(1/N)$-squares, at least one of which has $\Gamma$-split.
Note that some of these $9$ special $(1/N)$-squares may lie on neighboring $yz$-parallel unit squares.

We start with a fixed $\Gamma$-split $(1/N)$-square $A$ on some $yz$-parallel unit square.
Assume that the bad direction vector has length between $2^{\ell-1}$ and~$2^\ell$, where the positive integer $\ell\le m$.

(i)
Clearly $A+j^*\bfv$ is contained in a $3\times3$ array of $9$ exceptional $(1/N)$-squares contained
in at most $4$ distinct but adjoining $yz$-parallel unit squares.
Each such exceptional $(1/N)$-square contributes a set of bad directions for $A$ with surface area measure at most
\textcolor{white}{xxxxxxxxxxxxxxxxxxxxxxxxxxxxxx}
\begin{displaymath}
C_1\left(\frac{1/N}{2^\ell}\right)^2
\end{displaymath}
on~$\bfS^2_+$.

(ii)
The number of distinct $yz$-parallel unit squares reachable from the center of $A$ by a vector
of length between $2^{\ell-1}$ and~$2^\ell$ is at most $C_2(2^\ell)^3$.

(iii)
The total number of exceptional $(1/N)$-squares on any $yz$-parallel unit square is at most $C_3N$.

Combining (i)--(iii), we see immediately that the total measure of bad directions for $A$
that are characterized by bad direction vectors of length between $2^{\ell-1}$ and~$2^\ell$ is at most
\textcolor{white}{xxxxxxxxxxxxxxxxxxxxxxxxxxxxxx}
\begin{displaymath}
C_1\left(\frac{1/N}{2^\ell}\right)^2
\cdot
C_2(2^\ell)^3
\cdot
C_3N
=C_4\frac{2^\ell}{N}.
\end{displaymath}
It follows, in view of \eqref{eq4.6}, that the total measure of bad directions for $A$
that are characterized by bad direction vectors of length up to $\eta N$ is at most
\begin{equation}\label{eq4.7}
\sum_{\ell=1}^mC_4\frac{2^\ell}{N}<C_5\eta.
\end{equation}

Next, note that $A$ is contained in a $3\times3$ array of $9$ exceptional $(1/N)$-squares contained
in at most $4$ distinct but adjoining $yz$-parallel unit squares.
Recall that the total number of exceptional $(1/N)$-squares on any $yz$-parallel unit square is at most~$C_3N$.
Let $\Psi(A(\bfs_0);\bfv^*;N^2)$ denote the number of $\Gamma$-split-free chains in \eqref{eq4.5} of the form
\textcolor{white}{xxxxxxxxxxxxxxxxxxxxxxxxxxxxxx}
\begin{displaymath}
A+j_1\bfv,\ldots,A+j_2\bfv,
\quad
0\le j_1<j_2\le\eta N.
\end{displaymath}
Then we have the inequality
\begin{equation}\label{eq4.8}
\int_{[0,1)^2}\int_{\bfS^2}\Psi(A(\bfs_0);\bfv^*;N^2)\,\dd\bfv^*\,\dd\bfs_0
\le C_3N\cdot C_5\eta
=C_6\eta N.
\end{equation}

\begin{remark}
An inequality of the form
\begin{displaymath}
\int_{\bfS^2}\Psi(A(\bfs_0);\bfv^*;N^2)\,\dd\bfv^*
\le C_5\eta
\end{displaymath}
for every $A(\bfs^*)$ would be ideal, but cannot be deduced from the upper bound \eqref{eq4.7},
in view of the possibility that the integrand may exceed~$1$.
Thus we need to average over $\bfs_0$ as well.
The inequality \eqref{eq4.8} is valid with the integrand $\Psi(A(\bfs^*);\bfv^*;N^2)$,
since if the sequence has more than one $\Gamma$-split-free chain,
the multiplicity is taken care of by the parameter~$\bfs_0$,
as the terms of the sequence \eqref{eq4.5} correspond to $N^2$ distinct values of~$\bfs_0$.
\end{remark}

Let $\Phi(\eps_2N;\eta N)$ denote the set of pairs $(\bfs_0,\bfv^*)\in[0,1)^2\times\bfS^2_+$
such that the sequence \eqref{eq4.5} has more than $\eps_2N$ $\Gamma$-split-free chains with length at most~$\eta N$.
Then writing $\meas$ for the product measure, we have
\begin{equation}\label{eq4.9}
\meas(\Phi(\eps_2N;\eta N))\eps_2N\le\int_{[0,1)^2}\int_{\bfS^2}\Psi(A(\bfs_0);\bfv^*;N^2)\,\dd\bfv^*\,\dd\bfs_0.
\end{equation}
Combining \eqref{eq4.8} and \eqref{eq4.9}, we conclude that
\begin{displaymath}
\meas(\Phi(\eps_2N;\eta N))\le\frac{C_6\eta}{\eps_2}.
\end{displaymath}
The proof is now complete if we choose $\eta=c_3\eps_2^2$ with a suitable constant~$c_3$.
\end{proof}

\begin{proof}[Proof of Lemma~\ref{lem36}]
For any integer $k=1,\ldots,N^2$, write
\begin{equation}\label{eq4.10}
\Omega(N;k)=\left\{\bfv_0=(\alpha_1,\alpha_2)\in[0,1]^2:\Vert k\alpha_1\Vert<\frac{2}{N}\mbox{ and }\Vert k\alpha_2\Vert<\frac{2}{N}\right\}.
\end{equation}
It is not difficult to see that
\begin{displaymath}
\Omega(N;k)\cap[0,1)^2
=\left(
\left[0,\frac{2}{kN}\right)
\cup\bigcup_{j=1}^{k-1}\left(\frac{j}{k}-\frac{2}{kN},\frac{j}{k}+\frac{2}{kN}\right)
\cup\left(1-\frac{2}{kN},1\right)
\right)^2,
\end{displaymath}
so that
\textcolor{white}{xxxxxxxxxxxxxxxxxxxxxxxxxxxxxx}
\begin{equation}\label{eq4.11}
\lambda_2(\Omega(N;k))
=\left(\frac{4}{N}\right)^2
=\frac{16}{N^2}.
\end{equation}
For every $(\alpha_1,\alpha_2)\in[0,1]^2$, consider the counting function
\begin{equation}\label{eq4.12}
\omega_N(\alpha_1,\alpha_2)=\vert\{k=1,\ldots,N^2:(\alpha_1,\alpha_2)\in\Omega(N;k)\}\vert.
\end{equation}
Combining \eqref{eq4.10}--\eqref{eq4.12}, we see that
\begin{displaymath}
\int_0^1\int_0^1\omega_N(\alpha_1,\alpha_2)\,\dd\alpha_1\,\dd\alpha_2
=\sum_{k=1}^{N^2}\lambda_2(\Omega(N;k))
=16.
\end{displaymath}

Given any $\eps_2>0$, write
\begin{displaymath}
\Omega(\eps_2)=\left\{(\alpha_1,\alpha_2)\in[0,1]^2:\omega_N(\alpha_1,\alpha_2)\ge\frac{16}{\eps_2}\right\}.
\end{displaymath}
Since the function $\omega_N(\alpha_1,\alpha_2)$ is non-negative, we clearly have
\begin{displaymath}
16=\int_0^1\int_0^1\omega_N(\alpha_1,\alpha_2)\,\dd\alpha_1\,\dd\alpha_2\ge\frac{16}{\eps_2}\lambda_2(\Omega(\eps_2)),
\end{displaymath}
so that $\lambda_2(\Omega(\eps_2))\le\eps_2$.
It follows that
\begin{displaymath}
\lambda_2([0,1]^2\setminus\Omega(\eps_2))>1-\eps_2,
\end{displaymath}
so that the collection of vectors
\begin{equation}\label{eq4.13}
\bfv_0=(\alpha_1,\alpha_2)\in[0,1]^2\setminus\Omega(\eps_2)
\end{equation}
represents at least $(1-\eps_2)$-proportion of the set $[0,1]^2$.

We shall show that the lemma holds with the choice
\begin{displaymath}
C^*=C^*(\eps_2)=1+\frac{16}{\eps_2}.
\end{displaymath}
Suppose on the contrary that some axis-parallel square $Q$ with side length $2/N$ contains modulo one
more than $1+16/\eps_2$ elements of some sequence
\begin{displaymath}
\bfs_0+j\bfv_0,
\quad
j=0,1,\ldots,N^2-1.
\end{displaymath}
In other words, suppose that there exists a subset $J\subset\{0,1,\ldots,N^2-1\}$ with $\vert J\vert>1+16/\eps_2$ such that
\textcolor{white}{xxxxxxxxxxxxxxxxxxxxxxxxxxxxxx}
\begin{displaymath}
\bfs_0+j\bfv_0\in Q,
\quad
j\in J.
\end{displaymath}
Let $j_0$ be the smallest element of~$J$, and let
\begin{displaymath}
J^*=J\setminus\{j_0\}.
\end{displaymath}
It is not difficult to see that
\begin{displaymath}
J^*\subset\{k=1,\ldots,N^2:(\alpha_1,\alpha_2)\in\Omega(N;k)\},
\end{displaymath}
and so
\textcolor{white}{xxxxxxxxxxxxxxxxxxxxxxxxxxxxxx}
\begin{displaymath}
\omega_N(\bfv_0)
\ge\vert J^*\vert
>\frac{16}{\eps_2}.
\end{displaymath}
This implies that $\bfv_0\in\Omega(\eps_2)$, contradicting \eqref{eq4.13}.
The lemma now follows.
\end{proof}

%
%

\section{Establishing Lemma~\ref{lem34}}\label{sec5}

In this section, we establish an intermediate result from which Lemma~\ref{lem34} follows easily.

\begin{lemma}\label{lem51}
Let the integer $H$ be even and positive, and let $\Delta=\Delta(u_1,u_2;\theta)$ be an arbitrary rectangle with side lengths $2u_1$ and~$2u_2$,
and tilted by an angle $\theta$ in the anticlockwise direction, where $0<u_1,u_2<1/2$.
Let
\begin{displaymath}
F(\bfs_0;\alpha_1,\alpha_2;H)
=\vert\{j=0,1,\ldots,H-1:\bfs_0+j(\alpha_1,\alpha_2)\in\Delta\}\vert.
\end{displaymath}
Then for any parameter $\kappa>1$, we have
\begin{displaymath}
\lambda_2\left(\left\{(\alpha_1,\alpha_2)\in[0,1]^2:F(\bfs_0;\alpha_1,\alpha_2;H)\ge\frac{u_1u_2H}{2}-\kappa\psi(u_1,u_2)\right\}\right)\ge1-\frac{1}{\kappa},
\end{displaymath}
where
\textcolor{white}{xxxxxxxxxxxxxxxxxxxxxxxxxxxxxx}
\begin{equation}\label{eq5.1}
\psi(u_1,u_2)=2^{23}\max\left\{\frac{u_1}{u_2^{1/2}},\frac{u_2}{u_1^{1/2}}\right\}+\frac{2^{16}}{u_1^{1/3}u_2^{1/3}}.
\end{equation}
\end{lemma}

\begin{proof}[Proof of Lemma~\ref{lem34}]
We take $H=N^2$.
For the buffer zone $\Delta=\BBB_N$, we have
\begin{displaymath}
u_1=\frac{c_1}{2},
\quad\mbox{and}\quad
u_2=\frac{1}{4N}.
\end{displaymath}
Then
\textcolor{white}{xxxxxxxxxxxxxxxxxxxxxxxxxxxxxx}
\begin{equation}\label{eq5.2}
\frac{u_1u_2H}{2}=\frac{c_1N}{16}
\quad\mbox{and}\quad
\psi(u_1,u_2)=c_4N^{1/2}+c_5N^{1/3}
\end{equation}
for some positive absolute constants $c_4$ and~$c_5$.
For large values of~$N$, the right hand side in \eqref{eq5.2} is much smaller than the left hand side.
Thus given any $\eps_2>0$, we can choose a sufficiently large $\kappa=\kappa(\eps_2)>1$
to guarantee that the inequality \eqref{eq3.2} holds for some constant $c_2$ satisfying $0<c_2<c_1/16$.
\end{proof}

\begin{proof}[Proof of Lemma~\ref{lem51}]
We proceed by a number of steps.


\begin{step1}
We aim to give a good description of the term $F(\bfs_0;\alpha_1,\alpha_2;H)$.

Let $\LLL(\alpha_1,\alpha_2)$ denote the lattice in $\Rr^3$ generated by the vectors
\begin{displaymath}
\bfe_1=(\alpha_1,\alpha_2,1),
\quad
\bfe_2=(-1,0,0),
\quad
\bfe_3=(0,-1,0),
\end{displaymath}
and consider the $3\times3$ matrix
\begin{equation}\label{eq5.3}
M=\left(\begin{array}{ccc}
\bfe_1&\bfe_2&\bfe_3
\end{array}\right)
=\left(\begin{array}{ccc}
\alpha_1&-1&0\\
\alpha_2&0&-1\\
1&0&0
\end{array}\right),
\quad\mbox{with $\det(M)=1$}.
\end{equation}
Then, with $\bfn\in\Zz^3$ expressed as column vectors, we have
\begin{equation}\label{eq5.4}
\LLL(\alpha_1,\alpha_2)=\{M\bfn:\bfn\in\Zz^3\}.
\end{equation}
Writing $B=(\Delta-\bfs_0)\times[0,H)\subset\Rr^3$, we then have
\begin{equation}\label{eq5.5}
F(\bfs_0;\alpha_1,\alpha_2;H)
=\vert\LLL(\alpha_1,\alpha_2)\cap B\vert
=\sum_{\bfn\in\Zz^3}\chi_B(M\bfn),
\end{equation}
where $\chi_B$ is the characteristic function of~$B$.

We next use the Poisson summation formula, that under some mild condition on a function $f:\Rr^3\to\Rr$, we have
\begin{equation}\label{eq5.6}
\sum_{\bfn\in\Zz^3}f(\bfn)=\sum_{\bfm\in\Zz^3}\int_{\Rr^3}f(\bfx)\ee^{-2\pi\ii\bfx\cdot\bfm}\,\dd\bfx.
\end{equation}
Applying this formula with $f=\chi_B$, it follows from \eqref{eq5.5} and noting \eqref{eq5.3} that
\begin{align}\label{eq5.7}
&
F(\bfs_0;\alpha_1,\alpha_2;H)
=\sum_{\bfn\in\Zz^3}\chi_B(M\bfn)
=\sum_{\bfm\in\Zz^3}\int_{\Rr^3}\chi_B(M\bfx)\ee^{-2\pi\ii\bfx\cdot\bfm}\,\dd\bfx
\nonumber
\\
&\quad
=\frac{1}{\det(M)}\sum_{\bfm\in\Zz^3}\int_{\Rr^3}\chi_B(\bfz)\ee^{-2\pi\ii M^{-1}\bfz\cdot\bfm}\,\dd\bfz
=\sum_{\bfm\in\Zz^3}\int_B\ee^{-2\pi\ii M^{-1}\bfz\cdot\bfm}\,\dd\bfz.
\end{align}
Note, in particular, that for $\bfm=\bzero$, we have
\begin{equation}\label{eq5.8}
\int_B\ee^{-2\pi\ii M^{-1}\bfz\cdot\bfm}\,\dd\bfz
=\int_B\ee^{-2\pi\ii M^{-1}\bfz\cdot\bzero}\,\dd\bfz
=\lambda_3(B)
=\lambda_2(\Delta)H
=4u_1u_2H.
\end{equation}
On the other hand, it is easy to see that the inverse matrix
\begin{displaymath}
M^{-1}=\left(\begin{array}{ccc}
0&0&1\\
-1&0&\alpha_1\\
0&-1&\alpha_2
\end{array}\right),
\end{displaymath}
so that $M^{-1}\bfz=(z_3,\alpha_1z_3-z_1,\alpha_2z_3-z_2)$.
Thus with $\bfv=(1,\alpha_1,\alpha_2)$, we have
\begin{align}\label{eq5.9}
M^{-1}\bfz\cdot\bfm
&
=z_3m_1+(\alpha_1z_3-z_1)m_2+(\alpha_2z_3-z_2)m_3
\nonumber
\\
&
=-z_1m_2-z_2m_3+z_3(m_1+\alpha_1m_2+\alpha_2m_3)
\nonumber
\\
&
=-z_1m_2-z_2m_3+z_3\bfv\cdot\bfm.
\end{align}
Then for every $\bfm\in\Zz^3\setminus\{\bzero\}$, we can write
\begin{equation}\label{eq5.10}
\int_B\ee^{-2\pi\ii M^{-1}\bfz\cdot\bfm}\,\dd\bfz
=\III(\bfm;\Delta-\bfs_0)\JJJ(\bfm;\bfv;H),
\end{equation}
where
\textcolor{white}{xxxxxxxxxxxxxxxxxxxxxxxxxxxxxx}
\begin{equation}\label{eq5.11}
\III(\bfm;\Delta-\bfs_0)
=\int_{\Delta-\bfs_0}\ee^{2\pi\ii(z_1m_2+z_2m_3)}\,\dd z_1\,\dd z_2
\end{equation}
and
\textcolor{white}{xxxxxxxxxxxxxxxxxxxxxxxxxxxxxx}
\begin{equation}\label{eq5.12}
\JJJ(\bfm;\bfv;H)
=\int_0^H\ee^{-2\pi\ii z_3\bfv\cdot\bfm}\,\dd z_3
=\frac{1-\ee^{-2\pi\ii H\bfv\cdot\bfm}}{2\pi\ii\bfv\cdot\bfm}.
\end{equation}

Write
\textcolor{white}{xxxxxxxxxxxxxxxxxxxxxxxxxxxxxx}
\begin{displaymath}
\bfz^*=(z_1,z_2)
\quad\mbox{and}\quad
\bfm^*=(m_2,m_3),
\end{displaymath}
and let
\textcolor{white}{xxxxxxxxxxxxxxxxxxxxxxxxxxxxxx}
\begin{displaymath}
\Delta=\rho([-u_1,u_1]\times[-u_2,u_2])+\bfw\subset[0,1)^2,
\end{displaymath}
where $\rho$ denotes an anticlockwise rotation.
Then
\begin{displaymath}
\bfz^*\in\Delta-\bfs_0=\rho([-u_1,u_1]\times[-u_2,u_2])+\bfw-\bfs_0
\end{displaymath}
if and only if
\textcolor{white}{xxxxxxxxxxxxxxxxxxxxxxxxxxxxxx}
\begin{displaymath}
\bfz'=(z'_1,z'_2)=\rho^{-1}(\bfz^*-\bfw+\bfs_0)\in[-u_1,u_1]\times[-u_2,u_2].
\end{displaymath}
Note that
\textcolor{white}{xxxxxxxxxxxxxxxxxxxxxxxxxxxxxx}
\begin{displaymath}
\bfz'=\rho^{-1}(\bfz^*-\bfw+\bfs_0)
\quad\mbox{if and only if}\quad
\bfz^*=\rho\bfz'+\bfw-\bfs_0.
\end{displaymath}
It then follows from \eqref{eq5.11} that
\begin{align}\label{eq5.13}
\III(\bfm;\Delta-\bfs_0)
&
=\int_{\Delta-\bfs_0}\ee^{2\pi\ii\bfz^*\cdot\bfm^*}\,\dd\bfz^*
=\int_{[-u_1,u_1]\times[-u_2,u_2]}\ee^{2\pi\ii(\rho\bfz'+\bfw-\bfs_0)\cdot\bfm^*}\,\dd\bfz'
\nonumber
\\
&
=\ee^{2\pi\ii(\bfw-\bfs_0)\cdot\bfm^*}\int_{[-u_1,u_1]\times[-u_2,u_2]}\ee^{2\pi\ii\rho\bfz'\cdot\bfm^*}\,\dd\bfz'
\nonumber
\\
&
=\ee^{2\pi\ii(\bfw-\bfs_0)\cdot\bfm^*}\int_{[-u_1,u_1]\times[-u_2,u_2]}\ee^{2\pi\ii\bfz'\cdot\rho^{-1}\bfm^*}\,\dd\bfz',
\end{align}
since $\rho\bfz'\cdot\bfm^*=\bfz'\cdot\rho^{-1}\bfm^*$ as the scalar product remains unchanged
under identical rotation for both constituent vectors.
Furthermore, if $\rho$ is anticlockwise rotation by an angle~$\theta$, then
\begin{equation}\label{eq5.14}
\rho^{-1}\bfm^*=\rho^{-1}(m_2,m_3)=(m_2\cos\theta+m_3\sin\theta,m_3\cos\theta-m_2\sin\theta).
\end{equation}
Combining \eqref{eq5.13} and \eqref{eq5.14}, we conclude that
\begin{equation}\label{eq5.15}
\III(\bfm;\Delta-\bfs_0)=\ee^{2\pi\ii(\bfw-\bfs_0)\cdot\bfm^*}\III(m_2,m_3;\theta;u_1)\III(m_3,-m_2;\theta;u_2),
\end{equation}
where
\textcolor{white}{xxxxxxxxxxxxxxxxxxxxxxxxxxxxxx}
\begin{align}\label{eq5.16}
\III(m_2,m_3;\theta;u_1)
&
=\int_{-u_1}^{u_1}\ee^{2\pi\ii z'_1(m_2\cos\theta+m_3\sin\theta)}\,\dd z'_1
\nonumber
\\
&
=2\int_0^{u_1}\cos(2\pi z'_1(m_2\cos\theta+m_3\sin\theta))\,\dd z'_1
\nonumber
\\
&
=\frac{\sin(2\pi u_1(m_2\cos\theta+m_3\sin\theta))}{\pi(m_2\cos\theta+m_3\sin\theta)}
\end{align}
and
\textcolor{white}{xxxxxxxxxxxxxxxxxxxxxxxxxxxxxx}
\begin{align}\label{eq5.17}
\III(m_3,-m_2;\theta;u_2)
&
=\int_{-u_2}^{u_2}\ee^{2\pi\ii z'_2(m_3\cos\theta-m_2\sin\theta)}\,\dd z'_2
\nonumber
\\
&
=2\int_0^{u_2}\cos(2\pi z'_2(m_3\cos\theta-m_2\sin\theta))\,\dd z'_2
\nonumber
\\
&
=\frac{\sin(2\pi u_2(m_3\cos\theta-m_2\sin\theta))}{\pi(m_3\cos\theta-m_2\sin\theta)}.
\end{align}
Finally, combining \eqref{eq5.7}, \eqref{eq5.8}, \eqref{eq5.10} and \eqref{eq5.15}, we conclude that
\begin{align}\label{eq5.18}
&
F(\bfs_0;\alpha_1,\alpha_2;H)-4u_1u_2H
\nonumber
\\
&\quad
=\sum_{\bfm\in\Zz^3\setminus\{\bzero\}}\ee^{2\pi\ii(\bfw-\bfs_0)\cdot\bfm^*}
\III(m_2,m_3;\theta;u_1)\III(m_3,-m_2;\theta;u_2)\JJJ(\bfm;\bfv;H),
\end{align}
where the details for the various factors in the summand are given by \eqref{eq5.12}, \eqref{eq5.16} and \eqref{eq5.17}.
\end{step1}


\begin{step2}
Here we contract the interval $[0,H]$ in the third direction and average over all contractions.
More precisely, for every $h\in\Rr$ satisfying $0\le h\le H/2$, consider the smaller set
\textcolor{white}{xxxxxxxxxxxxxxxxxxxxxxxxxxxxxx}
\begin{displaymath}
B(h)=(\Delta-\bfs_0)\times[h,H-h)\subset B.
\end{displaymath}
Clearly, in view of an analog of \eqref{eq5.8}, we have
\begin{align}\label{eq5.19}
&
F(\bfs_0;\alpha_1,\alpha_2;H)
\ge\frac{2}{H}\int_0^{H/2}\left(\sum_{\bfn\in\Zz^3}\chi_{B(h)}(M\bfn)\right)\dd h
\nonumber
\\
&\quad
=\frac{2}{H}\int_0^{H/2}4u_1u_2(H-2h)\,\dd h
+\frac{2}{H}\int_0^{H/2}\left(\sum_{\bfm\in\Zz^3\setminus\{\bzero\}}\int_{B(h)}\ee^{-2\pi\ii M^{-1}\bfz\cdot\bfm}\,\dd\bfz\right)\dd h
\nonumber
\\
&\quad
=2u_1u_2H
+\frac{2}{H}\int_0^{H/2}\left(\sum_{\bfm\in\Zz^3\setminus\{\bzero\}}\int_{B(h)}\ee^{-2\pi\ii M^{-1}\bfz\cdot\bfm}\,\dd\bfz\right)\dd h.
\end{align}
For any $\bfm\in\Zz^3\setminus\{\bzero\}$, analogous to \eqref{eq5.10}, we have
\begin{equation}\label{eq5.20}
\int_{B(h)}\ee^{-2\pi\ii M^{-1}\bfz\cdot\bfm}\,\dd\bfz
=\III(\bfm;\Delta-\bfs_0)\JJJ(\bfm;\bfv;H;h)
\end{equation}
where the first two directions are unaffected and
\begin{equation}\label{eq5.21}
\JJJ(\bfm;\bfv;H;h)
=\int_h^{H-h}\ee^{-2\pi\ii z_3\bfv\cdot\bfm}\,\dd z_3
=\frac{\ee^{-2\pi\ii h\bfv\cdot\bfm}-\ee^{-2\pi\ii(H-h)\bfv\cdot\bfm}}{2\pi\ii\bfv\cdot\bfm},
\end{equation}
so that
\begin{displaymath}
\int_0^{H/2}\JJJ(\bfm;\bfv;H;h)\,\dd h
=\frac{1-2\ee^{-\pi\ii H\bfv\cdot\bfm}+\ee^{-2\pi\ii H\bfv\cdot\bfm}}{(2\pi\ii\bfv\cdot\bfm)^2}
=-\frac{(1-\ee^{-\pi\ii H\bfv\cdot\bfm})^2}{4\pi^2(\bfv\cdot\bfm)^2},
\end{displaymath}
and so
\textcolor{white}{xxxxxxxxxxxxxxxxxxxxxxxxxxxxxx}
\begin{equation}\label{eq5.22}
\widetilde{\JJJ}(\bfm;\bfv;H)
=\frac{2}{H}\int_0^{H/2}\JJJ(\bfm;\bfv;H;h)\,\dd h
=-\frac{(1-\ee^{-\pi\ii H\bfv\cdot\bfm})^2}{2\pi^2H(\bfv\cdot\bfm)^2}.
\end{equation}
Finally, combining \eqref{eq5.15}, \eqref{eq5.19}, \eqref{eq5.20} and \eqref{eq5.22}, we conclude that
\begin{align}\label{eq5.23}
&
F(\bfs_0;\alpha_1,\alpha_2;H)-2u_1u_2H
\nonumber
\\
&\quad
\ge\sum_{\bfm\in\Zz^3\setminus\{\bzero\}}\ee^{2\pi\ii(\bfw-\bfs_0)\cdot\bfm^*}
\III(m_2,m_3;\theta;u_1)\III(m_3,-m_2;\theta;u_2)\widetilde{\JJJ}(\bfm;\bfv;H),
\end{align}
the analog of \eqref{eq5.18}.
\end{step2}


\begin{step3}
Here we also contract the rectangle $[-u_1,u_1]\times[-u_2,u_2]$ in the first two directions and average over all contractions.
More precisely, for every $\gamma_1$ and $\gamma_2$ satisfying $0\le\gamma_1,\gamma_2\le1$, consider the smaller set
\begin{displaymath}
\Delta(\gamma_1,\gamma_2)=\rho([-\gamma_1u_1,\gamma_1u_1]\times[-\gamma_2u_2,\gamma_2u_2])+\bfw\subset[0,1)^2,
\end{displaymath}
and the corresponding smaller set
\begin{displaymath}
B(\gamma_1,\gamma_2;h)=(\Delta(\gamma_1,\gamma_2)-\bfs_0)\times[h,H-h)\subset B(h)\subset B.
\end{displaymath}
Clearly, in view of an analog of \eqref{eq5.8}, we have
\begin{align}\label{eq5.24}
&
F(\bfs_0;\alpha_1,\alpha_2;H)
\ge\frac{2}{H}\int_0^{H/2}\int_0^1\int_0^1\left(\sum_{\bfn\in\Zz^3}\chi_{B(\gamma_1,\gamma_2,h)}(M\bfn)\right)\dd\gamma_1\,\dd\gamma_2\,\dd h
\nonumber
\\
&\quad
=\frac{2}{H}\int_0^{H/2}\int_0^1\int_0^14\gamma_1\gamma_2u_1u_2(H-2h)\,\dd\gamma_1\,\dd\gamma_2\,\dd h
\nonumber
\\
&\quad\qquad
+\frac{2}{H}\int_0^{H/2}\int_0^1\int_0^1
\left(\sum_{\bfm\in\Zz^3\setminus\{\bzero\}}\int_{B(\gamma_1,\gamma_2,h)}\ee^{-2\pi\ii M^{-1}\bfz\cdot\bfm}\,\dd\bfz\right)
\,\dd\gamma_1\,\dd\gamma_2\,\dd h,
\end{align}
where the first term
\begin{equation}\label{eq5.25}
\frac{2}{H}\int_0^{H/2}\int_0^1\int_0^14\gamma_1\gamma_2u_1u_2(H-2h)\,\dd\gamma_1\,\dd\gamma_2\,\dd h=\frac{u_1u_2H}{2}.
\end{equation}
For any $\bfm\in\Zz^3\setminus\{\bzero\}$, analogous to \eqref{eq5.10} and \eqref{eq5.20}, we have
\begin{equation}\label{eq5.26}
\int_{B(\gamma_1,\gamma_2,h)}\ee^{-2\pi\ii M^{-1}\bfz\cdot\bfm}\,\dd\bfz
=\III(\bfm;\Delta(\gamma_1,\gamma_2)-\bfs_0)\JJJ(\bfm;\bfv;H;h)
\end{equation}
Then, analogous to \eqref{eq5.15}, we have
\begin{equation}\label{eq5.27}
\III(\bfm;\Delta(\gamma_1,\gamma_2)-\bfs_0)=\ee^{2\pi\ii(\bfw-\bfs_0)\cdot\bfm^*}\III(m_2,m_3;\theta;\gamma_1u_1)\III(m_3,-m_2;\theta;\gamma_2u_2).
\end{equation}
Simple calculations now give
\begin{align}\label{eq5.28}
\widetilde{\III}(m_2,m_3;\theta;u_1)
&
=\int_0^1\III(m_2,m_3;\theta;\gamma_1u_1)\,\dd\gamma_1
\nonumber
\\
&
=\int_0^1\frac{\sin(2\pi\gamma_1u_1(m_2\cos\theta+m_3\sin\theta))}{\pi(m_2\cos\theta+m_3\sin\theta)}\,\dd\gamma_1
\nonumber
\\
&
=\frac{1-\cos(2\pi u_1(m_2\cos\theta+m_3\sin\theta))}{2\pi^2u_1(m_2\cos\theta+m_3\sin\theta)^2}
\nonumber
\\
&
=\frac{\sin^2(\pi u_1(m_2\cos\theta+m_3\sin\theta))}{\pi^2u_1(m_2\cos\theta+m_3\sin\theta)^2}
\end{align}
and
\textcolor{white}{xxxxxxxxxxxxxxxxxxxxxxxxxxxxxx}
\begin{align}\label{eq5.29}
\widetilde{\III}(m_3,-m_2;\theta;u_2)
&
=\int_0^1\III(m_3,-m_2;\theta;\gamma_2u_2)\,\dd\gamma_2
\nonumber
\\
&
=\int_0^1\frac{\sin(2\pi\gamma_2u_2(m_3\cos\theta-m_2\sin\theta))}{\pi(m_3\cos\theta-m_2\sin\theta)}\,\dd\gamma_2
\nonumber
\\
&
=\frac{1-\cos(2\pi u_2(m_3\cos\theta-m_2\sin\theta))}{2\pi^2u_2(m_3\cos\theta-m_2\sin\theta)^2}
\nonumber
\\
&
=\frac{\sin^2(\pi u_2(m_3\cos\theta-m_2\sin\theta))}{\pi^2u_2(m_3\cos\theta-m_2\sin\theta)^2}.
\end{align}
It now follows from \eqref{eq5.22} and \eqref{eq5.24}--\eqref{eq5.29} that
\begin{align}\label{eq5.30}
&
F(\bfs_0;\alpha_1,\alpha_2;H)-\frac{u_1u_2H}{2}
\nonumber
\\
&\quad
\ge\sum_{\bfm\in\Zz^3\setminus\{\bzero\}}\ee^{2\pi\ii(\bfw-\bfs_0)\cdot\bfm^*}
\widetilde{\III}(m_2,m_3;\theta;u_1)\widetilde{\III}(m_3,-m_2;\theta;u_2)\widetilde{\JJJ}(\bfm;\bfv;H),
\end{align}
the analog of \eqref{eq5.18} and \eqref{eq5.23}.
\end{step3}


\begin{step4}
It follows from \eqref{eq5.30} that
\begin{displaymath}
F(\bfs_0;\alpha_1,\alpha_2;H)\ge\frac{u_1u_2H}{2}+\sum_{\bfm\in\Zz^3\setminus\{\bzero\}}\Lambda(\theta;u_1,u_2;H;\bfv;\bfm),
\end{displaymath}
where
\begin{equation}\label{eq5.31}
\Lambda(\theta;u_1,u_2;H;\bfv;\bfm)
=\widetilde{\III}(m_2,m_3;\theta;u_1)\widetilde{\III}(m_3,-m_2;\theta;u_2)\vert\widetilde{\JJJ}(\bfm;\bfv;H)\vert,
\end{equation}
since it is clear from \eqref{eq5.28} and \eqref{eq5.29} that $\widetilde{\III}(m_2,m_3;\theta;u_1)$ and $\widetilde{\III}(m_3,-m_2;\theta;u_2)$
are real and non-negative.

The inequalities in Steps 1--3 are at this stage only formal inequalities, as we have not considered the question of convergence.
Write
\begin{equation}\label{eq5.32}
\Xi(\alpha_1,\alpha_2)=\sum_{\bfm\in\Zz^3\setminus\{\bzero\}}\Lambda(\theta;u_1,u_2;H;\bfv;\bfm),
\end{equation}
where $\bfv=(1,\alpha_1,\alpha_2)$.
We shall use the first-moment method and analyze the average
\begin{displaymath}
\int_0^1\int_0^1\Xi(\alpha_1,\alpha_2)\,\dd\alpha_1\,\dd\alpha_2,
\end{displaymath}
and remove those directions $(\alpha_1,\alpha_2)$ for which $\Xi(\alpha_1,\alpha_2)$ is substantially larger than the average.
In this step, we do some further preparation.

First of all, the inequalities
\begin{equation}\label{eq5.33}
\widetilde{\III}(m_2,m_3;\theta;u_1)
\le\min\left\{u_1,\frac{1}{\pi^2u_1(m_2\cos\theta+m_3\sin\theta)^2}\right\}
\end{equation}
and
\textcolor{white}{xxxxxxxxxxxxxxxxxxxxxxxxxxxxxx}
\begin{equation}\label{eq5.34}
\widetilde{\III}(m_3,-m_2;\theta;u_2)
\le\min\left\{u_2,\frac{1}{\pi^2u_2(m_3\cos\theta-m_2\sin\theta)^2}\right\}
\end{equation}
follow on applying the inequality $\vert\sin y\vert\le\min\{\vert y\vert,1\}$, which holds for every $y\in\Rr$,
to \eqref{eq5.28} and \eqref{eq5.29} respectively.
Next, the inequality
\begin{equation}\label{eq5.35}
\vert\widetilde{\JJJ}(\bfm;\bfv;H)\vert
=\frac{\vert1-\ee^{-\pi\ii H\bfv\cdot\bfm}\vert^2}{2\pi^2H(\bfv\cdot\bfm)^2}
\le\min\left\{\frac{H}{2},\frac{2}{\pi^2H(\bfv\cdot\bfm)^2}\right\}
\end{equation}
follows on applying the inequality $\vert1-\ee^{\ii y}\vert\le\min\{\vert y\vert,2\}$, which also holds for every $y\in\Rr$, to \eqref{eq5.22}.
Since we consider only those $\bfm\ne\bzero$, it follows that in the trivial case when $\bfm^*=(m_2,m_3)=(0,0)$, we must have $m_1\ne0$
and, since $H$ is even, also
\begin{equation}\label{eq5.36}
\vert\widetilde{\JJJ}(\bfm;\bfv;H)\vert=\frac{\vert1-\ee^{-\pi\ii H\bfv\cdot\bfm}\vert^2}{2\pi^2H(\bfv\cdot\bfm)^2}=0.
\end{equation}
Hence we can assume that $\bfm^*=(m_2,m_3)\in\Zz^2\setminus\{(0,0)\}$.
We can write
\begin{equation}\label{eq5.37}
\Zz^2\setminus\{(0,0)\}=\bigcup_{j=0}^\infty\ZZZ_j
\end{equation}
as a disjoint union of subsets, where for every $j=0,1,2,3,\ldots,$ the subset
\begin{align}\label{eq5.38}
\ZZZ_j
&
=\{\bfm^*=(m_2,m_3)\in\Zz^2\setminus\{(0,0)\}:2^{j-1}<\max\{\vert m_2\vert,\vert m_3\vert\}\le2^j\}
\nonumber
\\
&
=\Zz^2\cap([-2^j,2^j]^2\setminus[-2^{j-1},2^{j-1}]^2).
\end{align}

For every $\bfm^*=(m_2,m_3)\in\ZZZ_j$, let
\begin{equation}\label{eq5.39}
\Omega_j(\bfm^*;0)=\left\{(\alpha_1,\alpha_2)\in[0,1]^2:\vert\bfv\cdot\bfm\vert\le\frac{1}{H}\mbox{ for some $m_1\in\Zz$}\right\},
\end{equation}
and for every integer $\ell\ge1$, let
\begin{equation}\label{eq5.40}
\Omega_j(\bfm^*;\ell)=\left\{(\alpha_1,\alpha_2)\in[0,1]^2:\frac{2^{\ell-1}}{H}<\vert\bfv\cdot\bfm\vert\le\frac{2^\ell}{H}\mbox{ for some $m_1\in\Zz$}\right\}.
\end{equation}
We claim that the inequality
\begin{equation}\label{eq5.41}
\lambda_2(\Omega_j(\bfm^*;\ell))
\le\min\left\{\frac{2^{\ell+5}}{H}+\frac{2^{2\ell+3-j}}{H^2},1\right\}
\le\min\left\{\frac{2^{\ell+6}}{H},1\right\}
\end{equation}
holds for every $\bfm^*=(m_2,m_3)\in\ZZZ_j$ and every integer $\ell\ge0$.

Suppose first of all that $\ell=0$.
To estimate $\lambda_2(\Omega_j(\bfm^*;0))$, we may assume that $\vert m_2\vert\le\vert m_3\vert$.
Then for fixed $m_1,\alpha_1,m_2$, the variable $\alpha_2\in[0,1]$ must satisfy the inequality
\textcolor{white}{xxxxxxxxxxxxxxxxxxxxxxxxxxxxxx}
\begin{equation}\label{eq5.42}
-\frac{1}{H}\le m_1+\alpha_1m_2+\alpha_2m_3\le\frac{1}{H},
\end{equation}
and so falls into an interval of length
\begin{displaymath}
\frac{2}{Hm_3}\le\frac{2^{2-j}}{H}.
\end{displaymath}
Note next that the inequality \eqref{eq5.42} has no solution if
\begin{displaymath}
\vert m_1\vert>2^{j+1}+\frac{1}{H}.
\end{displaymath}
It follows that
\begin{displaymath}
\lambda_2(\Omega_j(\bfm^*;0))
\le\min\left\{\left(2^{j+2}+1+\frac{2}{H}\right)\frac{2^{2-j}}{H},1\right\}
\le\min\left\{\frac{2^5}{H}+\frac{2^{3-j}}{H^2},1\right\},
\end{displaymath}
and this establishes \eqref{eq5.41} for $\ell=0$.

Suppose next that $\ell\ge1$.
To estimate $\lambda_2(\Omega_j(\bfm^*;\ell))$, we may assume again that $\vert m_2\vert\le\vert m_3\vert$.
Then for fixed $m_1,\alpha_1,m_2$, the variable $\alpha_2\in[0,1]$ must satisfy one of the inequalities
\begin{equation}\label{eq5.43}
-\frac{2^\ell}{H}\le m_1+\alpha_1m_2+\alpha_2m_3<-\frac{2^{\ell-1}}{H}
\quad\mbox{or}\quad
\frac{2^{\ell-1}}{H}<m_1+\alpha_1m_2+\alpha_2m_3\le\frac{2^\ell}{H},
\end{equation}
and so falls into two intervals of total length
\begin{displaymath}
\frac{2^\ell}{Hm_3}\le\frac{2^{\ell+1-j}}{H}.
\end{displaymath}
Note next that the inequalities in \eqref{eq5.43} have no solution if
\begin{displaymath}
\vert m_1\vert>2^{j+1}+\frac{2^\ell}{H}.
\end{displaymath}
It follows that
\begin{displaymath}
\lambda_2(\Omega_j(\bfm^*;\ell))
\le\min\left\{\left(2^{j+2}+1+\frac{2^{\ell+1}}{H}\right)\frac{2^{\ell+1-j}}{H},1\right\}
\le\min\left\{\frac{2^{\ell+4}}{H}+\frac{2^{2\ell+2-j}}{H^2},1\right\},
\end{displaymath}
and this more than establishes \eqref{eq5.41} for every $\ell\ge1$.
\end{step4}


\begin{step5}
For every $(m_2,m_3)\in\ZZZ_j$, we use the identity
\begin{equation}\label{eq5.44}
(m_2\cos\theta+m_3\sin\theta)^2+(m_3\cos\theta-m_2\sin\theta)^2=m_2^2+m_3^2
\end{equation}
to study the product $\widetilde{\III}(m_2,m_3;\theta;u_1)\widetilde{\III}(m_3,-m_2;\theta;u_2)$.

Case~1: Assume that
\begin{equation}\label{eq5.45}
\min\{(m_2\cos\theta+m_3\sin\theta)^2,(m_3\cos\theta-m_2\sin\theta)^2\}<2^{j-3}.
\end{equation}
Then
\textcolor{white}{xxxxxxxxxxxxxxxxxxxxxxxxxxxxxx}
\begin{align}\label{eq5.46}
&
\max\{(m_2\cos\theta+m_3\sin\theta)^2,(m_3\cos\theta-m_2\sin\theta)^2\}
\nonumber
\\
&\quad
=(m_2^2+m_3^2)-\min\{(m_2\cos\theta+m_3\sin\theta)^2,(m_3\cos\theta-m_2\sin\theta)^2\}
\nonumber
\\
&\quad
>2^{2j-2}-2^{j-3}>2^{2j-3}.
\end{align}
Since $\bfm^*=(m_2,m_3)\in\ZZZ_j$, it follows from \eqref{eq5.38} that the vector $(m_2,m_3)$ has length greater than~$2^{j-1}$,
while \eqref{eq5.45} shows that one of the two dot products
\begin{displaymath}
(m_2,m_3)\cdot(\cos\theta,\sin\theta)
\quad\mbox{and}\quad
(m_2,m_3)\cdot(-\sin\theta,\cos\theta)
\end{displaymath}
has size less than $2^{(j-3)/2}$ which is substantially smaller than $2^{j-1}$ if $j$ is large, meaning that the vector
$(m_2,m_3)$ is nearly perpendicular to one of the two vectors $(\cos\theta,\sin\theta)$ and $(-\sin\theta,\cos\theta)$
and nearly parallel to the other one, where nearly parallel means that the angle in between is in the range~$2^{-j/2}$, so
\begin{displaymath}
\ZZZ_j^\dagger=\{(m_2,m_3)\in\ZZZ_j:\mbox{\eqref{eq5.45} holds}\}
\end{displaymath}
is the set of points in $\Zz^2$ which fall into the dark shaded part in Figure~5.1,
and we have generous upper bound
\begin{equation}\label{eq5.47}
\vert\ZZZ_j^\dagger\vert\le2^{10}2^{2j}2^{-j/2}=2^{10+3j/2}.
\end{equation}
\begin{displaymath}
\begin{array}{c}
\includegraphics[scale=0.8]{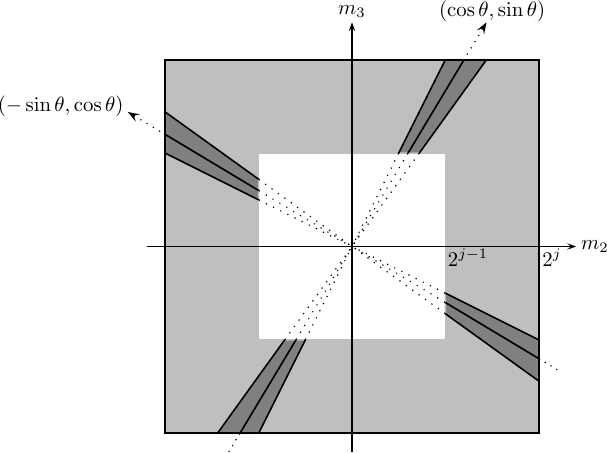}
\\
\mbox{Figure 5.1: the location of $\bfm^*=(m_2,m_3)\in\ZZZ_j^\dagger$}
\end{array}
\end{displaymath}

For every $\bfm^*=(m_2,m_3)\in\ZZZ_j^\dagger$, it then follows from the bounds \eqref{eq5.33}, \eqref{eq5.34} and \eqref{eq5.46} that
\begin{align}\label{eq5.48}
&
\widetilde{\III}(m_2,m_3;\theta;u_1)\widetilde{\III}(m_3,-m_2;\theta;u_2)\nonumber
\\
&\quad
\le\min\left\{u_1u_2,\frac{u_1}{\pi^2u_2(m_3\cos\theta-m_2\sin\theta)^2},\frac{u_2}{\pi^2u_1(m_2\cos\theta+m_3\sin\theta)^2}\right\}
\nonumber
\\
&\quad
\le\min\left\{u_1u_2,\max\left\{\frac{u_1}{u_2},\frac{u_2}{u_1}\right\}\frac{1}{\pi^22^{2j-3}}\right\}
\nonumber
\\
&\quad
\le\min\left\{u_1u_2,\frac{1}{2^{2j}}\max\left\{\frac{u_1}{u_2},\frac{u_2}{u_1}\right\}\right\}.
\end{align}

Case~2: Suppose that \eqref{eq5.45} fails, so that $\bfm^*=(m_2,m_3)\in\ZZZ_j^\ddagger=\ZZZ_j\setminus\ZZZ_j^\dagger$.
Then
\begin{equation}\label{eq5.49}
\min\{(m_2\cos\theta+m_3\sin\theta)^2,(m_3\cos\theta-m_2\sin\theta)^2\}\ge2^{j-3}.
\end{equation}
On the other hand, using the identity \eqref{eq5.44} and the inequality \eqref{eq5.38}, we see that
\begin{equation}\label{eq5.50}
\max\{(m_2\cos\theta+m_3\sin\theta)^2,(m_3\cos\theta-m_2\sin\theta)^2\}
\ge\frac{m_2^2+m_3^2}{2}
\ge2^{2j-3}.
\end{equation}
It then follows from the bounds \eqref{eq5.33}, \eqref{eq5.34}, \eqref{eq5.49} and \eqref{eq5.50} that
\begin{align}\label{eq5.51}
&
\widetilde{\III}(m_2,m_3;\theta;u_1)\widetilde{\III}(m_3,-m_2;\theta;u_2)\nonumber
\\
&\quad
\le\min\left\{u_1u_2,\frac{1}{\pi^4u_1u_2(m_2\cos\theta+m_3\sin\theta)^2(m_3\cos\theta-m_2\sin\theta)^2}\right\}
\nonumber
\\
&\quad
\le\min\left\{u_1u_2,\frac{1}{\pi^4u_1u_22^{3j-6}}\right\}
\le\min\left\{u_1u_2,\frac{1}{2^{3j}u_1u_2}\right\}.
\end{align}
\end{step5}


\begin{step6}
Using \eqref{eq5.32} and taking into account that the identity \eqref{eq5.36} holds whenever $\bfm^*=(m_2,m_3)=(0,0)$,
we now have the trivial upper bound
\begin{equation}\label{eq5.52}
\int_0^1\int_0^1\Xi(\alpha_1,\alpha_2)\,\dd\alpha_1\,\dd\alpha_2\le\frakI_1+\frakI_2,
\end{equation}
where
\textcolor{white}{xxxxxxxxxxxxxxxxxxxxxxxxxxxxxx}
\begin{align}
\frakI_1
&
=\sum_{j=0}^\infty\sum_{\ell=0}^\infty\sum_{\bfm^*\in\ZZZ_j^\dagger}\int_{\Omega_j(\bfm^*;\ell)}
\Lambda(\theta;u_1,u_2;H;\bfv;\bfm)\,\dd\alpha_1\,\dd\alpha_2,
\label{eq5.53}
\\
\frakI_2
&
=\sum_{j=0}^\infty\sum_{\ell=0}^\infty\sum_{\bfm^*\in\ZZZ_j^\ddagger}\int_{\Omega_j(\bfm^*;\ell)}
\Lambda(\theta;u_1,u_2;H;\bfv;\bfm)\,\dd\alpha_1\,\dd\alpha_2.
\label{eq5.54}
\end{align}

\begin{lemma}\label{lem52}
We have
\textcolor{white}{xxxxxxxxxxxxxxxxxxxxxxxxxxxxxx}
\begin{equation}\label{eq5.55}
\frakI_1\le2^{23}\max\left\{\frac{u_1}{u_2^{1/2}},\frac{u_2}{u_1^{1/2}}\right\}.
\end{equation}
\end{lemma}

\begin{proof}
(i)
Suppose first that $\ell=0$.
It follows from \eqref{eq5.31}, \eqref{eq5.35} and \eqref{eq5.48} that for every $\bfm^*=(m_2,m_3)\in\ZZZ_j^\dagger$
and $(\alpha_1,\alpha_2)\in\Omega_j(\bfm^*;0)$, we have
\begin{displaymath}
0\le\Lambda(\theta;u_1,u_2;H;\bfv;\bfm)\le\frac{H}{2}\min\left\{u_1u_2,\frac{1}{2^{2j}}\max\left\{\frac{u_1}{u_2},\frac{u_2}{u_1}\right\}\right\}.
\end{displaymath}
Combining this with \eqref{eq5.41} and \eqref{eq5.47}, we deduce that
\begin{align}\label{eq5.56}
&
\sum_{\bfm^*\in\ZZZ_j^\dagger}\int_{\Omega_j(\bfm^*;0)}\Lambda(\theta;u_1,u_2;H;\bfv;\bfm)\,\dd\alpha_1\,\dd\alpha_2
\nonumber
\\
&\quad
\le2^{10+3j/2}\lambda_2(\Omega_j(\bfm^*;0))\frac{H}{2}\min\left\{u_1u_2,\frac{1}{2^{2j}}\max\left\{\frac{u_1}{u_2},\frac{u_2}{u_1}\right\}\right\}
\nonumber
\\
&\quad
\le2^{15+3j/2}\min\left\{u_1u_2,\frac{1}{2^{2j}}\max\left\{\frac{u_1}{u_2},\frac{u_2}{u_1}\right\}\right\}.
\end{align}
We claim that
\begin{equation}\label{eq5.57}
\sum_{j=0}^\infty2^{3j/2}\min\left\{u_1u_2,\frac{1}{2^{2j}}\max\left\{\frac{u_1}{u_2},\frac{u_2}{u_1}\right\}\right\}
\le2^{3}\max\left\{\frac{u_1}{u_2^{1/2}},\frac{u_2}{u_1^{1/2}}\right\}.
\end{equation}
To see this, let $J^\dagger$ denote the largest non-negative integer $j$ such that
\begin{equation}\label{eq5.58}
u_1u_2\le\frac{1}{2^{2j}}\max\left\{\frac{u_1}{u_2},\frac{u_2}{u_1}\right\},
\end{equation}
so that
\textcolor{white}{xxxxxxxxxxxxxxxxxxxxxxxxxxxxxx}
\begin{equation}\label{eq5.59}
2^{2J^\dagger}\le\frac{1}{u_1u_2}\max\left\{\frac{u_1}{u_2},\frac{u_2}{u_1}\right\}
\quad\mbox{and}\quad
2^{2J^\dagger+2}>\frac{1}{u_1u_2}\max\left\{\frac{u_1}{u_2},\frac{u_2}{u_1}\right\}.
\end{equation}
Then it follows from \eqref{eq5.58} and \eqref{eq5.59} that
\begin{align}\label{eq5.60}
&
\sum_{j=0}^{J^\dagger}2^{3j/2}\min\left\{u_1u_2,\frac{1}{2^{2j}}\max\left\{\frac{u_1}{u_2},\frac{u_2}{u_1}\right\}\right\}
\le u_1u_2\sum_{j=0}^{J^\dagger}2^{3j/2}
\nonumber
\\
&\quad
\le2u_1u_2\left(\frac{1}{u_1u_2}\max\left\{\frac{u_1}{u_2},\frac{u_2}{u_1}\right\}\right)^{3/4}
=2\max\left\{\frac{u_1}{u_2^{1/2}},\frac{u_2}{u_1^{1/2}}\right\}
\end{align}
and
\begin{align}\label{eq5.61}
&
\sum_{j=J^\dagger+1}^\infty2^{3j/2}\min\left\{u_1u_2,\frac{1}{2^{2j}}\max\left\{\frac{u_1}{u_2},\frac{u_2}{u_1}\right\}\right\}
\le\max\left\{\frac{u_1}{u_2},\frac{u_2}{u_1}\right\}\sum_{j=J^\dagger+1}^\infty2^{-j/2}
\nonumber
\\
&\quad
\le4\max\left\{\frac{u_1}{u_2},\frac{u_2}{u_1}\right\}\left(\frac{1}{u_1u_2}\max\left\{\frac{u_1}{u_2},\frac{u_2}{u_1}\right\}\right)^{-1/4}
=4\max\left\{\frac{u_1}{u_2^{1/2}},\frac{u_2}{u_1^{1/2}}\right\}.
\end{align}
The inequality \eqref{eq5.57} follows on combining \eqref{eq5.60} and \eqref{eq5.61}.
Combining \eqref{eq5.56} and \eqref{eq5.57} now leads to the inequality
\begin{equation}\label{eq5.62}
\sum_{j=0}^\infty\sum_{\bfm^*\in\ZZZ_j^\dagger}\int_{\Omega_j(\bfm^*;0)}\Lambda(\theta;u_1,u_2;H;\bfv;\bfm)\,\dd\alpha_1\,\dd\alpha_2
\le2^{18}\max\left\{\frac{u_1}{u_2^{1/2}},\frac{u_2}{u_1^{1/2}}\right\}.
\end{equation}

(ii)
Suppose next that $\ell\ge1$.
It follows from \eqref{eq5.31}, \eqref{eq5.35}, \eqref{eq5.40} and \eqref{eq5.48} that for every $\bfm^*=(m_2,m_3)\in\ZZZ_j^\dagger$
and $(\alpha_1,\alpha_2)\in\Omega_j(\bfm^*;\ell)$, we have
\begin{equation}\label{eq5.63}
0\le\Lambda(\theta;u_1,u_2;H;\bfv;\bfm)\le\frac{H}{2^{2\ell}}\min\left\{u_1u_2,\frac{1}{2^{2j}}\max\left\{\frac{u_1}{u_2},\frac{u_2}{u_1}\right\}\right\}.
\end{equation}
Combining \eqref{eq5.41}, \eqref{eq5.47} and \eqref{eq5.63}, we deduce that
\begin{align}
&
\sum_{\bfm^*\in\ZZZ_j^\dagger}\int_{\Omega_j(\bfm^*;\ell)}\Lambda(\theta;u_1,u_2;H;\bfv;\bfm)\,\dd\alpha_1\,\dd\alpha_2
\nonumber
\\
&\quad
\le2^{10+3j/2}\lambda_2(\Omega_j(\bfm^*;\ell))\frac{H}{2^{2\ell}}\min\left\{u_1u_2,\frac{1}{2^{2j}}\max\left\{\frac{u_1}{u_2},\frac{u_2}{u_1}\right\}\right\}
\nonumber
\\
&\quad
\le2^{10+3j/2}\min\left\{\frac{1}{2^{\ell-8}},\frac{H}{2^{2\ell}}\right\}
\min\left\{u_1u_2,\frac{1}{2^{2j}}\max\left\{\frac{u_1}{u_2},\frac{u_2}{u_1}\right\}\right\},
\nonumber
\end{align}
so that
\textcolor{white}{xxxxxxxxxxxxxxxxxxxxxxxxxxxxxx}
\begin{align}\label{eq5.64}
&
\sum_{\ell=1}^\infty\sum_{\bfm^*\in\ZZZ_j^\dagger}\int_{\Omega_j(\bfm^*;\ell)}\Lambda(\theta;u_1,u_2;H;\bfv;\bfm)\,\dd\alpha_1\,\dd\alpha_2
\nonumber
\\
&\quad
\le2^{10+3j/2}\min\left\{u_1u_2,\frac{1}{2^{2j}}\max\left\{\frac{u_1}{u_2},\frac{u_2}{u_1}\right\}\right\}
\sum_{\ell=1}^\infty\min\left\{\frac{1}{2^{\ell-8}},\frac{H}{2^{2\ell}}\right\}
\nonumber
\\
&\quad
\le2^{19+3j/2}\min\left\{u_1u_2,\frac{1}{2^{2j}}\max\left\{\frac{u_1}{u_2},\frac{u_2}{u_1}\right\}\right\}.
\end{align}
Combining \eqref{eq5.57} and \eqref{eq5.64} now leads to the inequality
\begin{equation}\label{eq5.65}
\sum_{j=0}^\infty\sum_{\ell=1}^\infty\sum_{\bfm^*\in\ZZZ_j^\dagger}\int_{\Omega_j(\bfm^*;\ell)}\Lambda(\theta;u_1,u_2;H;\bfv;\bfm)\,\dd\alpha_1\,\dd\alpha_2
\le2^{22}\max\left\{\frac{u_1}{u_2^{1/2}},\frac{u_2}{u_1^{1/2}}\right\}.
\end{equation}

The inequality \eqref{eq5.55} now follows on combining \eqref{eq5.53}, \eqref{eq5.62} and \eqref{eq5.65}.
\end{proof}

\begin{lemma}\label{lem53}
We have
\textcolor{white}{xxxxxxxxxxxxxxxxxxxxxxxxxxxxxx}
\begin{equation}\label{eq5.66}
\frakI_2\le\frac{2^{16}}{u_1^{1/3}u_2^{1/3}}.
\end{equation}
\end{lemma}

\begin{proof}
(i)
Suppose first that $\ell=0$.
It follows from \eqref{eq5.31}, \eqref{eq5.35} and \eqref{eq5.51} that for every $\bfm^*=(m_2,m_3)\in\ZZZ_j^\ddagger$
and $(\alpha_1,\alpha_2)\in\Omega_j(\bfm^*;0)$, we have
\begin{equation}\label{eq5.67}
0\le\Lambda(\theta;u_1,u_2;H;\bfv;\bfm)\le\frac{H}{2}\min\left\{u_1u_2,\frac{1}{2^{3j}u_1u_2}\right\}.
\end{equation}
Note that it follows from \eqref{eq5.38} that $\vert\ZZZ_0^\ddagger\vert\le8$ and $\vert\ZZZ_j^\ddagger\vert\le4^{j+1}$ for $j\ge1$,
and so
\begin{equation}\label{eq5.68}
\vert\ZZZ_j^\ddagger\vert\le4^{j+2}
\end{equation}
for all $j\ge0$.
Combining this with \eqref{eq5.41} and \eqref{eq5.67}, we deduce that
\begin{align}\label{eq5.69}
&
\sum_{\bfm^*\in\ZZZ_j^\ddagger}\int_{\Omega_j(\bfm^*;0)}\Lambda(\theta;u_1,u_2;H;\bfv;\bfm)\,\dd\alpha_1\,\dd\alpha_2
\nonumber
\\
&\quad
\le4^{j+2}\lambda_2(\Omega_j(\bfm^*;0))\frac{H}{2}\min\left\{u_1u_2,\frac{1}{2^{3j}u_1u_2}\right\}
\nonumber
\\
&\quad
\le2^{2j+9}\min\left\{u_1u_2,\frac{1}{2^{3j}u_1u_2}\right\}.
\end{align}
We claim that
\textcolor{white}{xxxxxxxxxxxxxxxxxxxxxxxxxxxxxx}
\begin{equation}\label{eq5.70}
\sum_{j=0}^\infty2^{2j}\min\left\{u_1u_2,\frac{1}{2^{3j}u_1u_2}\right\}
\le\frac{4}{u_1^{1/3}u_2^{1/3}}.
\end{equation}
To see this, let $J^\ddagger$ denote the largest non-negative integer $j$ such that
\begin{equation}\label{eq5.71}
u_1u_2\le\frac{1}{2^{3j}u_1u_2},
\end{equation}
so that
\textcolor{white}{xxxxxxxxxxxxxxxxxxxxxxxxxxxxxx}
\begin{equation}\label{eq5.72}
2^{3J^\ddagger}\le\frac{1}{u_1^2u_2^2}
\quad\mbox{and}\quad
2^{3J^\ddagger+3}>\frac{1}{u_1^2u_2^2}.
\end{equation}
Then it follows from \eqref{eq5.71} and \eqref{eq5.72} that
\begin{align}\label{eq5.73}
&
\sum_{j=0}^{J^\ddagger}2^{2j}\min\left\{u_1u_2,\frac{1}{2^{3j}u_1u_2}\right\}
\le u_1u_2\sum_{j=0}^{J^\ddagger}2^{2j}
\nonumber
\\
&\quad
\le2u_1u_2\left(\frac{1}{u_1^2u_2^2}\right)^{2/3}
=\frac{2}{u_1^{1/3}u_2^{1/3}}
\end{align}
and
\textcolor{white}{xxxxxxxxxxxxxxxxxxxxxxxxxxxxxx}
\begin{align}\label{eq5.74}
&
\sum_{j=J^\ddagger+1}^\infty2^{2j}\min\left\{u_1u_2,\frac{1}{2^{3j}u_1u_2}\right\}
\le\frac{1}{u_1u_2}\sum_{j=J^\ddagger+1}^\infty\frac{1}{2^j}
\nonumber
\\
&\quad
\le\frac{2}{u_1u_2}(u_1^2u_2^2)^{1/3}
=\frac{2}{u_1^{1/3}u_2^{1/3}}.
\end{align}
Combining \eqref{eq5.73} and \eqref{eq5.74} leads to \eqref{eq5.70} which, together with \eqref{eq5.69},
leads to the inequality
\begin{equation}\label{eq5.75}
\sum_{j=0}^\infty\sum_{\bfm^*\in\ZZZ_j^\ddagger}\int_{\Omega_j(\bfm^*;0)}\Lambda(\theta;u_1,u_2;H;\bfv;\bfm)\,\dd\alpha_1\,\dd\alpha_2
\le\frac{2^{11}}{u_1^{1/3}u_2^{1/3}}.
\end{equation}

(ii)
Suppose next that $\ell\ge1$.
It follows from \eqref{eq5.31}, \eqref{eq5.35}, \eqref{eq5.40} and \eqref{eq5.51} that for every $\bfm^*=(m_2,m_3)\in\ZZZ_j^\ddagger$
and $(\alpha_1,\alpha_2)\in\Omega_j(\bfm^*;\ell)$, we have
\begin{equation}\label{eq5.76}
0\le\Lambda(\theta;u_1,u_2;H;\bfv;\bfm)\le\frac{H}{2^{2\ell-2}}\min\left\{u_1u_2,\frac{1}{2^{3j}u_1u_2}\right\}.
\end{equation}
Combining \eqref{eq5.41}, \eqref{eq5.68} and \eqref{eq5.76}, we deduce that
\begin{align}
&
\sum_{\bfm^*\in\ZZZ_j^\ddagger}\int_{\Omega_j(\bfm^*;\ell)}\Lambda(\theta;u_1,u_2;H;\bfv;\bfm)\,\dd\alpha_1\,\dd\alpha_2
\nonumber
\\
&\quad
\le4^{j+2}\lambda_2(\Omega_j(\bfm^*;\ell))\frac{H}{2^{2\ell-2}}\min\left\{u_1u_2,\frac{1}{2^{3j}u_1u_2}\right\}
\nonumber
\\
&\quad
\le2^{2j+4}\min\left\{\frac{1}{2^{\ell-8}},\frac{H}{2^{2\ell-2}}\right\}\min\left\{u_1u_2,\frac{1}{2^{3j}u_1u_2}\right\},
\nonumber
\end{align}
so that
\textcolor{white}{xxxxxxxxxxxxxxxxxxxxxxxxxxxxxx}
\begin{align}\label{eq5.77}
&
\sum_{\ell=1}^\infty\sum_{\bfm^*\in\ZZZ_j^\ddagger}\int_{\Omega_j(\bfm^*;\ell)}\Lambda(\theta;u_1,u_2;H;\bfv;\bfm)\,\dd\alpha_1\,\dd\alpha_2
\nonumber
\\
&\quad
\le2^{2j+4}\min\left\{u_1u_2,\frac{1}{2^{3j}u_1u_2}\right\}\sum_{\ell=1}^\infty\min\left\{\frac{1}{2^{\ell-8}},\frac{H}{2^{2\ell-2}}\right\}
\nonumber
\\
&\quad
\le2^{2j+13}\min\left\{u_1u_2,\frac{1}{2^{3j}u_1u_2}\right\}.
\end{align}
Combining \eqref{eq5.70} and \eqref{eq5.77} now leads to the inequality
\begin{equation}\label{eq5.78}
\sum_{j=0}^\infty\sum_{\ell=1}^\infty\sum_{\bfm^*\in\ZZZ_j^\ddagger}\int_{\Omega_j(\bfm^*;0)}\Lambda(\theta;u_1,u_2;H;\bfv;\bfm)\,\dd\alpha_1\,\dd\alpha_2
\le\frac{2^{15}}{u_1^{1/3}u_2^{1/3}}.
\end{equation}

The inequality \eqref{eq5.66} now follows on combining \eqref{eq5.54}, \eqref{eq5.75} and \eqref{eq5.78}.
\end{proof}

It now follows from \eqref{eq5.52}, \eqref{eq5.55} and \eqref{eq5.66} that
\begin{displaymath}
\int_0^1\int_0^1\Xi(\alpha_1,\alpha_2)\,\dd\alpha_1\,\dd\alpha_2\le\psi(u_1,u_2),
\end{displaymath}
where $\psi(u_1,u_2)$ is given by \eqref{eq5.1}.
Then for every real parameter $\kappa>1$, we have
\begin{equation}\label{eq5.79}
\lambda_2(\{(\alpha_1,\alpha_2)\in[0,1]^2:\Xi(\alpha_1,\alpha_2)\ge\kappa\psi(u_1,u_2)\})\le\frac{1}{\kappa}.
\end{equation}
\end{step6}


\begin{step7}
It remains to justify all the steps.
Choosing any real parameter $\kappa>1$, we remove the set
\textcolor{white}{xxxxxxxxxxxxxxxxxxxxxxxxxxxxxx}
\begin{displaymath}
\frakB(\kappa)=\{(\alpha_1,\alpha_2)\in[0,1]^2:\Xi(\alpha_1,\alpha_2)\ge\kappa\psi(u_1,u_2)\}
\end{displaymath}
of $\kappa$-bad direction vectors $\bfv_0=(\alpha_1,\alpha_2)\in[0,1]^2$.
In view of the estimate \eqref{eq5.79}, the remaining set
\textcolor{white}{xxxxxxxxxxxxxxxxxxxxxxxxxxxxxx}
\begin{displaymath}
\frakG(\kappa)=[0,1]^2\setminus\frakB(\kappa)
\end{displaymath}
of $\kappa$-good direction vectors satisfies
\begin{displaymath}
\lambda_2(\frakG(\kappa))\ge1-\frac{1}{\kappa}.
\end{displaymath}
Furthermore, for every $\bfv_0=(\alpha_1,\alpha_2)\in\frakG(\kappa)$, we have the identity
\begin{align}\label{eq5.80}
&
\frac{2}{H}\int_0^{H/2}\int_0^1\int_0^1\left(\sum_{\bfn\in\Zz^3}\chi_{B(\gamma_1,\gamma_2,h)}(M\bfn)\right)\dd\gamma_1\,\dd\gamma_2\,\dd h
\nonumber
\\
&\quad
=\frac{2}{H}\int_0^{H/2}\int_0^1\int_0^1
\left(\sum_{\bfm\in\Zz^3}\int_{B(\gamma_1,\gamma_2,h)}\ee^{-2\pi\ii M^{-1}\bfz\cdot\bfm}\,\dd\bfz\right)
\,\dd\gamma_1\,\dd\gamma_2\,\dd h,
\end{align}
and the quantitative argument in Steps 4--6 shows that the infinite sum on the right hand side of \eqref{eq5.80}
is absolutely convergent.
This justifies the use of the Poisson summation formula and makes the argument precise.
\end{step7}

This completes the proof of Lemma~\ref{lem51}.
\end{proof}

%
%

\section{Piecewise smooth barriers}\label{sec6}

Recall that in Theorems \ref{thm1} and~\ref{thm4}, we have made a rather specific restriction on the $2$-coloring on the $yz$-parallel square faces
of the $n$-cube $3$-manifold, that each of the red and green parts is the union of finitely many polygons.
Here we investigate how any polygon of the $2$-colorings can be replaced by a circle, an ellipse, or any other piecewise smooth closed curve.

\begin{theorem}\label{thm5}
Let $n\ge2$ be an integer, and let $\MMM$ be any $n$-cube $3$-manifold with barriers,
where each $yz$-parallel square face has a $2$-coloring such that each of the red and green parts is the union of finitely many regions
with edges that are piecewise smooth closed curves, satisfying some mild technical requirements as stated in Lemma~\ref{lem61},
and where the green part has positive area.
Suppose further that the Restriction on Red Coloring holds,
and that there is a local repetition color-split neighborhood on the $yz$-parallel square faces.
Then for almost every starting point and almost every direction $\bfv=(1,\alpha_1,\alpha_2)\in\Rr^3$,
the corresponding half-infinite $1$-direction geodesic is equidistributed in~$\MMM$.
\end{theorem}

Since the edges of a polygon have zero curvature, the buffer zone $\BBB_N$ within the local repetition color-split neighborhood
which is central to our argument can be taken as a long and narrow rectangular color-split strip.
To establish Theorem~\ref{thm5}, it is sufficient to deal with some buffer zone color-split strip which involves \textit{proper bending},
with non-zero curvature.
We therefore need to establish a suitable analog of Lemma~\ref{lem34}.
Clearly it suffices to establish a suitable analog of Lemma~\ref{lem51}.

The first step is to describe such a strip with proper bending.

Let $f:[a,b]\to\Rr$ be a monotonic and sufficiently smooth function.
More precisely, we require that its derivative $f'$ satisfies
\begin{displaymath}
\min_{x\in[a,b]}\vert f'(x)\vert>0,
\end{displaymath}
the curvature is non-zero in $[a,b]$, together with some mild technical requirements on higher order derivatives
which we shall describe later.

The strip with proper bending is then a $u$-neighborhood
\begin{equation}\label{eq6.1}
\CCC(f;u)=\CCC(f;a,b;u)=\{(x,y)\in\Rr^2:x\in[a,b]\mbox{ and }y\in[f(x)-u,f(x)+u]\}
\end{equation}
of a color-split curved edge
\begin{displaymath}
\CCC(f)=\CCC(f;a,b)=\{(x,y)\in\Rr^2:x\in[a,b]\mbox{ and }y=f(x)\},
\end{displaymath}
as illustrated in Figure~6.1.

\begin{displaymath}
\begin{array}{c}
\includegraphics[scale=0.8]{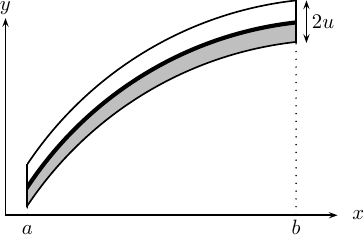}
\\
\mbox{Figure 6.1: the $u$-neighborhood $\CCC(f;u)$ of the color-split curved edge $\CCC(f)$}
\end{array}
\end{displaymath}

Clearly the $u$-neighborhood $\CCC(f;u)$ has area
\begin{displaymath}
\lambda_2(\CCC(f;u))=2u(b-a).
\end{displaymath}

We consider functions $f:[a,b]\to\Rr$ that are $3$-times continuously differentiable,
with positive constants $c_6,\ldots,c_{11}$, depending at most on the function $f$ in the interval $[a,b]$, such that
\textcolor{white}{xxxxxxxxxxxxxxxxxxxxxxxxxxxxxx}
\begin{align}
0<c_6=\min_{x\in[a,b]}\vert f'(x)\vert
&\le\max_{x\in[a,b]}\vert f'(x)\vert=c_7,
\label{eq6.2}
\\
0<c_8=\min_{x\in[a,b]}\vert f''(x)\vert
&\le\max_{x\in[a,b]}\vert f''(x)\vert=c_9,
\label{eq6.3}
\\
0<c_{10}=\min_{x\in[a,b]}\vert f'''(x)\vert
&\le\max_{x\in[a,b]}\vert f'''(x)\vert=c_{11}.
\label{eq6.4}
\end{align}
We also require the function $f$ to satisfy a mild technical condition, that there exists a constant $c_{12}=c_{12}(f;a,b)>0$,
depending at most on the function $f$ in the interval $[a,b]$, such that for any integer pair $(m_2,m_3)\in\Zz^2\setminus\{(0,0)\}$,
the interval $[a,b]$ is the union of at most $c_{12}$ subintervals such that the functions
\begin{equation}\label{eq6.5}
\frac{(m_2+m_3f'(x))^2}{f''(x)}
\quad\mbox{and}\quad
\frac{(m_2+m_3f'(x))^3}{f'(x)f''(x)}
\end{equation}
are monotonic in each of the subintervals.

\begin{lemma}\label{lem61}
Let $f:[a,b]\to\Rr$ be a $3$-times continuously differentiable function for which the conditions \eqref{eq6.2}--\eqref{eq6.4}
as well as the technical condition concerning the functions \eqref{eq6.5} hold.
For any even positive integer~$H$, let
\begin{displaymath}
G(\bfs_0;\alpha_1,\alpha_2;H)
=\vert\{j=0,1,\ldots,H-1:\bfs_0+j(\alpha_1,\alpha_2)\in\CCC(f;u)\}\vert,
\end{displaymath}
where the $u$-neighborhood $\CCC(f;u)$ is defined by \eqref{eq6.1}.
Then for any parameter $\kappa>1$, we have
\begin{displaymath}
\lambda_2\left(\left\{(\alpha_1,\alpha_2)\in[0,1]^2:G(\bfs_0;\alpha_1,\alpha_2;H)
\ge\frac{u(b-a)H}{4}-\kappa\psi(u)\right\}\right)\ge1-\frac{1}{\kappa},
\end{displaymath}
where
\textcolor{white}{xxxxxxxxxxxxxxxxxxxxxxxxxxxxxx}
\begin{equation}\label{eq6.6}
\psi(u)=c_{13}\left(u^{-1/2}+c_{14}\right),
\end{equation}
and the constants $c_{13}=c_{13}(f;a,b)>0$ and $c_{14}=c_{14}(f;a,b)>0$
depend at most on the function $f$ in the interval $[a,b]$.
\end{lemma}

\begin{remark}
Note that $(b-a)/2$ and $u$ play the roles of $u_1$ and $u_2$ in Lemma~\ref{lem51},
and that $G(\bfs_0;\alpha_1,\alpha_2;H)$ is the analog of the counting function $F(\bfs_0;\alpha_1,\alpha_2;H)$ there.
We also have non-zero curvature instead of rotation there.
\end{remark}

\begin{proof}[Proof of Lemma~\ref{lem61}]
We proceed by a number of steps corresponding to those in the proof of Lemma~\ref{lem51}.
We also use similar notation as much as possible.


\begin{step1}
Here we give a good description of $G(\bfs_0;\alpha_1,\alpha_2;H)$.
As before, we consider the lattice $\LLL(\alpha_1,\alpha_2)$ defined by \eqref{eq5.3} and \eqref{eq5.4}.
Then writing
\begin{displaymath}
B=(\CCC(f;u)-\bfs_0)\times[0,H)\subset\Rr^3,
\end{displaymath}
we then have
\textcolor{white}{xxxxxxxxxxxxxxxxxxxxxxxxxxxxxx}
\begin{displaymath}
G(\bfs_0;\alpha_1,\alpha_2;H)
=\vert\LLL(\alpha_1,\alpha_2)\cap B\vert
=\sum_{\bfn\in\Zz^3}\chi_B(M\bfn).
\end{displaymath}
Applying the Poisson summation formula \eqref{eq5.6} with $f=\chi_B$, we have, analogous to \eqref{eq5.7}, the formal identity
\begin{equation}\label{eq6.7}
G(\bfs_0;\alpha_1,\alpha_2;H)
=\sum_{\bfm\in\Zz^3}\int_B\ee^{-2\pi\ii M^{-1}\bfz\cdot\bfm}\,\dd\bfz.
\end{equation}
Note, in particular, that for $\bfm=\bzero$, we have
\begin{equation}\label{eq6.8}
\int_B\ee^{-2\pi\ii M^{-1}\bfz\cdot\bfm}\,\dd\bfz
=\lambda_3(B)
=\lambda_2(\CCC(f;u))H
=2u(b-a)H.
\end{equation}
Then, using \eqref{eq5.9} and analogous to \eqref{eq5.10}--\eqref{eq5.12}, for every $\bfm\in\Zz^3\setminus\{\bzero\}$, we can write
\textcolor{white}{xxxxxxxxxxxxxxxxxxxxxxxxxxxxxx}
\begin{equation}\label{eq6.9}
\int_B\ee^{-2\pi\ii M^{-1}\bfz\cdot\bfm}\,\dd\bfz
=\III(\bfm;\CCC(f;u)-\bfs_0)\JJJ(\bfm;\bfv;H),
\end{equation}
where
\textcolor{white}{xxxxxxxxxxxxxxxxxxxxxxxxxxxxxx}
\begin{equation}\label{eq6.10}
\III(\bfm;\CCC(f;u)-\bfs_0)
=\int_{\CCC(f;u)-\bfs_0}\ee^{2\pi\ii(z_1m_2+z_2m_3)}\,\dd z_1\,\dd z_2
\end{equation}
and
\textcolor{white}{xxxxxxxxxxxxxxxxxxxxxxxxxxxxxx}
\begin{equation}\label{eq6.11}
\JJJ(\bfm;\bfv;H)
=\int_0^H\ee^{-2\pi\ii z_3\bfv\cdot\bfm}\,\dd z_3
=\frac{1-\ee^{-2\pi\ii H\bfv\cdot\bfm}}{2\pi\ii\bfv\cdot\bfm}.
\end{equation}
To study the term \eqref{eq6.10}, note from \eqref{eq6.1} that a typical point in $\CCC(f;u)$ is of the form $(x,f(x)+t)$,
where $x\in[a,b]$ and $-u\le t\le u$.
We therefore use the substitution
\begin{displaymath}
(z_1,z_2)=(x,f(x)+t)-\bfs_0,
\quad
\mbox{with Jacobian}\quad
f'(x),
\end{displaymath}
so that, writing $\bfm^*=(m_2,m_3)$, we have
\begin{align}\label{eq6.12}
\III(\bfm;\CCC(f;u)-\bfs_0)
&
=\int_a^b\int_{-u}^u\ee^{2\pi\ii((x,f(x)+t)-\bfs_0)\cdot\bfm^*}f'(x)\,\dd t\,\dd x
\nonumber
\\
&
=\ee^{-2\pi\ii\bfs_0\cdot\bfm^*}\int_a^b\int_{-u}^u\ee^{2\pi\ii(x,f(x)+t)\cdot\bfm^*}f'(x)\,\dd t\,\dd x
\nonumber
\\
&
=\ee^{-2\pi\ii\bfs_0\cdot\bfm^*}\III_1(\bfm;\CCC(f;u))\III_2(\bfm;\CCC(f;u)),
\end{align}
where
\textcolor{white}{xxxxxxxxxxxxxxxxxxxxxxxxxxxxxx}
\begin{equation}\label{eq6.13}
\III_1(\bfm;\CCC(f;u))=\int_a^b\ee^{2\pi\ii(m_2x+m_3f(x))}f'(x)\,\dd x
\end{equation}
and
\textcolor{white}{xxxxxxxxxxxxxxxxxxxxxxxxxxxxxx}
\begin{equation}\label{eq6.14}
\III_2(\bfm;\CCC(f;u))=\int_{-u}^u\ee^{2\pi\ii tm_3}\,\dd t=\frac{\sin(2\pi um_3)}{\pi m_3},
\end{equation}
provided that $m_3\ne0$.
If $m_3=0$, then clearly
\begin{equation}\label{eq6.15}
\III_2(\bfm;\CCC(f;u))=\int_{-u}^u\dd t=2u.
\end{equation}

The analysis of the integral \eqref{eq6.13} is complicated.
However, the smoothness of the function $f$ makes it possible to have effective estimates.
We observe that it is possible to establish good estimates for complex exponential integrals
whenever the exponential function exhibits \textit{rapid fluctuations}.
The worst case scenario is when the exponential function $\ee^{2\pi\ii(m_2x+m_3f(x))}$ is almost constant
in a small neighborhood of some point~$x_0\in[a,b]$.

Let us consider some heuristics.
Write $h(x)=m_2x+m_3f(x)$.
Then we have the finite Taylor expansion
\begin{align}\label{eq6.16}
&
h(x)-h(x_0)
=m_2(x-x_0)+m_3(f(x)-f(x_0))
\nonumber
\\
&\quad
=(x-x_0)(m_2+m_3f'(x_0))+\frac{(x-x_0)^2m_3f''(x_0)}{2}+\frac{(x-x_0)^3m_3f'''(y)}{6}
\end{align}
for some appropriate $y$ between $x_0$ and~$x$.
Suppose further that there exists some $x_0\in[a,b]$ such that $h'(x_0)=m_2+m_3f'(x_0)=0$.
Then for every $x\in[a,b]$, the Taylor expansion \eqref{eq6.16} simplifies to
\begin{equation}\label{eq6.17}
h(x)-h(x_0)=\frac{(x-x_0)^2m_3f''(x_0)}{2}+\frac{(x-x_0)^3m_3f'''(y)}{6}
\end{equation}
for some $y\in[a,b]$.
The Taylor expansions \eqref{eq6.16} and \eqref{eq6.17} show that the complex exponential function
\textcolor{white}{xxxxxxxxxxxxxxxxxxxxxxxxxxxxxx}
\begin{displaymath}
\ee^{2\pi\ii h(x)}=\ee^{2\pi\ii(m_2x+m_3f(x))}
\end{displaymath}
starts to exhibit increasing fluctuations if $x$ is \textit{relatively far} from~$x_0$.

The estimation of \textit{fluctuating integrals} is a well known general problem in number theory
and especially in the theory of the Riemann zeta-function.
The following result can be found in the treatise of Titchmarsh \cite[Lemma~4.3]{titchmarsh86} on the latter.

\begin{lemma}\label{lem62}
Suppose that for real functions $F$ and~$G$, the quotient $F'(x)/G(x)$ is monotonic in the interval $[a,b]$,
and $\vert F'(x)/G(x)\vert\ge\mu>0$ for every $x\in[a,b]$.
Then
\begin{displaymath}
\left\vert\int_a^b\ee^{\ii F(x)}G(x)\,\dd x\right\vert\le\frac{4}{\mu}.
\end{displaymath}
\end{lemma}

To estimate the integral \eqref{eq6.13}, we consider the functions
\begin{displaymath}
F(x)=2\pi(m_2x+m_3f(x))
\quad\mbox{and}\quad
G(x)=f'(x)
\end{displaymath}
in the interval $[a,b]$.
Note that
\begin{equation}\label{eq6.18}
\frac{F'(x)}{G(x)}=\frac{2\pi(m_2+m_3f'(x))}{f'(x)}=\frac{2\pi m_2}{f'(x)}+2\pi m_3,
\end{equation}
and the derivative
\textcolor{white}{xxxxxxxxxxxxxxxxxxxxxxxxxxxxxx}
\begin{displaymath}
\frac{\dd}{\dd x}\frac{F'(x)}{G(x)}=-\frac{2\pi m_2f''(x)}{(f'(x))^2}
\end{displaymath}
has constant sign in $[a,b]$, since $f'(x)$ and $f''(x)$ are continuous and non-zero in $[a,b]$, in view of \eqref{eq6.2} and \eqref{eq6.3}.
Hence the function \eqref{eq6.18} is monotonic in $[a,b]$.
To apply Lemma~\ref{lem62}, we need to ensure that the function \eqref{eq6.18} is also bounded away from~$0$.
If no $x_0$ is in or near the interval $[a,b]$ such that $F'(x_0)=2\pi(m_2+m_3f'(x_0))=0$,
then we can find a good value for $\mu$ in Lemma~\ref{lem62}.
However, we clearly have a problem if there exists $x_0\in[a,b]$ such that $F'(x_0)=2\pi(m_2+m_3f'(x_0))=0$.

The following is the worst case scenario.
Suppose that $x_0\in(a,b)$.
Let $\delta>0$ be a \textit{small} parameter, to be specified later, such that $a<x_0-\delta<x_0+\delta<b$.
Then
\begin{align}\label{eq6.19}
&
\III_1(\bfm;\CCC(f;u))
=\int_a^b\ee^{2\pi\ii(m_2x+m_3f(x))}f'(x)\,\dd x
\nonumber
\\
&\quad
=\III_1^{(-)}(\bfm;\CCC(f;u))+\III_1^{(0)}(\bfm;\CCC(f;u))+\III_1^{(+)}(\bfm;\CCC(f;u)),
\end{align}
where
\textcolor{white}{xxxxxxxxxxxxxxxxxxxxxxxxxxxxxx}
\begin{align}
\III_1^{(-)}(\bfm;\CCC(f;u))
&=\int_a^{x_0-\delta}\ee^{2\pi\ii(m_2x+m_3f(x))}f'(x)\,\dd x,
\nonumber
\\
\III_1^{(0)}(\bfm;\CCC(f;u))
&=\int_{x_0-\delta}^{x_0+\delta}\ee^{2\pi\ii(m_2x+m_3f(x))}f'(x)\,\dd x,
\nonumber
\\
\III_1^{(+)}(\bfm;\CCC(f;u))
&=\int_{x_0+\delta}^b\ee^{2\pi\ii(m_2x+m_3f(x))}f'(x)\,\dd x.
\nonumber
\end{align}
For the integral $\III_1^{(0)}(\bfm;\CCC(f;u))$, we have the trivial estimate
\begin{equation}\label{eq6.20}
\vert\III_1^{(0)}(\bfm;\CCC(f;u))\vert\le2\delta\max_{x\in[a,b]}\vert f'(x)\vert=2c_7\delta,
\end{equation}
in view of \eqref{eq6.2}.
For the integral $\III_1^{(+)}(\bfm;\CCC(f;u))$, we use the Taylor expansion
\begin{align}
F'(x_0+\delta)
&
=2\pi(m_2+m_3f'(x_0+\delta))
\nonumber
\\
&
=2\pi(m_2+m_3f'(x_0))+2\pi m_3\left(\delta f''(x_0)+\frac{\delta^2f'''(y_+)}{2}\right)
\nonumber
\\
&
=2\pi m_3\left(\delta f''(x_0)+\frac{\delta^2f'''(y_+)}{2}\right)
\nonumber
\end{align}
for some $y_+\in(x_0,x_0+\delta)$.
Since $F'(x_0)/G(x_0)=0$ and $F'(x)/G(x)$ is monotonic in $[a,b]$, it follows that for every $x\in[x_0+\delta,b]$, we have
\begin{align}\label{eq6.21}
\left\vert\frac{F'(x)}{G(x)}\right\vert
&
\ge\left\vert\frac{F'(x_0+\delta)}{G(x_0+\delta)}\right\vert
=\frac{2\pi\vert m_3\vert}{\vert f'(x_0+\delta)\vert}\left\vert\delta f''(x_0)+\frac{\delta^2f'''(y_+)}{2}\right\vert
\nonumber
\\
&
\ge\frac{2\pi\vert m_3\vert}{\vert f'(x_0+\delta)\vert}\left(\delta\vert f''(x_0)\vert-\frac{\delta^2\vert f'''(y_+)\vert}{2}\right).
\end{align}
For the integral $\III_1^{(-)}(\bfm;\CCC(f;u))$, we use the Taylor expansion
\begin{align}
F'(x_0-\delta)
&
=2\pi(m_2+m_3f'(x_0-\delta))
\nonumber
\\
&
=2\pi(m_2+m_3f'(x_0))+2\pi m_3\left(-\delta f''(x_0)+\frac{\delta^2f'''(y_-)}{2}\right)
\nonumber
\\
&
=2\pi m_3\left(-\delta f''(x_0)+\frac{\delta^2f'''(y_-)}{2}\right)
\nonumber
\end{align}
for some $y_-\in(x_0-\delta,x_0)$.
Since $F'(x_0)/G(x_0)=0$ and $F'(x)/G(x)$ is monotonic in $[a,b]$, it follows that for every $x\in[a,x_0-\delta]$, we have
\begin{align}\label{eq6.22}
\left\vert\frac{F'(x)}{G(x)}\right\vert
&
\ge\left\vert\frac{F'(x_0-\delta)}{G(x_0-\delta)}\right\vert
=\frac{2\pi\vert m_3\vert}{\vert f'(x_0-\delta)\vert}\left\vert\delta f''(x_0)-\frac{\delta^2f'''(y_-)}{2}\right\vert
\nonumber
\\
&
\ge\frac{2\pi\vert m_3\vert}{\vert f'(x_0-\delta)\vert}\left(\delta\vert f''(x_0)\vert-\frac{\delta^2\vert f'''(y_-)\vert}{2}\right).
\end{align}
Combining \eqref{eq6.2}--\eqref{eq6.4}, \eqref{eq6.21} and \eqref{eq6.22},
we conclude that for any $x\in[a,x_0-\delta]$ or any $x\in[x_0+\delta,b]$, we have the lower bound
\begin{displaymath}
\left\vert\frac{F'(x)}{G(x)}\right\vert\ge\frac{2\pi\vert m_3\vert}{c_7}\left(c_8\delta-\frac{c_{11}\delta^2}{2}\right)\ge c_{15}\delta\vert m_3\vert,
\end{displaymath}
provided that $0<\delta\le c_{16}$, where $c_{15}=c_{15}(f;a,b)>0$ and $c_{16}=c_{16}(f;a,b)>0$
are constants depending at most on the function $f$ in the interval $[a,b]$.
Using this bound, Lemma~\ref{lem62} then gives the estimates
\begin{equation}\label{eq6.23}
\vert\III_1^{(\pm)}(\bfm;\CCC(f;u))\vert\le\frac{4}{c_{15}\delta\vert m_3\vert}.
\end{equation}
Combining \eqref{eq6.19}, \eqref{eq6.20} and \eqref{eq6.23}, we arrive at the estimate
\begin{equation}\label{eq6.24}
\vert\III_1(\bfm;\CCC(f;u))\vert\le2c_7\delta+\frac{8}{c_{15}\delta\vert m_3\vert}\le\frac{c_{17}}{\vert m_3\vert^{1/2}}
\end{equation}
if we choose $\delta$ in the range of $\vert m_3\vert^{-1/2}$,
where the constant $c_{17}=c_{17}(f;a,b)>0$ depends at most on the function $f$ in the interval $[a,b]$.

The above analysis and the estimate \eqref{eq6.24} are only valid provided that $m_3\ne0$.
If $m_3=0$, then \eqref{eq6.13} becomes
\begin{displaymath}
\III_1(\bfm;\CCC(f;u))=\int_a^b\ee^{2\pi\ii m_2x}f'(x)\,\dd x,
\end{displaymath}
and it follows from \eqref{eq6.2} that $\vert\III_1(\bfm;\CCC(f;u))\vert\le(b-a)c_7$.
For convenience, we can choose $c_{17}=c_{17}(f;a,b)>0$ sufficiently large so that
\begin{displaymath}
c_{17}\ge(b-a)c_7.
\end{displaymath}
Then
\textcolor{white}{xxxxxxxxxxxxxxxxxxxxxxxxxxxxxx}
\begin{equation}\label{eq6.25}
\vert\III_1(\bfm;\CCC(f;u))\vert\le c_{17}.
\end{equation}

We emphasize again that the above represents the worst case scenario.
Of course, the estimates \eqref{eq6.24} and \eqref{eq6.25} remain valid if $m_2+m_3f'(x)\ne0$ for any $x\in[a,b]$.

Combining \eqref{eq6.7}--\eqref{eq6.9} and \eqref{eq6.12}, we conclude that
\begin{align}
&
G(\bfs_0;\alpha_1,\alpha_2;H)-2u(b-a)H
\nonumber
\\
&\quad
=\sum_{\bfm\in\Zz^3\setminus\{\bzero\}}\ee^{-2\pi\ii\bfs_0\cdot\bfm^*}
\III_1(\bfm;\CCC(f;u))\III_2(\bfm;\CCC(f;u))\JJJ(\bfm;\bfv;H),
\nonumber
\end{align}
where estimates for the various factors in the summand are given by \eqref{eq6.11}, \eqref{eq6.14}, \eqref{eq6.15}, \eqref{eq6.24} and \eqref{eq6.25}.
\end{step1}


\begin{step2}
Here we mimic Step~2 in Section~\ref{sec5}, contract the interval $[0,H]$ in the third direction and average over all contractions
in a similar manner.
We then conclude, analogous to \eqref{eq5.23}, that
\begin{align}
&
G(\bfs_0;\alpha_1,\alpha_2;H)-u(b-a)H
\nonumber
\\
&\quad
\ge\sum_{\bfm\in\Zz^3\setminus\{\bzero\}}\ee^{-2\pi\ii\bfs_0\cdot\bfm^*}
\III_1(\bfm;\CCC(f;u))\III_2(\bfm;\CCC(f;u))\widetilde{\JJJ}(\bfm;\bfv;H),
\nonumber
\end{align}
where the factor $\widetilde{\JJJ}(\bfm;\bfv;H)$ is given by \eqref{eq5.22}.
\end{step2}


\begin{step3}
Here we also contract the intervals $[a,b]$ and $[f(x)-u,f(x)+u]$ in $\CCC(f;u)$ and average over all contractions.
More precisely, for every $\gamma_1$ and $\gamma_2$ satisfying $0\le\gamma_1\le(b-a)/2$ and $0\le\gamma_2\le1$,
consider the smaller set
\begin{align}
&
\Delta(\gamma_1,\gamma_2)
=\CCC(f;a+\gamma_1,b-\gamma_1;\gamma_2u)
\nonumber
\\
&\quad
=\{(x,y)\in\Rr^2:x\in[a+\gamma_1,b-\gamma_1]\mbox{ and }y\in[f(x)-\gamma_2u,f(x)+\gamma_2u]\}.
\nonumber
\end{align}
Then, analogous to \eqref{eq5.24} and \eqref{eq5.25}, we have
\begin{align}\label{eq6.26}
&
G(\bfs_0;\alpha_1,\alpha_2;H)-\frac{u(b-a)H}{4}
\nonumber
\\
&\quad
\ge\frac{4}{(b-a)H}\int_0^{H/2}\int_0^1\int_0^{(b-a)/2}
\left(\sum_{\bfm\in\Zz^3\setminus\{\bzero\}}\int_{B(\gamma_1,\gamma_2,h)}\ee^{-2\pi\ii M^{-1}\bfz\cdot\bfm}\,\dd\bfz\right)
\,\dd\gamma_1\,\dd\gamma_2\,\dd h,
\end{align}
where
\textcolor{white}{xxxxxxxxxxxxxxxxxxxxxxxxxxxxxx}
\begin{displaymath}
B(\gamma_1,\gamma_2,h)=(\Delta(\gamma_1,\gamma_2)-\bfs_0)\times[h,H-h),
\end{displaymath}
and, analogous to \eqref{eq5.26}, we have
\begin{equation}\label{eq6.27}
\int_{B(\gamma_1,\gamma_2,h)}\ee^{-2\pi\ii M^{-1}\bfz\cdot\bfm}\,\dd\bfz
=\III(\bfm;\Delta(\gamma_1,\gamma_2)-\bfs_0)\JJJ(\bfm;\bfv;H;h).
\end{equation}
Here $\JJJ(\bfm;\bfv;H;h)$ is given by \eqref{eq5.21} and
\begin{equation}\label{eq6.28}
\III(\bfm;\Delta(\gamma_1,\gamma_2)-\bfs_0)
=\ee^{-2\pi\ii\bfs_0\cdot\bfm^*}\III_1(\bfm;\Delta(\gamma_1,\gamma_2))\III_2(\bfm;\Delta(\gamma_1,\gamma_2)),
\end{equation}
where
\textcolor{white}{xxxxxxxxxxxxxxxxxxxxxxxxxxxxxx}
\begin{equation}\label{eq6.29}
\III_1(\bfm;\Delta(\gamma_1,\gamma_2))=\int_{a+\gamma_1}^{b-\gamma_1}\ee^{2\pi\ii(m_2x+m_3f(x))}f'(x)\,\dd x
\end{equation}
and
\textcolor{white}{xxxxxxxxxxxxxxxxxxxxxxxxxxxxxx}
\begin{equation}\label{eq6.30}
\III_2(\bfm;\Delta(\gamma_1,\gamma_2))=\int_{-\gamma_2u}^{\gamma_2u}\ee^{2\pi\ii tm_3}\,\dd t=\frac{\sin(2\pi\gamma_2um_3)}{\pi m_3},
\end{equation}
provided that $m_3\ne0$.
If $m_3=0$, then clearly
\begin{equation}\label{eq6.31}
\III_2(\bfm;\Delta(\gamma_1,\gamma_2))=2\gamma_2u.
\end{equation}
Write
\textcolor{white}{xxxxxxxxxxxxxxxxxxxxxxxxxxxxxx}
\begin{equation}\label{eq6.32}
\widetilde{\III}_1(\bfm;\CCC(f;u))
=\frac{2}{b-a}\int_0^{(b-a)/2}\int_{a+\gamma_1}^{b-\gamma_1}\ee^{2\pi\ii(m_2x+m_3f(x))}f'(x)\,\dd x\,\dd\gamma_1
\end{equation}
and
\textcolor{white}{xxxxxxxxxxxxxxxxxxxxxxxxxxxxxx}
\begin{align}\label{eq6.33}
&
\widetilde{\III}_2(\bfm;\CCC(f;u))
=\int_0^1\int_{-\gamma_2u}^{\gamma_2u}\ee^{2\pi\ii tm_3}\,\dd t\,\dd\gamma_2
\nonumber
\\
&\quad
=\int_0^1\frac{\sin(2\pi\gamma_2um_3)}{\pi m_3}\,\dd\gamma_2
=\frac{2\sin^2(\pi um_3)}{u(\pi m_3)^2},
\end{align}
provided that $m_3\ne0$.
If $m_3=0$, then clearly
\begin{equation}\label{eq6.34}
\widetilde{\III}_2(\bfm;\CCC(f;u))
=\int_0^1\int_{-\gamma_2u}^{\gamma_2u}\dd t\,\dd\gamma_2
=\int_0^12\gamma_2u\,\dd\gamma_2
=u.
\end{equation}
It now follows from \eqref{eq5.22} and \eqref{eq6.26}--\eqref{eq6.34} that
\begin{align}\label{eq6.35}
&
G(\bfs_0;\alpha_1,\alpha_2;H)-\frac{u(b-a)H}{4}
\nonumber
\\
&\quad
\ge\sum_{\bfm\in\Zz^3\setminus\{\bzero\}}\ee^{-2\pi\ii\bfs_0\cdot\bfm^*}
\widetilde{\III}_1(\bfm;\CCC(f;u))\widetilde{\III}_2(\bfm;\CCC(f;u))\widetilde{\JJJ}(\bfm;\bfv;H).
\end{align}
\end{step3}


\begin{step4}
It follows from \eqref{eq6.35} that
\begin{displaymath}
G(\bfs_0;\alpha_1,\alpha_2;H)\ge\frac{u(b-a)H}{4}+\sum_{\bfm\in\Zz^3\setminus\{\bzero\}}\Lambda(\CCC(f;u);H;\bfv;\bfm),
\end{displaymath}
where
\begin{equation}\label{eq6.36}
\Lambda(\CCC(f;u);H;\bfv;\bfm)
=\vert\widetilde{\III}_1(\bfm;\CCC(f;u))\vert\,\vert\widetilde{\III}_2(\bfm;\CCC(f;u))\vert\,\vert\widetilde{\JJJ}(\bfm;\bfv;H)\vert.
\end{equation}

As before in Section~\ref{sec5}, the inequalities in Steps 1--3 are at this stage only formal inequalities,
as we have not considered the question of convergence.
Write
\begin{equation}\label{eq6.37}
\Xi(\alpha_1,\alpha_2)=\sum_{\Zz^3\setminus\{\bzero\}}\Lambda(\CCC(f;u);H;\bfv;\bfm),
\end{equation}
where $\bfv=(1,\alpha_1,\alpha_2)$.
We shall use the first-moment method and analyze the average
\begin{displaymath}
\int_0^1\int_0^1\Xi(\alpha_1,\alpha_2)\,\dd\alpha_1\,\dd\alpha_2,
\end{displaymath}
and remove those directions $(\alpha_1,\alpha_2)$ for which $\Xi(\alpha_1,\alpha_2)$ is substantially larger than the average.
In this step, we do some further preparation.

First of all, note that the estimates \eqref{eq6.24} and \eqref{eq6.25} are independent of $a$ and~$b$.
Combining these with \eqref{eq6.29} and \eqref{eq6.32}, we obtain the bound
\begin{equation}\label{eq6.38}
\vert\widetilde{\III}_1(\bfm;\CCC(f;u))\vert\le c_{17}\min\left\{\frac{1}{\vert m_3\vert^{1/2}},1\right\}.
\end{equation}
On the other hand, using \eqref{eq6.30}, \eqref{eq6.31} and the inequality $\vert\sin y\vert\le\min\{\vert y\vert,1\}$, which holds for every $y\in\Rr$,
we obtain the bound
\begin{equation}\label{eq6.39}
\vert\widetilde{\III}_2(\bfm;\CCC(f;u))\vert\le2\min\left\{u,\frac{1}{u(\pi m_3)^2}\right\}.
\end{equation}
Meanwhile, the estimate for the term $\vert\widetilde{\JJJ}(\bfm;\bfv;H)\vert$ is given by \eqref{eq5.35}.
Furthermore, since $H$ is even, in view of \eqref{eq5.36}, we may assume that $\bfm^*=(m_2,m_3)\in\Zz^2\setminus\{(0,0)\}$.

Finally, the definitions and estimates \eqref{eq5.37}--\eqref{eq5.41} concerning $\ZZZ_j$ and $\Omega_j(\bfm^*;\ell)$,
where $j=0,1,2,3,\ldots$ and $\ell=0,1,2,3,\ldots,$ remain valid.
\end{step4}


\begin{step5}
Recall our comment in Step~1 that the complex exponential function
\begin{displaymath}
\ee^{2\pi\ii h(x)}=\ee^{2\pi\ii(m_2x+m_3f(x))}
\end{displaymath}
starts to exhibit increasing fluctuations if $x$ is \textit{relatively far} from any root $x_0$ of the equation $m_2+m_3f'(x)=0$.
If $x_0$ is far from the interval $[a,b]$, then the lack of fluctuation near $x_0$ is not a problem.
However, it is a serious problem if $x_0\in[a,b]$.

Motivated by this observation, we split the argument into two cases, the first of which
and some of its consequences are summarized by the following lemma.

\begin{lemma}\label{lem63}
Suppose that $(m_2,m_3)\in\Zz^2\setminus\{(0,0\}$ and the condition
\begin{equation}\label{eq6.40}
\min_{a\le x\le b}\vert m_2+m_3f'(x)\vert
<\frac{1}{4}\max\{\vert m_2\vert,\vert m_3\vert\}\min\left\{\min_{a\le x\le b}\vert f'(x)\vert,1\right\}
\end{equation}
is satisfied.
Then

\emph{(i)} $m_2m_3\ne0$; and

\emph{(ii)} the inequality
\textcolor{white}{xxxxxxxxxxxxxxxxxxxxxxxxxxxxxx}
\begin{equation}\label{eq6.41}
\frac{1}{c_{18}}\le\left\vert\frac{m_2}{m_3}\right\vert\le c_{18}
\end{equation}
holds, where the constant $c_{18}=c_{18}(f;a,b)>1$ depends at most on the function $f$ in the interval $[a,b]$.
\end{lemma}

\begin{proof}
(i)
If $m_2=0$, then the inequality \eqref{eq6.40} becomes
\begin{displaymath}
\vert m_3\vert\min_{a\le x\le b}\vert f'(x)\vert
<\frac{1}{4}\vert m_3\vert\min\left\{\min_{a\le x\le b}\vert f'(x)\vert,1\right\}
\end{displaymath}
which is clearly absurd.
If $m_3=0$, then the inequality \eqref{eq6.40} becomes
\begin{displaymath}
\vert m_2\vert
<\frac{1}{4}\vert m_2\vert\min\left\{\min_{a\le x\le b}\vert f'(x)\vert,1\right\}
\end{displaymath}
which is also clearly absurd.
Thus $m_2m_3\ne0$ as claimed.

(ii)
Dividing both sides by $m_3\ne0$, the inequality \eqref{eq6.40} becomes
\begin{equation}\label{eq6.42}
\min_{a\le x\le b}\left\vert\frac{m_2}{m_3}+f'(x)\right\vert
<\frac{1}{4}\max\left\{\left\vert\frac{m_2}{m_3}\right\vert,1\right\}\min\left\{\min_{a\le x\le b}\vert f'(x)\vert,1\right\}.
\end{equation}
Suppose first that
\textcolor{white}{xxxxxxxxxxxxxxxxxxxxxxxxxxxxxx}
\begin{displaymath}
\frac{m_2}{m_3}=C_L,
\end{displaymath}
where $\vert C_L\vert$ is \textit{large}.
Then using \eqref{eq6.2}, the left and right sides of \eqref{eq6.42} are
\begin{displaymath}
\min_{a\le x\le b}\left\vert\frac{m_2}{m_3}+f'(x)\right\vert
=\min_{a\le x\le b}\vert C_L+f'(x)\vert
\ge\vert C_L\vert-\max_{a\le x\le b}\vert f'(x)\vert
\ge\vert C_L\vert-c_7
\end{displaymath}
and
\begin{displaymath}
\frac{1}{4}\max\left\{\left\vert\frac{m_2}{m_3}\right\vert,1\right\}\min\left\{\min_{a\le x\le b}\vert f'(x)\vert,1\right\}
=\frac{1}{4}\vert C_L\vert\min\left\{\min_{a\le x\le b}\vert f'(x)\vert,1\right\}
\le\frac{1}{4}\vert C_L\vert
\end{displaymath}
respectively.
Clearly
\textcolor{white}{xxxxxxxxxxxxxxxxxxxxxxxxxxxxxx}
\begin{displaymath}
\vert C_L\vert-c_7>\frac{1}{4}\vert C_L\vert
\end{displaymath}
if $\vert C_L\vert$ is sufficiently large in terms of the function $f$ in the interval $[a,b]$, so that \eqref{eq6.42} fails.
This establishes the upper bound in \eqref{eq6.41}.
Suppose next that
\begin{displaymath}
\frac{m_2}{m_3}=C_S,
\end{displaymath}
where $\vert C_S\vert$ is \textit{small}.
Then the left and right sides of \eqref{eq6.42} are
\begin{displaymath}
\min_{a\le x\le b}\left\vert\frac{m_2}{m_3}+f'(x)\right\vert
=\min_{a\le x\le b}\vert C_S+f'(x)\vert
\ge\min_{a\le x\le b}\vert f'(x)\vert-\vert C_S\vert
\end{displaymath}
and
\begin{displaymath}
\frac{1}{4}\max\left\{\left\vert\frac{m_2}{m_3}\right\vert,1\right\}\min\left\{\min_{a\le x\le b}\vert f'(x)\vert,1\right\}
=\frac{1}{4}\min\left\{\min_{a\le x\le b}\vert f'(x)\vert,1\right\}
\le\frac{1}{4}\min_{a\le x\le b}\vert f'(x)\vert
\end{displaymath}
respectively.
Clearly
\textcolor{white}{xxxxxxxxxxxxxxxxxxxxxxxxxxxxxx}
\begin{displaymath}
\min_{a\le x\le b}\vert f'(x)\vert-\vert C_S\vert>\frac{1}{4}\min_{a\le x\le b}\vert f'(x)\vert
\end{displaymath}
if $\vert C_S\vert$ is sufficiently small in terms of the function $f$ in the interval $[a,b]$, so that \eqref{eq6.42} fails.
This establishes the lower bound in \eqref{eq6.41}.
\end{proof}

We also need to consider the case when the condition \eqref{eq6.40} is violated, so that the opposite condition
\begin{equation}\label{eq6.43}
\min_{a\le x\le b}\vert m_2+m_3f'(x)\vert
\ge\frac{1}{4}\max\{\vert m_2\vert,\vert m_3\vert\}\min\left\{\min_{a\le x\le b}\vert f'(x)\vert,1\right\}
\end{equation}
is satisfied.

Accordingly, for every $j=1,2,3,\ldots,$ write $\ZZZ_j=\ZZZ_j^\dagger\cup\ZZZ_j^\ddagger$, where
\begin{align}
\ZZZ_j^\dagger
&
=\{(m_2,m_3)\in\ZZZ_j:\mbox{\eqref{eq6.40} holds}\},
\nonumber
\\
\ZZZ_j^\ddagger
&
=\{(m_2,m_3)\in\ZZZ_j:\mbox{\eqref{eq6.43} holds}\}.
\nonumber
\end{align}
\end{step5}


\begin{step6}
Using \eqref{eq6.37} and taking into account that the identity \eqref{eq5.36} holds whenever $\bfm^*=(m_2,m_3)\ne(0,0)$,
we now have the trivial upper bound
\begin{equation}\label{eq6.44}
\int_0^1\int_0^1\Xi(\alpha_1,\alpha_2)\,\dd\alpha_1\,\dd\alpha_2\le\frakI_1+\frakI_2,
\end{equation}
where
\textcolor{white}{xxxxxxxxxxxxxxxxxxxxxxxxxxxxxx}
\begin{align}
\frakI_1
&
=\sum_{j=0}^\infty\sum_{\ell=0}^\infty\sum_{\bfm^*\in\ZZZ_j^\dagger}\int_{\Omega_j(\bfm^*;\ell)}
\Lambda(\CCC(f;u);H;\bfv;\bfm)\,\dd\alpha_1\,\dd\alpha_2,
\label{eq6.45}
\\
\frakI_2
&
=\sum_{j=0}^\infty\sum_{\ell=0}^\infty\sum_{\bfm^*\in\ZZZ_j^\ddagger}\int_{\Omega_j(\bfm^*;\ell)}
\Lambda(\CCC(f;u);H;\bfv;\bfm)\,\dd\alpha_1\,\dd\alpha_2.
\label{eq6.46}
\end{align}

\begin{lemma}\label{lem64}
We have
\textcolor{white}{xxxxxxxxxxxxxxxxxxxxxxxxxxxxxx}
\begin{equation}\label{eq6.47}
\frakI_1\le\frac{c_{19}}{u^{1/2}},
\end{equation}
where the constant $c_{19}=c_{19}(f;a,b)>0$ depends at most on the function $f$ in the interval $[a,b]$.
\end{lemma}

\begin{proof}
Suppose that $(m_2,m_3)\in\ZZZ_j^\dagger$.
It follows from \eqref{eq5.38} and \eqref{eq6.41} that
\begin{equation}\label{eq6.48}
c_{20}2^j\le\min\{\vert m_2\vert,\vert m_3\vert\}\le\max\{\vert m_2\vert,\vert m_3\vert\}\le2^j,
\end{equation}
where the constant $c_{20}=c_{20}(f;a,b)>0$ satisfies $2c_{18}c_{20}=1$.

If $(\alpha_1,\alpha_2)\in\Omega_j(\bfm^*;0)$, then it follows from \eqref{eq5.35}, \eqref{eq6.36}, \eqref{eq6.38} and \eqref{eq6.39} that
\begin{equation}\label{eq6.49}
\Lambda(\CCC(f;u);H;\bfv;\bfm)\le\frac{c_{17}H}{\vert m_3\vert^{1/2}}\min\left\{u,\frac{1}{u(\pi m_3)^2}\right\}.
\end{equation}
Using \eqref{eq5.41}, \eqref{eq6.48} and \eqref{eq6.49}, we deduce that
\begin{align}\label{eq6.50}
&
\sum_{\bfm^*\in\ZZZ_j^\dagger}\int_{\Omega_j(\bfm^*;0)}\Lambda(\CCC(f;u);H;\bfv;\bfm)\,\dd\alpha_1\,\dd\alpha_2
\nonumber
\\
&\quad
\le4\sum_{c_{20}2^j\le m_2\le2^j}\sum_{c_{20}2^j\le m_3\le2^j}\lambda_2(\Omega_j(\bfm^*;0))\Lambda(\CCC(f;u);H;\bfv;\bfm)
\nonumber
\\
&\quad
\le2^8c_{17}\sum_{c_{20}2^j\le m_3\le2^j}\frac{2^j}{m_3^{1/2}}\min\left\{u,\frac{1}{u(\pi m_3)^2}\right\}.
\end{align}
Note next that
\begin{equation}\label{eq6.51}
\sum_{c_{20}2^j\le m_3\le2^j}\frac{2^j}{m_3^{1/2}}\min\left\{u,\frac{1}{u(\pi m_3)^2}\right\}
\le\sum_{c_{20}2^j\le m_3\le2^j}\frac{2^j}{u\pi^2m_3^{5/2}}
\le\frac{2^{2j}}{u\pi^2(c_{20}2^j)^{5/2}}.
\end{equation}
Combining \eqref{eq6.50} and \eqref{eq6.51}, we deduce that
\begin{equation}\label{eq6.52}
\sum_{\bfm^*\in\ZZZ_j^\dagger}\int_{\Omega_j(\bfm^*;0)}\Lambda(\CCC(f;u);H;\bfv;\bfm)\,\dd\alpha_1\,\dd\alpha_2\le\frac{c_{21}}{u2^{j/2}},
\end{equation}
where the constant $c_{21}=c_{21}(f;a,b)>0$ satisfies $c_{21}c_{20}^{5/2}\pi^2=2^8c_{17}$.
Note also that
\begin{equation}\label{eq6.53}
\sum_{c_{20}2^j\le m_3\le2^j}\frac{2^j}{m_3^{1/2}}\min\left\{u,\frac{1}{u(\pi m_3)^2}\right\}
\le\frac{u2^{2j}}{(c_{20}2^j)^{1/2}}
=\frac{u2^{3j/2}}{c_{20}^{1/2}}.
\end{equation}
Combining \eqref{eq6.50} and \eqref{eq6.53}, we deduce that
\begin{equation}\label{eq6.54}
\sum_{\bfm^*\in\ZZZ_j^\dagger}\int_{\Omega_j(\bfm^*;0)}\Lambda(\CCC(f;u);H;\bfv;\bfm)\,\dd\alpha_1\,\dd\alpha_2\le c_{22}u2^{3j/2},
\end{equation}
where the constant $c_{22}=c_{22}(f;a,b)>0$ satisfies $c_{22}c_{20}^{1/2}=2^8c_{17}$.

If $(\alpha_1,\alpha_2)\in\Omega_j(\bfm^*;\ell)$, where $\ell\ge1$,
then it follows from \eqref{eq5.35}, \eqref{eq5.40}, \eqref{eq6.36}, \eqref{eq6.38} and \eqref{eq6.39} that
\begin{equation}\label{eq6.55}
\Lambda(\CCC(f;u);H;\bfv;\bfm)\le\frac{2^4c_{17}H}{\pi^24^\ell\vert m_3\vert^{1/2}}\min\left\{u,\frac{1}{u(\pi m_3)^2}\right\}.
\end{equation}
Using \eqref{eq5.41}, \eqref{eq6.48}, \eqref{eq6.51} and \eqref{eq6.55}, we deduce that
\begin{align}
&
\sum_{\bfm^*\in\ZZZ_j^\dagger}\int_{\Omega_j(\bfm^*;\ell)}\Lambda(\CCC(f;u);H;\bfv;\bfm)\,\dd\alpha_1\,\dd\alpha_2
\nonumber
\\
&\quad
\le4\sum_{c_{20}2^j\le m_2\le2^j}\sum_{c_{20}2^j\le m_3\le2^j}\lambda_2(\Omega_j(\bfm^*;\ell))\Lambda(\CCC(f;u);H;\bfv;\bfm)
\nonumber
\\
&\quad
\le\frac{2^{12}c_{17}}{\pi^22^\ell}\sum_{c_{20}2^j\le m_3\le2^j}\frac{2^j}{m_3^{1/2}}\min\left\{u,\frac{1}{u(\pi m_3)^2}\right\}
\le\frac{c_{23}}{u2^\ell2^{j/2}},
\nonumber
\end{align}
where the constant $c_{23}=c_{23}(f;a,b)>0$ depends at most on the function $f$ in the interval $[a,b]$,
so that
\begin{equation}\label{eq6.56}
\sum_{\ell=1}^\infty\sum_{\bfm^*\in\ZZZ_j^\dagger}\int_{\Omega_j(\bfm^*;\ell)}\Lambda(\CCC(f;u);H;\bfv;\bfm)\,\dd\alpha_1\,\dd\alpha_2
\le\frac{c_{23}}{u2^{j/2}}.
\end{equation}
Meanwhile, using \eqref{eq6.53} instead of \eqref{eq6.51}, we deduce that
\begin{displaymath}
\sum_{\bfm^*\in\ZZZ_j^\dagger}\int_{\Omega_j(\bfm^*;\ell)}\Lambda(\CCC(f;u);H;\bfv;\bfm)\,\dd\alpha_1\,\dd\alpha_2
\le\frac{c_{24}u2^{3j/2}}{2^\ell},
\end{displaymath}
where the constant $c_{24}=c_{24}(f;a,b)>0$ depends at most on the function $f$ in the interval $[a,b]$,
so that
\begin{equation}\label{eq6.57}
\sum_{\ell=1}^\infty\sum_{\bfm^*\in\ZZZ_j^\dagger}\int_{\Omega_j(\bfm^*;\ell)}\Lambda(\CCC(f;u);H;\bfv;\bfm)\,\dd\alpha_1\,\dd\alpha_2
\le c_{24}u2^{3j/2}.
\end{equation}

Let $J^\dagger$ denote the largest non-negative integer such that
\begin{displaymath}
2^j\le\frac{1}{u},
\end{displaymath}
so that
\textcolor{white}{xxxxxxxxxxxxxxxxxxxxxxxxxxxxxx}
\begin{equation}\label{eq6.58}
2^{J^\dagger}\le\frac{1}{u}
\quad\mbox{and}\quad
2^{J^\dagger+1}>\frac{1}{u}.
\end{equation}

Combining \eqref{eq6.45}, \eqref{eq6.52}, \eqref{eq6.54}, \eqref{eq6.56}, and \eqref{eq6.57}, we conclude that
\begin{equation}\label{eq6.59}
\frakI_1\le c_{25}u\sum_{j=0}^{J^\dagger}2^{3j/2}+\frac{c_{26}}{u}\sum_{j=J^\dagger+1}^\infty\frac{1}{2^{j/2}},
\end{equation}
where the constant $c_{25}=c_{25}(f;a,b)>0$ satisfies $c_{25}=\max\{c_{22},c_{24}\}$
while the constant $c_{26}=c_{26}(f;a,b)>0$ satisfies $c_{26}=\max\{c_{21},c_{23}\}$.
The inequality \eqref{eq6.47} follows on combining \eqref{eq6.58} and \eqref{eq6.59}.
\end{proof}

\begin{lemma}\label{lem65}
We have
\textcolor{white}{xxxxxxxxxxxxxxxxxxxxxxxxxxxxxx}
\begin{equation}\label{eq6.60}
\frakI_2\le c_{27},
\end{equation}
where the constant $c_{27}=c_{27}(f;a,b)>0$ depends at most on the function $f$ in the interval $[a,b]$.
\end{lemma}

\begin{proof}
The proof is in three parts.

\begin{part1}
Here we obtain new bounds for the term $\widetilde{\III}_1(\bfm;C(f;u))$ in the case when the inequality \eqref{eq6.43} holds.
Indeed, combining \eqref{eq6.2} and \eqref{eq6.43}, we obtain the lower bound
\textcolor{white}{xxxxxxxxxxxxxxxxxxxxxxxxxxxxxx}
\begin{equation}\label{eq6.61}
\min_{a\le x\le b}\vert m_2+m_3f'(x)\vert
\ge\frac{1}{4}\min\{c_6,1\}\max\{\vert m_2\vert,\vert m_3\vert\}.
\end{equation}
We make use the simple identity
\begin{equation}\label{eq6.62}
\frac{\dd}{\dd x}\left(\frac{\ee^{2\pi\ii(m_2x+m_3f(x))}}{2\pi\ii(m_2+m_3f'(x))}\right)
=\ee^{2\pi\ii(m_2x+m_3f(x))}+\frac{\ii m_3\ee^{2\pi\ii(m_2x+m_3f(x))}f''(x)}{2\pi(m_2+m_3f'(x))^2}.
\end{equation}
Combining \eqref{eq6.62} and integration by parts, we obtain
\begin{align}\label{eq6.63}
&
\int_{a+\gamma_1}^{b-\gamma_1}\ee^{2\pi\ii(m_2x+m_3f(x))}f'(x)\,\dd x
+\ii\int_{a+\gamma_1}^{b-\gamma_1}\frac{m_3\ee^{2\pi\ii(m_2x+m_3f(x))}f''(x)}{2\pi(m_2+m_3f'(x))^2}f'(x)\,\dd x
\nonumber
\\
&\quad
=\int_{a+\gamma_1}^{b-\gamma_1}\frac{\dd}{\dd x}\left(\frac{\ee^{2\pi\ii(m_2x+m_3f(x))}}{2\pi\ii(m_2+m_3f'(x))}\right)f'(x)\,\dd x
\nonumber
\\
&\quad
=\left[\frac{\ee^{2\pi\ii(m_2x+m_3f(x))}}{2\pi\ii(m_2+m_3f'(x))}f'(x)\right]_{a+\gamma_1}^{b-\gamma_1}
+\ii\int_{a+\gamma_1}^{b-\gamma_1}\frac{\ee^{2\pi\ii(m_2x+m_3f(x))}}{2\pi(m_2+m_3f'(x))}f''(x)\,\dd x.
\end{align}
Combining \eqref{eq6.32} and \eqref{eq6.63}, we have
\begin{equation}\label{eq6.64}
\widetilde{\III}_1(\bfm;C(f;u))=\ii\,\widetilde{\III}_1^{(1)}-\ii\,\widetilde{\III}_1^{(2)}+\ii\,\widetilde{\III}_1^{(3)}-\ii\,\widetilde{\III}_1^{(4)},
\end{equation}
where
\textcolor{white}{xxxxxxxxxxxxxxxxxxxxxxxxxxxxxx}
\begin{align}
\widetilde{\III}_1^{(1)}
&
=\frac{2}{b-a}\int_0^{(b-a)/2}
\int_{a+\gamma_1}^{b-\gamma_1}\frac{\ee^{2\pi\ii(m_2x+m_3f(x))}}{2\pi(m_2+m_3f'(x))}f''(x)\,\dd x\,\dd\gamma_1,
\label{eq6.65}
\\
\widetilde{\III}_1^{(2)}
&
=\frac{2}{b-a}\int_0^{(b-a)/2}
\int_{a+\gamma_1}^{b-\gamma_1}\frac{m_3\ee^{2\pi\ii(m_2x+m_3f(x))}f''(x)}{2\pi(m_2+m_3f'(x))^2}f'(x)\,\dd x\,\dd\gamma_1,
\label{eq6.66}
\\
\widetilde{\III}_1^{(3)}
&
=\frac{2}{b-a}\int_0^{(b-a)/2}
\frac{\ee^{2\pi\ii(m_2(a+\gamma_1)+m_3f(a+\gamma_1))}}{2\pi(m_2+m_3f'(a+\gamma_1))}f'(a+\gamma_1)\,\dd\gamma_1
\nonumber
\\
&
=\frac{2}{b-a}\int_a^{(a+b)/2}
\frac{\ee^{2\pi\ii(m_2(x)+m_3f(x))}}{2\pi(m_2+m_3f'(x))}f'(x)\,\dd x,
\label{eq6.67}
\\
\widetilde{\III}_1^{(4)}
&
=\frac{2}{b-a}\int_0^{(b-a)/2}
\frac{\ee^{2\pi\ii(m_2(b-\gamma_1)+m_3f(b-\gamma_1))}}{2\pi(m_2+m_3f'(b-\gamma_1))}f'(b-\gamma_1)\,\dd\gamma_1
\nonumber
\\
&
=\frac{2}{b-a}\int_{(a+b)/2}^b
\frac{\ee^{2\pi\ii(m_2(x)+m_3f(x))}}{2\pi(m_2+m_3f'(x))}f'(x)\,\dd x.
\label{eq6.68}
\end{align}

To study the integral $\widetilde{\III}_1^{(1)}$, note that
\begin{displaymath}
\int_{a+\gamma_1}^{b-\gamma_1}\frac{\ee^{2\pi\ii(m_2x+m_3f(x))}}{2\pi(m_2+m_3f'(x))}f''(x)\,\dd x=\int_{a+\gamma_1}^{b-\gamma_1}\ee^{\ii F(x)}G(x)\,\dd x,
\end{displaymath}
where
\textcolor{white}{xxxxxxxxxxxxxxxxxxxxxxxxxxxxxx}
\begin{equation}\label{eq6.69}
F(x)=2\pi(m_2x+m_3f(x))
\quad\mbox{and}\quad
G(x)=\frac{f''(x)}{2\pi(m_2+m_3f'(x))},
\end{equation}
so that
\textcolor{white}{xxxxxxxxxxxxxxxxxxxxxxxxxxxxxx}
\begin{displaymath}
\frac{F'(x)}{G(x)}
=\frac{4\pi^2(m_2+m_3f'(x))^2}{f''(x)}
=4\pi^2H(x),
\end{displaymath}
where
\textcolor{white}{xxxxxxxxxxxxxxxxxxxxxxxxxxxxxx}
\begin{displaymath}
H(x)=\frac{(m_2+m_3f'(x))^2}{f''(x)}.
\end{displaymath}
Using \eqref{eq6.3}, \eqref{eq6.61} and \eqref{eq6.69}, we deduce the lower bound
\begin{displaymath}
\left\vert\frac{F'(x)}{G(x)}\right\vert
=\frac{4\pi^2\vert m_2+m_3f'(x)\vert^2}{\vert f''(x)\vert}
\ge c_{28}(\max\{\vert m_2\vert,\vert m_3\vert\})^2,
\end{displaymath}
where the constant $c_{28}=c_{28}(f;a,b)>0$ depends at most on the function $f$ in the interval $[a,b]$.
Furthermore, the technical condition implies that we can apply Lemma~\ref{lem62} to each of the at most $c_{12}$
subintervals of $[a+\gamma_1,b-\gamma_1]$, and this leads to the bound
\begin{displaymath}
\left\vert\int_{a+\gamma_1}^{b-\gamma_1}\frac{\ee^{2\pi(m_2x+m_3f(x))}}{2\pi\ii(m_2+m_3f'(x))}f''(x)\,\dd x\right\vert
\le\frac{c_{29}}{(\max\{\vert m_2\vert,\vert m_3\vert\})^2},
\end{displaymath}
where the constant $c_{29}=c_{29}(f;a,b)>0$ depends at most on the function $f$ in the interval $[a,b]$.
It then follows trivially from \eqref{eq6.65} that
\begin{equation}\label{eq6.70}
\vert\widetilde{\III}_1^{(1)}\vert\le\frac{c_{29}}{(\max\{\vert m_2\vert,\vert m_3\vert\})^2}.
\end{equation}

To study the integral $\widetilde{\III}_1^{(2)}$, note that
\begin{displaymath}
\int_{a+\gamma_1}^{b-\gamma_1}\frac{m_3\ee^{2\pi\ii(m_2x+m_3f(x))}f''(x)}{2\pi(m_2+m_3f'(x))^2}f'(x)\,\dd x
=\int_{a+\gamma_1}^{b-\gamma_1}\ee^{\ii F(x)}G(x)\,\dd x,
\end{displaymath}
where
\textcolor{white}{xxxxxxxxxxxxxxxxxxxxxxxxxxxxxx}
\begin{equation}\label{eq6.71}
F(x)=2\pi(m_2x+m_3f(x))
\quad\mbox{and}\quad
G(x)=\frac{m_3f'(x)f''(x)}{2\pi(m_2+m_3f'(x))^2},
\end{equation}
so that
\textcolor{white}{xxxxxxxxxxxxxxxxxxxxxxxxxxxxxx}
\begin{displaymath}
\frac{F'(x)}{G(x)}
=\frac{4\pi^2(m_2+m_3f'(x))^3}{m_3f'(x) f''(x)}=\frac{4\pi^2}{m_3}H(x),
\end{displaymath}
where
\textcolor{white}{xxxxxxxxxxxxxxxxxxxxxxxxxxxxxx}
\begin{displaymath}
H(x)=\frac{(m_2+m_3f'(x))^3}{f'(x)f''(x)}.
\end{displaymath}
Using \eqref{eq6.3}, \eqref{eq6.61} and \eqref{eq6.71}, we deduce the lower bound
\begin{displaymath}
\left\vert\frac{F'(x)}{G(x)}\right\vert
=\frac{4\pi^2\vert m_2+m_3f'(x)\vert^3}{\vert m_3\vert\,\vert f'(x)\vert\,\vert f''(x)\vert}
\ge c_{30}(\max\{\vert m_2\vert,\vert m_3\vert\})^2
\end{displaymath}
where the constant $c_{30}=c_{30}(f;a,b)>0$ depends at most on the function $f$ in the interval $[a,b]$.
Furthermore, the technical condition implies that we can apply Lemma~\ref{lem62} to each of the at most $c_{12}$
subintervals of $[a+\gamma_1,b-\gamma_1]$, and this leads to the bound
\begin{displaymath}
\left\vert\int_{a+\gamma_1}^{b-\gamma_1}\frac{m_3\ee^{2\pi\ii(m_2x+m_3f(x))}f''(x)}{2\pi(m_2+m_3f'(x))^2}f'(x)\,\dd x\right\vert
\le\frac{c_{31}}{(\max\{\vert m_2\vert,\vert m_3\vert\})^2},
\end{displaymath}
where the constant $c_{31}=c_{31}(f;a,b)>0$ depends at most on the function $f$ in the interval $[a,b]$.
It then follows trivially from \eqref{eq6.66} that
\begin{equation}\label{eq6.72}
\vert\widetilde{\III}_1^{(2)}\vert\le\frac{c_{31}}{(\max\{\vert m_2\vert,\vert m_3\vert\})^2}.
\end{equation}

To study the integrals $\widetilde{\III}_1^{(3)}$ and $\widetilde{\III}_1^{(4)}$ given by \eqref{eq6.67} and \eqref{eq6.68}, note that
\begin{align}
\int_a^{(a+b)/2}\frac{\ee^{2\pi\ii(m_2(x)+m_3f(x))}}{2\pi(m_2+m_3f'(x))}f'(x)\,\dd x
&
=\int_a^{(a+b)/2}\ee^{\ii F(x)}G(x)\,\dd x,
\nonumber
\\
\int_{(a+b)/2}^b\frac{\ee^{2\pi\ii(m_2(x)+m_3f(x))}}{2\pi(m_2+m_3f'(x))}f'(x)\,\dd x
&
=\int_{(a+b)/2}^b\ee^{\ii F(x)}G(x)\,\dd x,
\nonumber
\end{align}
where
\textcolor{white}{xxxxxxxxxxxxxxxxxxxxxxxxxxxxxx}
\begin{equation}\label{eq6.73}
F(x)=2\pi(m_2x+m_3f(x))
\quad\mbox{and}\quad
G(x)=\frac{f'(x)}{2\pi(m_2+m_3f'(x))},
\end{equation}
so that
\textcolor{white}{xxxxxxxxxxxxxxxxxxxxxxxxxxxxxx}
\begin{displaymath}
\frac{F'(x)}{G(x)}
=\frac{4\pi^2(m_2+m_3f'(x))^2}{f'(x)}=4\pi^2H(x),
\end{displaymath}
where
\textcolor{white}{xxxxxxxxxxxxxxxxxxxxxxxxxxxxxx}
\begin{displaymath}
H(x)=\frac{(m_2+m_3f'(x))^2}{f'(x)}.
\end{displaymath}
Using \eqref{eq6.3}, \eqref{eq6.61} and \eqref{eq6.73}, we deduce the lower bound
\begin{displaymath}
\left\vert\frac{F'(x)}{G(x)}\right\vert
=\frac{4\pi^2\vert m_2+m_3f'(x)\vert^2}{\vert f'(x)\vert}
\ge c_{32}(\max\{\vert m_2\vert,\vert m_3\vert\})^2,
\end{displaymath}
where the constant $c_{32}=c_{32}(f;a,b)>0$ depends at most on the function $f$ in the interval $[a,b]$.
Furthermore, the derivative $H'(x)$ is equal to
\begin{displaymath}
\frac{2m_3f'(x)f''(x)(m_2+m_3f'(x))-f''(x)(m_2+m_3f'(x))^2}{(f'(x))^2},
\end{displaymath}
where the denominator is non-zero and the numerator is equal to
\begin{displaymath}
f''(x)(m_2+m_3f'(x))(m_3f'(x)-m_2).
\end{displaymath}
Here, in view of \eqref{eq6.3}, the factor $f''(x)$ is non-zero, and the factors $m_2+m_3f'(x)$ and $m_3f'(x)-m_2$ can each have at most one zero in $[a,b]$.
Hence the interval $[a,b]$ is a union of at most $3$ subintervals such that $F'(x)/G(x)$ is monotonic in each subinterval.
We can apply Lemma~\ref{lem62} to each of these subintervals, and this leads to the bounds
\begin{equation}\label{eq6.74}
\vert\widetilde{\III}_1^{(3)}\vert\le\frac{c_{33}}{(\max\{\vert m_2\vert,\vert m_3\vert\})^2}
\quad\mbox{and}\quad
\vert\widetilde{\III}_1^{(4)}\vert\le\frac{c_{33}}{(\max\{\vert m_2\vert,\vert m_3\vert\})^2},
\end{equation}
where the constant $c_{33}=c_{33}(f;a,b)>0$ depends at most on the function $f$ in the interval $[a,b]$.

Combining \eqref{eq6.64}, \eqref{eq6.70}, \eqref{eq6.72} and \eqref{eq6.74}, we obtain the bound
\begin{equation}\label{eq6.75}
\vert\widetilde{\III}_1(\bfm;C(f;u))\vert\le\frac{c_{34}}{(\max\{\vert m_2\vert,\vert m_3\vert\})^2},
\end{equation}
where the constant $c_{34}=c_{34}(f;a,b)>0$ satisfies $c_{34}=4\max\{c_{29},c_{31},c_{33}\}$.
\end{part1}

\begin{part2}
The idea is to combine the bound \eqref{eq6.75} for the term $\widetilde{\III}_1(\bfm;C(f;u))$
with our earlier bounds \eqref{eq6.39} for $\widetilde{\III}_2(\bfm;C(f;u))$ and \eqref{eq5.35} for $\widetilde{\JJJ}(\bfm;\bfv;H)$.
However, we need some care and make some refinement to the sets $\ZZZ_j$, $j=0,1,2,3,\ldots,$ discussed in Step~4
in Section~\ref{sec5} and defined by \eqref{eq5.38}.

For every $j=0,1,2,3,\ldots$ and $k=0,1,\ldots,j$, let
\begin{equation}\label{eq6.76}
\ZZZ_{j,k}=\left\{\begin{array}{ll}
\{(m_2,m_3)\in\ZZZ_j:2^{k-1}<\min\{\vert m_2\vert,\vert m_3\vert\}\le2^k\},&\mbox{if $k\ne0$},\\
\{(m_2,m_3)\in\ZZZ_j:0\le\min\{\vert m_2\vert,\vert m_3\vert\}\le1\},&\mbox{if $k=0$},
\end{array}\right.
\end{equation}
so that
\textcolor{white}{xxxxxxxxxxxxxxxxxxxxxxxxxxxxxx}
\begin{equation}\label{eq6.77}
\vert\ZZZ_{j,k}\vert
\le2^{j+k+3},
\end{equation}
and the set $\ZZZ_j$ can be written as a disjoint union
\begin{displaymath}
\ZZZ_j=\bigcup_{k=0}^j\ZZZ_{j,k}.
\end{displaymath}
\begin{displaymath}
\begin{array}{c}
\includegraphics[scale=0.8]{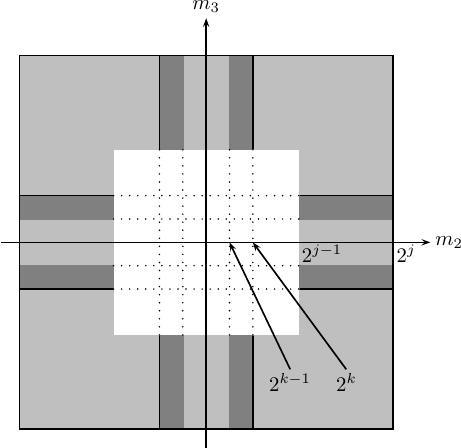}
\\
\mbox{Figure 6.2: the location of $\bfm^*=(m_2,m_3)\in\ZZZ_{j,k}$}
\end{array}
\end{displaymath}

It then follows from \eqref{eq6.46} that
\begin{equation}\label{eq6.78}
\frakI_2=\sum_{j=0}^\infty\sum_{k=0}^j\sum_{\ell=0}^\infty\sum_{\bfm^*\in\ZZZ_{j,k}}\int_{\Omega_j(\bfm^*;\ell)}
\Lambda(\CCC(f;u);H;\bfv;\bfm)\,\dd\alpha_1\,\dd\alpha_2.
\end{equation}
\end{part2}

\begin{part3}
Suppose that $(m_2,m_3)\in\ZZZ_{j,k}$, where $k\ne0$.
It follows from \eqref{eq5.38} and \eqref{eq6.76} that
\begin{equation}\label{eq6.79}
2^{k-1}<\min\{\vert m_2\vert,\vert m_3\vert\}<2^k
\quad\mbox{and}\quad
2^{j-1}<\max\{\vert m_2\vert,\vert m_3\vert\}<2^j.
\end{equation}
Suppose furthermore that $(\alpha_1,\alpha_2)\in\Omega_j(\bfm^*;0)$.
Then it follows from \eqref{eq5.35}, \eqref{eq6.36}, \eqref{eq6.39}, \eqref{eq6.75} and \eqref{eq6.79} that
\begin{align}
\Lambda(\CCC(f;u);H;\bfv;\bfm)
&
\le\frac{c_{34}H}{(\max\{\vert m_2\vert,\vert m_3\vert\})^2}\min\left\{u,\frac{1}{u(\pi m_3)^2}\right\}
\nonumber
\\
&
\le\frac{c_{34}H}{(\max\{\vert m_2\vert,\vert m_3\vert\})^2}\min\left\{u,\frac{1}{u\pi^2(\min\{\vert m_2\vert,\vert m_3\vert\})^2}\right\}
\nonumber
\\
&
\le\frac{c_{34}H}{2^{2j-2}}\min\left\{u,\frac{1}{u\pi^22^{2k-2}}\right\}.
\nonumber
\end{align}
Combining this with \eqref{eq5.41} and \eqref{eq6.77}, we deduce that
\begin{align}\label{eq6.80}
&
\sum_{\bfm^*\in\ZZZ_{j,k}}\int_{\Omega_j(\bfm^*;0)}
\Lambda(\CCC(f;u);H;\bfv;\bfm)\,\dd\alpha_1\,\dd\alpha_2
\nonumber
\\
&\quad
\le\frac{2^{j+k+9}}{H}\frac{c_{34}H}{2^{2j-2}}\min\left\{u,\frac{1}{u\pi^22^{2k-2}}\right\}
=\frac{2^{11}c_{34}}{2^j}\min\left\{2^ku,\frac{4}{u\pi^22^k}\right\}.
\end{align}
It is easily checked that the inequality \eqref{eq6.80} holds for $k=0$ also.

Let $K^\ddagger$ denote the largest integer $k$ such that
\begin{displaymath}
2^ku\le\frac{4}{u\pi^22^k},
\quad\mbox{equivalent to}\quad
2^k\le\frac{2}{u\pi}.
\end{displaymath}
Then
\textcolor{white}{xxxxxxxxxxxxxxxxxxxxxxxxxxxxxx}
\begin{equation}\label{eq6.81}
2^{K^\ddagger}\le\frac{2}{u\pi}<\frac{1}{u}
\quad\mbox{and}\quad
2^{K^\ddagger+1}>\frac{2}{u\pi},
\end{equation}
and
\textcolor{white}{xxxxxxxxxxxxxxxxxxxxxxxxxxxxxx}
\begin{align}\label{eq6.82}
&
\sum_{k=0}^j\min\left\{2^ku,\frac{4}{u\pi^22^k}\right\}
\le\sum_{k=0}^{K^\ddagger}2^ku+\frac{4}{u\pi^2}\left(\frac{1}{2^{K^\ddagger+1}}+\ldots+\frac{1}{2^j}\right)
\nonumber
\\
&\quad
\le2^{K^\ddagger+1}u+\frac{8}{u\pi^22^{K^\ddagger+1}}
<4,
\end{align}
in view of \eqref{eq6.81}.
It now follows from \eqref{eq6.80} and \eqref{eq6.82} that
\begin{equation}\label{eq6.83}
\sum_{j=0}^\infty\sum_{k=0}^j\sum_{\bfm^*\in\ZZZ_{j,k}}\int_{\Omega_j(\bfm^*;0)}
\Lambda(\CCC(f;u);H;\bfv;\bfm)\,\dd\alpha_1\,\dd\alpha_2
\le\sum_{j=0}^\infty\frac{2^{13}c_{34}}{2^j}
\le c_{35},
\end{equation}
where the constant $c_{35}=c_{35}(f;a,b)>0$ depends at most on the function $f$ in the interval $[a,b]$.

Suppose next that $(m_2,m_3)\in\ZZZ_{j,k}$, where $k\ne0$, and $(\alpha_1,\alpha_2)\in\Omega_j(\bfm^*;\ell)$,
where $\ell\ge1$.
Then it follows from \eqref{eq5.35}, \eqref{eq5.40}, \eqref{eq6.36}, \eqref{eq6.39}, \eqref{eq6.75} and \eqref{eq6.79} that
\begin{align}
\Lambda(\CCC(f;u);H;\bfv;\bfm)
&
\le\frac{2c_{34}}{(\max\{\vert m_2\vert,\vert m_3\vert\})^2}\min\left\{u,\frac{1}{u(\pi m_3)^2}\right\}\frac{2}{\pi^2H(\bfv\cdot\bfm)^2}
\nonumber
\\
&
\le\frac{2c_{34}H}{2^{2\ell}(\max\{\vert m_2\vert,\vert m_3\vert\})^2}\min\left\{u,\frac{1}{u\pi^2(\min\{\vert m_2\vert,\vert m_3\vert\})^2}\right\}
\nonumber
\\
&
\le\frac{2c_{34}H}{2^{2\ell}2^{2j-2}}\min\left\{u,\frac{1}{u\pi^22^{2k-2}}\right\}.
\nonumber
\end{align}
Combining this with \eqref{eq5.41} and \eqref{eq6.77}, we deduce that
\begin{align}\label{eq6.84}
&
\sum_{\bfm^*\in\ZZZ_{j,k}}\int_{\Omega_j(\bfm^*;\ell)}
\Lambda(\CCC(f;u);H;\bfv;\bfm)\,\dd\alpha_1\,\dd\alpha_2
\nonumber
\\
&\quad
\le\frac{2^{j+k+10}2^\ell}{H}\frac{c_{34}H}{2^{2\ell}2^{2j-2}}\min\left\{u,\frac{1}{u\pi^22^{2k-2}}\right\}
=\frac{2^{12}c_{34}}{2^\ell2^j}\min\left\{2^ku,\frac{4}{u\pi^22^k}\right\}.
\end{align}
It now follows from \eqref{eq6.82} and \eqref{eq6.84} that
\begin{align}\label{eq6.85}
&
\sum_{j=0}^\infty\sum_{k=0}^j\sum_{\ell=1}^\infty\sum_{\bfm^*\in\ZZZ_{j,k}}\int_{\Omega_j(\bfm^*;\ell)}
\Lambda(\CCC(f;u);H;\bfv;\bfm)\,\dd\alpha_1\,\dd\alpha_2
\nonumber
\\
&\quad
\le\sum_{j=0}^\infty\sum_{k=0}^j\sum_{\ell=1}^\infty\frac{2^{12}c_{34}}{2^\ell2^j}\min\left\{2^ku,\frac{4}{u\pi^22^k}\right\}
\le\sum_{j=0}^\infty\frac{2^{14}c_{34}}{2^j}
\le2c_{35}.
\end{align}
\end{part3}

The desired inequality \eqref{eq6.60} follows on combining \eqref{eq6.78}, \eqref{eq6.83} and \eqref{eq6.85}.
\end{proof}

It now follows from \eqref{eq6.44}, \eqref{eq6.47} and \eqref{eq6.60} that
\begin{displaymath}
\int_0^1\int_0^1\Xi(\alpha_1,\alpha_2)\,\dd\alpha_1\,\dd\alpha_2\le\psi(u),
\end{displaymath}
where $\psi(u)$ is given by \eqref{eq6.6}.
Then for every real parameter $\kappa>1$, we have
\begin{displaymath}
\lambda_2(\{(\alpha_1,\alpha_2)\in[0,1]^2:\Xi(\alpha_1,\alpha_2)\ge\kappa\psi(u)\})\le\frac{1}{\kappa}.
\end{displaymath}
\end{step6}


\begin{step7}
We can now justify all the steps in a similar way as in Step~7 in Section~\ref{sec5}.
\end{step7}

This completes the proof of Lemma~\ref{lem61}.
\end{proof}

%
%

\end{document}